\documentclass[twoside]{article}
\usepackage{psfig,epsf}
\usepackage{amsmath,amssymb,amsthm,hyperref}
\newcommand{\A}{ \mathbb{A}}
\newcommand{\bb}{ \overline{b} }
\newcommand{\CC}{ \mathcal{C} }
\newcommand{\D}{ \mathcal{D} }
\newcommand{\dd}{\mathsf{d}}
\newcommand{\eps}{\varepsilon }
\newcommand{\E}{\mathbb E }
\newcommand{\EE}{\mathcal{E}}

\newcommand{\F}{\mathcal F }

\newcommand{\G}{\mathbf{\Gamma}}
\newcommand{\GG}{\mathbb G }
\newcommand{\HH}{\mathbb{H}}
\newcommand{\I}{\mathbf I }
\newcommand{\overI}{\overline{\I}}
\newcommand{\partI}{\partial \I}
\newcommand{\II} {1\!\!1}
\newcommand{\LL}{\mathcal{L}}
\newcommand{\M}{\mathbb{M}}
\newcommand{\sfM}{\mathsf M}
\newcommand{\Om}{\Omega }
\newcommand{\PP}{\mathbb P }
\newcommand{\R}{ \mathbb R }
\newcommand{\s}{\sigma }
\newcommand{\sg}{\overline{\sigma} }

\newcommand{\T}{ \mathbf T }
\newcommand{\sfS}{\mathsf S}
\newcommand{\overU}{\overline{\U}}
\newcommand{\U}{ \mathbf U }
\newcommand{\partU}{ \partial \U }
\newcommand{\overV}{\overline{\V}}
\newcommand{\V}{ \mathbf V }
\newcommand{\partV}{ \partial \V }
\newcommand{\W}{ \mathbf W }
\newcommand{\oomega}{\overline{\omega}}
\newcommand{\x}{\mathbf x}
\newcommand{\X}{\mathbf X }
\newcommand{\y}{\mathbf y}
\newcommand{\Y}{ \mathbf Y }
\newcommand{\z}{\mathbf z}
\newcommand{\nul} {\mathbf{0}}

\newcommand{\rank}{\operatorname{rank}}
\newcommand{\supp}{\operatorname{supp}}

\catcode`\@=11
\long\def\@makefntext#1{
\protect\noindent \hbox to 3.2pt {\hskip-.9pt
$^{{\eightrm\@thefnmark}}$\hfil}#1\hfill}       

\def\@makefnmark{\hbox to 0pt{$^{\@thefnmark}$\hss}}    

\def\ps@myheadings{\let\@mkboth\@gobbletwo
\def\@oddhead{\hfill\hbox{}\rightmark}
\def\@oddfoot{}\def\@evenhead{\leftmark\hbox{}\hfill}\def\@evenfoot{}
\def\sectionmark##1{}\def\subsectionmark##1{}}

\oddsidemargin=\evensidemargin
\newcounter{sectionc}\newcounter{subsectionc}\newcounter{subsubsectionc}
\renewcommand{\section}[1] {\vspace{12pt}\addtocounter{sectionc}{1}
\setcounter{subsectionc}{0}\setcounter{subsubsectionc}{0}\noindent
    {\tenbf\thesectionc. #1}\par\vspace{5pt}}
\renewcommand{\subsection}[1] {\vspace{12pt}\addtocounter{subsectionc}{1}
    \setcounter{subsubsectionc}{0}\noindent
  {\bf\thesectionc.\thesubsectionc. {\kern1pt \bfit #1}}\par\vspace{5pt}}
\renewcommand{\subsubsection}[1] {\vspace{12pt}\addtocounter{subsubsectionc}{1}
    \noindent{\tenrm\thesectionc.\thesubsectionc.\thesubsubsectionc.
    {\kern1pt \tenit #1}}\par\vspace{5pt}}
\newcommand{\nonumsection}[1] {\vspace{12pt}\noindent{\tenbf #1}
    \par\vspace{5pt}}

\newcounter{appendixc}
\newcounter{subappendixc}[appendixc]
\newcounter{subsubappendixc}[subappendixc]
\renewcommand{\thesubappendixc}{\Alph{appendixc}.\arabic{subappendixc}}
\renewcommand{\thesubsubappendixc}
    {\Alph{appendixc}.\arabic{subappendixc}.\arabic{subsubappendixc}}

\renewcommand{\appendix}[1] {\vspace{12pt}
        \refstepcounter{appendixc}
        \setcounter{figure}{0}
        \setcounter{table}{0}
        \setcounter{lemma}{0}
        \setcounter{theorem}{0}
        \setcounter{corollary}{0}
        \setcounter{definition}{0}
        \setcounter{equation}{0}
        \renewcommand{\thefigure}{\Alph{appendixc}.\arabic{figure}}
        \renewcommand{\thetable}{\Alph{appendixc}.\arabic{table}}
        \renewcommand{\theappendixc}{\Alph{appendixc}}
        \renewcommand{\thelemma}{\Alph{appendixc}.\arabic{lemma}}
        \renewcommand{\thetheorem}{\Alph{appendixc}.\arabic{theorem}}
        \renewcommand{\thedefinition}{\Alph{appendixc}.\arabic{definition}}
        \renewcommand{\thecorollary}{\Alph{appendixc}.\arabic{corollary}}
        \renewcommand{\theequation}{\Alph{appendixc}.\arabic{equation}}
        \noindent{\tenbf Appendix#1}\par\vspace{5pt}}
\newcommand{\subappendix}[1] {\vspace{12pt}
        \refstepcounter{subappendixc}
        \noindent{\bf Appendix \thesubappendixc. {\kern1pt \bfit #1}}
    \par\vspace{5pt}}
\newcommand{\subsubappendix}[1] {\vspace{12pt}
        \refstepcounter{subsubappendixc}
        \noindent{\rm Appendix \thesubsubappendixc. {\kern1pt \tenit #1}}
    \par\vspace{5pt}}

\newcommand{\textlineskip}{\baselineskip=13pt}
\newcommand{\smalllineskip}{\baselineskip=10pt}

\newcommand{\copyrightheading}[1]
    {\vspace*{-2.5cm}\smalllineskip{\flushleft
    {\footnotesize Stochastics and Dynamics #1}\\
    {\footnotesize \copyright\kern2pt World Scientific Publishing
     Company}\\
     }}

\def\abstracts#1#2#3{{
    \centering{\begin{minipage}{4.5in}\footnotesize\baselineskip=10pt
    \parindent=0pt #1\par
    \parindent=15pt #2\par
    \parindent=15pt #3
    \end{minipage}}\par}}

\def\keywords#1{{
    \centering{\begin{minipage}{4.5in}\footnotesize\baselineskip=10pt
    {\footnotesize\it Keywords}\/: #1
     \end{minipage}}\par}}

\newcounter{itemlistc}
\newcounter{romanlistc}
\newcounter{alphlistc}
\newcounter{arabiclistc}

\newenvironment{romanlist}
    {\setcounter{romanlistc}{0}
     \begin{list}{$($\roman{romanlistc}$)$}
    {\usecounter{romanlistc}
     \setlength{\parsep}{0pt}
     \setlength{\itemsep}{0pt}}}{\end{list}}

\newcommand{\fcaption}[1]{
        \refstepcounter{figure}
        \setbox\@tempboxa = \hbox{\footnotesize Fig.~\thefigure. #1}
        \ifdim \wd\@tempboxa > 5in
           {\begin{center}
        \parbox{5in}{\footnotesize\smalllineskip Fig.~\thefigure. #1}
            \end{center}}
        \else
             {\begin{center}
             {\footnotesize Fig.~\thefigure. #1}
              \end{center}}
        \fi}

\newcommand{\tcaption}[1]{
        \refstepcounter{table}
        \setbox\@tempboxa = \hbox{\footnotesize Table~\thetable. #1}
        \ifdim \wd\@tempboxa > 5in
           {\begin{center}
        \parbox{5in}{\footnotesize\smalllineskip Table~\thetable. #1}
            \end{center}}
        \else
             {\begin{center}
             {\footnotesize Table~\thetable. #1}
              \end{center}}
        \fi}

\def\pmb#1{\setbox0=\hbox{#1}
    \kern-.025em\copy0\kern-\wd0
    \kern.05em\copy0\kern-\wd0
    \kern-.025em\raise.0433em\box0}

\def\fnt#1#2{\footnotetext{\kern-.3em
        {$^{\mbox{\scriptsize #1}}$}{#2}}}

\font\tenrm=cmr10
\font\tenit=cmti10
\font\tenbf=cmbx10
\font\bfit=cmbxti10 at 10pt
\font\ninerm=cmr9

\font\eightrm=cmr8

\newtheorem{theorem}{Theorem}[sectionc]

\newtheorem{lemma}{Lemma}[sectionc]

\newtheorem{proposition}{Proposition}[sectionc]

\theoremstyle{definition}
\newtheorem{definition}{Definition}[sectionc]

\theoremstyle{definition}
\newtheorem{example}{Example}[sectionc]

\theoremstyle{definition}
\newtheorem{remark}{Remark}[sectionc]

\def\qed{\hbox{${\vcenter{\vbox{         
   \hrule height 0.4pt\hbox{\vrule width 0.4pt height 6pt
   \kern5pt\vrule width 0.4pt}\hrule height 0.4pt}}}$}}

\def\theequation{\thesectionc.\arabic{equation}}  

        {\setcounter{itemlistc}{0}      
     \begin{list}{$\bullet$}        
    {\usecounter{itemlistc}         
     \leftmargin10pt           
     \setlength{\parsep}{0pt}
     \setlength{\itemsep}{0pt}     
    }}{\end{list}}

    {\setcounter{romanlistc}{0}     
     \begin{list}{$($\roman{romanlistc}$)$} 
    {\usecounter{romanlistc}        
     \leftmargin18pt 
     \setlength{\parsep}{0pt}
     \setlength{\itemsep}{0pt}  
     \settowidth{\labelwidth}{#1}
    }}{\end{list}}

    {\setcounter{enumii}{0}         
     \begin{list}{$($\alph{enumii}$)$}  
    {\usecounter{enumii}            
     \leftmargin18pt        
     \setlength{\parsep}{0pt}
     \setlength{\itemsep}{0pt}  
     \settowidth{\labelwidth}{#1}
    }}{\end{list}}

\def\bsc{{\sc a\kern-6.4pt\sc a\kern-6.4pt\sc a}}   
\def\bflatex{\bf L\kern-.30em\raise.3ex\hbox{\bsc}\kern-.14em
T\kern-.1667em\lower.7ex\hbox{E}\kern-.125em X}

\begin{document}
\setlength{\textheight}{7.7truein}    

\markboth{\protect{\footnotesize\it N. V. O'Bryant}}{\protect{\footnotesize\it Noisy system}}

\normalsize\textlineskip
\thispagestyle{empty}
\setcounter{page}{1}

\copyrightheading{}

\vspace*{1in}         

\centerline{\bf A NOISY SYSTEM WITH A FLATTENED HAMILTONIAN} \baselineskip=13pt \centerline{\bf AND MULTIPLE
TIME SCALES }

\vspace*{0.37truein} \centerline{\footnotesize NATELLA V. O'BRYANT} \baselineskip=12pt
\centerline{\footnotesize\it Department of Mathematics, University of California - Irvine} \baselineskip=10pt
\centerline{\footnotesize\it 272 Multipurpose Science \& Technology Building, Irvine, CA 92697}

\vspace*{0.225truein}

\vspace*{0.21truein} \abstracts{We consider a two-dimensional weakly dissipative dynamical system with
time-periodic drift and diffusion coefficients. The average of the drift is governed by a degenerate
Hamiltonian whose set of critical points has an interior. The dynamics of the system is studied in the
presence of three time scales. Using the martingale problem approach and separating the time scales, we
average the system to show convergence to a Markov process on a stratified space. The averaging combines the
deterministic time averaging of periodic coefficients, and the stochastic averaging of the resulting system.
The corresponding strata of the reduced space are a two-sphere, a point and a line segment. Special attention
is given to the description of the domain of the limiting generator, including the analysis of the gluing
conditions at the point where the strata meet. These gluing conditions, resulting from the effects of the
hierarchy of time scales, are similar to the conditions on the domain of skew Brownian motion and are related
to the description of spider martingales.}{}{}

\vspace*{5pt} \keywords{Hamiltonian systems, Markov processes, stochastic averaging, martingale problem.}

\vspace*{4pt}
\baselineskip=13pt              
\normalsize                 
\section{Introduction}    
\noindent Scientists have been studying Hamiltonian dynamical systems for centuries. One of the flourishing
directions emerging as a result of these studies is current research on problems involving oscillations and
vibrations of real mechanical systems. Generally, models that reflect such processes of real life are
nonlinear and stochastic, and they often require some model reduction. One specific technique of model
reduction is stochastic averaging which reduces more complicated random processes to simpler ones. We are
interested in studying both the means of approximation of the initial model by its reduced model and the
nature of the reduced model. One of the features of the reduced model is the structure of its state space,
which may acquire dimensional discontinuities becoming a stratified space \cite{G.-M.}. In cases like these,
the domains of the generator of the limiting reduced process may require gluing conditions. The dynamics of
the reduced process can be understood via the Fokker-Planck equation, with the gluing conditions reflected in
the corresponding conservation of flux and continuity equations.

Real-life dynamical systems are typically described with equations which have small parameters. As these
parameters vary, they can cause radical changes in the qualitative structure of solutions. Although invisible
on a usual time scale, in many cases these changes may become dramatically visible on a larger time scale.
Bogolyubov \cite{B.} proved a deterministic result for ordinary differential equations which resolves a
difficulty of this type, and  which is known as the averaging principle.

The fact that most physical problems are non-deterministic inspired the development of stochastic averaging.
The first rigorous arguments for stochastic averaging were given by Has'minski{\u\i} \cite{H.2}, \cite{H.3}.
The essence of the classical stochastic averaging is the asymptotic separation of time scales.

Consider the following diffusion on $\R^2:$
    \begin{align}
    \dot{\X}^\eps_t(\x) & = \nabla^{\perp} H \left(\X^\eps_t(\x)\right)+ \eps \dot{\W}_t \notag\\
    \X^\eps_0(\x) & =  \x,\notag
    \end{align}
where $\nabla^\perp$ denotes a vector orthogonal to $\nabla,$ $H$ is a sufficiently smooth function, $\eps$~
is a small parameter, and $\W$ is two-dimensional Brownian motion. Here $H$ is the first integral of the
corresponding deterministic system
\begin{equation*}
    \begin{split}
    \dot{\X}^0_t(\x) & = \nabla^{\perp} H \left(\X^0_t(\x)\right)\\
    \X^0_0(\x) & =  \x,
    \end{split}
    \end{equation*}
This system is Hamiltonian with one degree of freedom. If the energy function $H$ is a single-well
Hamiltonian, i.e., if all the level sets $L_h = \left\{\x \in \R^2: H(\x)=h \right\}$ of $H$ have a single
connected component formed by periodic trajectories of the system, then the invariant measure of this system
depends only on the height $h$. There are two time scales involved here. The movement of the deterministic
system is fast, but the energy of the stochastic system is varying slowly, and this slow transversal movement
is visible only on a larger time scale. It is natural to separate the two time scales in order to witness the
energy dynamics. From \cite{H.3}, on a large time scale, $H(\X^\eps)$ weakly converges to a Markov process on
a reduced space, which is the real line. In this case, the stratified space has one one-dimensional stratum.

Freidlin and Wentzell \cite{F.-W.}, and then Freidlin and Weber \cite{F.-Wb.}, considered the problem with a
multiple-well Hamiltonian $H$. One of the essential assumptions made on the system in their work was the
non-degenerate structure of the Hamiltonian. In the case of a multiple-well energy function, the invariant
measure of the system depends not only on the height, but also on the particular component of the
Hamiltonian. In the presence of saddle points the reduced space becomes a tree. Each vertex of the tree
corresponds either to a single critical point of $H$, or to a whole homoclinic trajectory of the system.
Every interior point of a segment on the tree corresponds to a particular periodic orbit. The energy of the
stochastic system is varying slowly, and this slow transversal movement is naturally visible on a large time
scale. Studying the energy dynamics as $\eps \rightarrow 0,$ the authors identified the limiting graph-valued
process. They showed that the classical calculations of Has'minski{\u\i} \cite{H.3} can be applied, as long
as the process stays away from the vertices of the graph. As for the vertices, the limiting Markov process
requires gluing conditions to be enforced on the domain of its generator. The dynamics of the limiting
process can be understood via the Fokker-Plank equation. The gluing reflects the two following properties.
The solution of the Fokker-Plank equation for the corresponding density must be continuous at vertices, and
the weighted fluxes trough vertices from different legs of a tree must sum to zero.

One way to describe a diffusion is to illustrate the dynamics of its trajectories by a corresponding
stochastic differential equation. Another way, introduced by Stroock and Varadhan \cite{S.-V.}, is to build a
martingale on the corresponding coordinate space. This alternative method is used to construct a limiting
process by constructing a martingale problem for which this process is a solution. This is known as the
martingale problem approach.

Sowers \cite{S.1} investigated a stochastic system with a flattened Hamiltonian, removing the assumption of
non-degeneracy. A canonical example of a Hamiltonian of this kind is
\begin{equation*}
H(\x) = \biggl(\max \left\{ 0, \|\x\|-1\right\}\biggr)^n, \ n > 2.
\end{equation*}
To regularize noise, the chain equivalence \cite{R.} was used. This general equivalence relation identifies
points if they can be taken from one to the other and back by a combination of small diffusive perturbation
and fast drift. In this problem the limiting space has a dimensional discontinuity. The author showed
convergence to a Markov process on the reduced space, whose strata are a two-sphere, a line segment and a
point. A step towards a general theory of Markov processes on stratified spaces was made in the recent work
of Evans and Sowers \cite{E.-S.}.

Recently, Namachchivaya and Sowers \cite{N.-S.} made an attempt to develop a unified approach to study the
energy dynamics of a non-degenerate single-degree-of-freedom system excited by both periodic and random
perturbations. As a prototype for one of the cases, they studied a two-dimensional weakly dissipative system
with time-periodic coefficients, and achieved model reduction through stochastic averaging.

We combine \cite{S.1} and \cite{N.-S.}. We consider a two-dimensional weakly dissipative system with
time-periodic coefficients whose time average is governed by a degenerate Hamiltonian whose set of critical
points has an interior. The dynamics of the system is studied in the presence of three time scales. The
periodic fluctuations of the coefficients occur on the time scale of order $1/\eps^{2}$, the effect of the
drift is visible on the time scale of order $1/\eps$, and the diffusion coefficients are of order one.

The main technical difficulty in this problem is stochastic averaging near the boundary of the critical set.
The main step of the averaging is to approximately solve a certain partial differential equation. This
equation involves drift of order one and small diffusion leading to a singular perturbation problem. We solve
this problem by introducing new coordinates in the vicinity of the boundary. This coordinate transformation
was inspired by Has'minski{\u\i} \cite{H.1}, and was previously used in \cite{S.1}. These coordinates were
introduced with the idea that equal angular displacement should correspond to equal amounts of diffusion
across the boundary of the critical set. They perform a role similar to action-angle components. Solving the
partial differential equation in Has'minski{\u\i}'s coordinates involves Fourier analysis and leads to an
inhomogeneous Bessel's ordinary differential equation. The unique solution is expressed via Bessel functions
of small order and purely imaginary argument.

We perform stochastic averaging using the martingale problem. This approach was developed by Papanicolaou and Kohler
\cite{P.-K.}. Using the martingale problem and separating the time scales, we average the system to show convergence to
a Markov process on a stratified space. The corresponding strata are, as in \cite{S.1}, a two-sphere, a point and a
line segment. Special attention is given to the description of the domain of the limiting generator, including the
analysis of the gluing condition at the point where the strata meet. The differential operator governing the limiting
process and the gluing describing the domain of its generator are given explicitly (Section 4.1). The gluing condition,
resulting from the effects of the hierarchy of time scales, is similar to the conditions described in \cite{S.1}. It is
also similar to the conditions on the domain of skew Brownian motion, and is related to the description of spider
martingales \cite{Y.1}.

In our problem, the particular interest is in a careful analysis of the gluing condition on the domain of the limiting
generator. Here, the multiple time scales affect gluing in a straightforward way. It would be interesting to see the
changes to the initial problem which result in various changes to the structure of the reduced stratified space and
which affect gluing in a different way.

\subsection{Notation}
\noindent Let $\alpha = (\alpha_1,\alpha_2)$ be an ordered tuple of nonnegative integers, i.e., let $\alpha $
be a multi-index. We set $|\alpha| = \alpha_1+\alpha_2$ and
    \begin{equation*}
    \partial^\alpha =
    \left(\frac{\partial}{\partial x_1}\right)^{\alpha_1}
    \left( \frac{\partial}{\partial x_2} \right)^{\alpha_2}.
    \end{equation*}

We denote by $C^k(A;B)$ the collection of all functions $f: A \rightarrow B$ which have continuous derivatives of
orders $0,1, \ldots, k$. When $B$ is clear from the context, we write simply $C^k(A)$. We say that a function $f$ is a
$C^k$-diffeomorphism if $f$ is a $C^k$-homeomorphism from $A$ to $B$ with a $C^k$ inverse. We say that $f(\x,t): A
\times T \rightarrow B$ is $C^{k,l}(A \times T; B)$ if $\partial^\alpha f$ is continuous for all $|\alpha| \leq k$ and
$\left( \frac{\partial f}{\partial t}\right)^\beta$ is continuous for all $0 \leq \beta\leq l.$ If $\ l=0$ we write
$C^k(A \times T;B),$ and if $l=k=0$ we write simply $C(A \times T;B)$. We denote by $C^k_c(A)$ the subset of $C^k(A)$
of functions having compact support, and by $C_b(A)$ the subset of continuous functions $C(A)$ bounded on the set $A$.

    Let $U$ be an open subset of $\R^2$. We denote by $\HH^k(U)$ the Sobolev space, which consist of all locally
integrable functions $h: U \rightarrow \R$ such that for each multi-index $\alpha$ with $ |\alpha| \leq k$,
$\partial^\alpha f$ exists in a weak sense and belongs to the space of square-integrable functions $L^2(U).$
For every $h \in \HH^k(U)$, we define
    \begin{equation*}
    \|h\|_{\HH^k(U)}:= \left( \sum_{|\alpha|\leq k } \int_U \left| \partial^\alpha h \right|^2 \right)^{1/2}.
    \end{equation*}
For any functions $f \in C^{k,l}\left(\R^2 \times \R; \R \right)$ and $g \in C^k\left(\R; \R \right) \
\left(k,l \geq 0 \right) $ and for any compact subsets $ A \subset \R^2,$ and $B, T \subset \R,$ we define
    \begin{equation*}
    \left\| f \right\|_{C^k \left(A \times T \right)}
    := \sup_{\substack{\x \in A \\t \in T \\|\alpha| \leq k}}
    \left| \partial^\alpha f(\x,t)\right|
    \end{equation*} and
    \begin{equation*}
    \left\| g \right\|_{C^k \left(B \right)}
    := \sup_{\substack{\y \in B \\|\alpha| \leq k}}
    \left| \partial^\alpha g(\y)\right|.
    \end{equation*}
We say that a function $f \in C^2(\R^2)$ is non-degenerate at a point $\x$ if the Hessian
    \begin{equation*}
    D^2 f(\x) =
    {\left|
    \begin{array}{cc}
      \frac{\partial^2 f}{\partial x_1 \partial x_1}(\x)
        & \frac{\partial^2 f}{\partial x_1 \partial x_2} (\x)\\
      \frac{\partial^2 f}{\partial x_2 \partial x_1}(\x)
        & \frac{\partial^2 f}{\partial x_2 \partial x_2}(\x)\
    \end{array}
    \right|}
    \end{equation*}
of $f$ at $\x$ is non-degenerate, that is, if $\rank D^2 f(\x) = 2.$ We call a function $f~\in~C^2(\R^2)$
degenerate at
 $\x$ if $\rank D^2 f(\x) < 2.$

\section{Statement of the Problem} \noindent
Since we will deal with a perturbed dynamical system, it makes sense to begin by describing our problem from
the point of view of the theory of dynamical systems. To show that a certain limiting system is Markovian, we
use the martingale problem, and we will need to restate our problem in the language more appropriate for this
approach.

\subsection{Mechanical description of the problem} \noindent
When we speak about small perturbations, we will mean that the perturbing effect depends on a small
parameter~$\eps$. Consider small random perturbations of a two-dimensional weekly-dissipative system with
time-periodic coefficients. In other words, our interest is a system which has both periodic and stochastic
perturbations,
    \begin{align}
   & d \widetilde{\X}_t^\eps  = \eps b(\widetilde{\X}_t^\eps, t) \ d t
    + \eps \s (\widetilde{\X}_t^\eps, t) \ d \W_t, \ t
    \geq 0 \notag\\
    & \widetilde{\X}^\eps_0  = \x  \text{ a.s. } \label{eqn:1}
    \end{align}
Here $\W$ is a two-dimensional Brownian motion. Coefficients $b$ and $\s$ in (\ref{eqn:1}) are $C^{2,1}
\left( \R^2 \times \R \right)$-functions and are $2 \pi$-periodic in time. We also assume that, given any
$\x=(x_1,x_2) \in \R^2,$
    \begin{equation} \label{eqn:perpH}
    (\M b) (\x) := \frac{1}{2 \pi} \int_0^{2\pi} b(\x,t) dt
        = \nabla^{\perp} H(\x),
    \end{equation}
where $\nabla^{\perp}H(\x) = \left(\frac{\partial H(\x)}{\partial x_2},-\frac{\partial H(\x)}{\partial
x_1}\right)$ is perpendicular to $\nabla H(\x)$. We assume that $H  \geq 0$ is smooth and that it has a
single well, that is, its level sets have only one connected component. Typically one assumes that $H$ is
non-degenerate at its critical points. Here we allow $H$ to be degenerate at its critical points, so that
$\rank D^2 H(\x) \leq 2,$ and that the set of the critical points of $H$ can have an interior.

We will refer to the quantity $H$ as the energy of the system (\ref{eqn:1}). To describe $H$ with some rigor,
we introduce a function $K \in C^\infty \left(\R^2;\R\right)$ with the following properties:
\begin{romanlist}
    \item The set $ \left\{\x \in \R^2: K(\x) \leq 0 \right\}$ is a diffeomorphism of the unit disc.
    \item $K$ is nonsingular on the set $\left\{\x \in \R^2: K(\x) \geq 0 \right\},$
    meaning that $\nabla K(\x) \neq \nul$ on this set.
    \item $\lim_{\|\x\| \rightarrow \infty}K(\x) = \infty$.
\end{romanlist}
Having in mind these properties of $K$, we set
    \begin{equation*}
    H(\x) := \biggl( \max \left\{ 0, K(\x) \right\}
        \biggr)^n, n> 2.
    \end{equation*}
 We denote
 \begin{equation*}
 \V := \left\{\x \in \R^2: K(\x) < 0 \right\}.
 \end{equation*}

\begin{example} \label{ex:H}
Let $K$ satisfy the above properties, $K(\x) = \|\x\|-1$ for $\|\x\| \geq 1$ and $K(\x) < 0 $ for $\|\x\| <
1$. Then $H(\x) = \biggl(\max \left\{ 0,\| \x\|-1 \right\} \biggr)^n, \ n>2,$ does not have any isolated
critical points ~(Fig. \ref{fig:H}).
\begin{figure}[h]
\centerline{\psfig{file=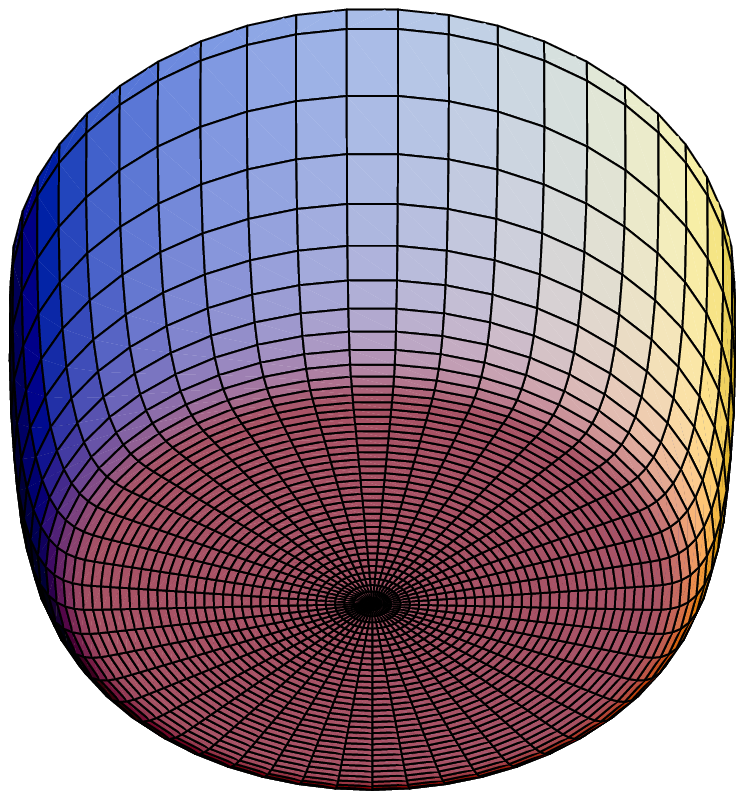,width=2.5in}} \vspace*{6pt} \fcaption{$H$: View from below.}
\label{fig:H} 
\end{figure}
 Instead, the set $\overV$ of its critical points
has an interior $ \V = \left\{\x \in \R^2: K(\x) < 0 \right\} = \left\{\x \in \R^2: \|\x\| < 1 \right\}. $
For any $\x \in \overV,$ $ \rank D^2 H(\x) = 0,$ so all the critical points of $H$ are degenerate. Any
positive level set $ L_h := \left\{\x \in \R^2: H(\x) = h \right\}, \ h>0, $ of $H$ is a connected closed
smooth curve in $\R^2.$
\end{example}

\begin{remark} \label{rem:K} The property (ii) is essential for our construction. Set
    \begin{equation*}
    \overline{K}(\x) =
    \begin{cases}
        0 & \text{ if  }\,\|\x\| \leq 1,\\
       \|\x\|\, \exp\left\{\frac{1}{1-\|\x\|^2}\right\} & \text{ if  }\, \|\x\| > 1.
    \end{cases}
    \end{equation*}
This smooth function satisfies the properties (i) and (iii), but not (ii), since it is singular on the set
$\left\{\x \in \R^2: \overline{K}(\x)=0 \right\}.$

The fact that $K$ is nonsingular on $\left\{\x \in \R^2: K(\x) \geq 0 \right\}$ implies that there exists $\
a>0$ such that the small set
    \begin{equation} \label{eqn:EE}
    \EE = \left\{\y \in \R^2: \left| K(\y) \right| < a \right\}
    \end{equation}
surrounding the boundary $\ \partV$ of the set $\ \overV = \left\{\x \in \R^2: K(\x) \leq 0 \right\}$ has one
component and $\nabla K(\x) \neq \nul$ for all $\x \in \overline{\EE}.$
\end{remark}

Our goal is to study the dynamics of the energy of the system (\ref{eqn:1}) as the small parameter $\eps$
tends to $0$ and to describe the limiting system. We will show that it has the Markovian property.

\begin{definition}
We write $\gamma = \s \s^T,$ and for any functions $f,g \in C^{2}(\R^2 ;\R)$  and for every $t \in \R$, we
define two differential operators: the generator
    \begin{equation} \label{eqn:L}
    (\LL_t f)(\x)
    := \frac{1}{2} \sum_{i,j=1}^2 \gamma_{ij}(\x,t)
        \frac{\partial^2 f}{\partial x_i \partial x_j}(\x);
    \end{equation} and the bracket
    \begin{equation} \label{eqn:bracket}
    \langle d f, d g \rangle_t (\x)
    := \sum_{i,j=1}^2 \gamma_{ij}(\x,t)
        \frac{\partial f}{\partial x_i}(\x)
        \frac{\partial g}{ \partial x_j}(\x).
    \end{equation}
If $f, g \in C^{2,1}(\R^2 \times \R ;\R),$ then $(\LL_t f)(\x,t)$ and $\langle d f, d g \rangle_t (\x,t)$ are
obtained by applying $\LL_t$ and $\langle d f, d g \rangle_t$ to $f( \cdot, t) $ and $g(\cdot, t),$ i.e.,
$\left(\LL_t f \right)(\x,t) = \left(\LL_t f( \cdot,t)\right)(\x)$ and $\langle d f, d g \rangle_t f( \x, t)
= \left(\langle d f, d g \rangle_t ( \cdot,t)\right)(\x).$
\end{definition}

\begin{remark}
 If $h \in C^2 \left(\R;\R \right)$ and $f \in C^2(\R^2;\R)$ then
    \begin{equation*}
    \biggl(\LL_t \left( h \circ f \right)\biggr) (\x)
    = h'(f(\x)) (\LL_t f)(\x) + \frac{1}{2} h''(f(\x))
        \langle d f, d f \rangle_t(\x).
    \end{equation*}
\end{remark}

Applying It\^o's formula, from (\ref{eqn:1}) we have that
    \begin{equation*}
    \begin{split}
    H(\widetilde{\X}_t^\eps) & = H(\x) + \eps \int_0^t
        \left( \nabla H(\widetilde{\X}^\eps_s),b(\widetilde{\X}_s^\eps, s)
        \right)d s \\ & \quad + \eps \int_0^t
        \left( \nabla H(\widetilde{\X}_s^\eps), \s (\widetilde{\X}_s^\eps, s)
        \right) d \W_s + \eps^2 \int_0^t (\LL_s H)(\widetilde{\X}^\eps_s) d s,
        \end{split}
    \end{equation*}
and we see that $H(\widetilde{\X}^\eps)$ is slowly varying (as opposed to being conserved). We want to choose a time
scale fine enough to see the effect of the periodic coefficient $b$ on the fluctuations of $H (\widetilde{\X}^\eps)$.
In order to do this, we rescale time by introducing the process $\X^\eps$ so that $\X^\eps_t
=\widetilde{\X}^\eps_{t/\eps^2}$ for every $t \geq 0.$ Using the scaling property of Brownian motion, without loss of
generality, we can rewrite~(\ref{eqn:1}) in the form
\begin{align}
    & d \X_t^\eps
        = \frac{1}{\eps}
            b \left(\X_t^\eps, \frac{t}{\eps^2}\right) \ d t
            + \s \left(\X_t^\eps, \frac{t}{\eps^2}\right)
            \  d \W_t, \ t \geq 0 \notag \\
    & \X^\eps_0  = \x \ \text{a.s.} \label{eqn:1-rescaled}
    \end{align}
We can distinguish between three time scales here. The motion described by~(\ref{eqn:1-rescaled}) consists of the very
fast (of order $1/\eps^{2}$) periodic oscillation of the coefficients, the fast (of order $1/\eps$) rotation along the
trajectories of the corresponding unperturbed system, and the slow diffusion (of order 1) across these trajectories.
From~(\ref{eqn:1-rescaled}), using It\^o's formula,
    \begin{equation*}
    \begin{split}
    H(\X_t^\eps) & = H(\x)+ \frac{1}{\eps} \int_0^t \left( \nabla H(\X^\eps_s),
            b \left(\X_s^\eps, \frac{s}{\eps^2}\right) \right)d s\\
        & \quad +  \int_0^t \left( \nabla H(\X_s^\eps),
            \s \left(\X_s^\eps, \frac{s}{\eps^2}\right) \right) d \W_s
            + \int_0^t (\LL_{s/\eps^2} H) \left(\X^\eps_s \right) d s.
        \end{split}
    \end{equation*}

After rescaling, we can rephrase our goal. We want to study the behavior of $H \left({\X}^\eps\right)$ as
$\eps \downarrow 0,$ and to show that the process describing the energy of the limiting system is Markovian.

\subsection{Topological identification and stratification of the quotient space}
\noindent
 Let $\xi_t$ be the flow generated by $\nabla^\perp H,$ i.e.,
    \begin{equation*}
    \dot{\xi}_t(\x)  = \nabla^{\perp} H(\xi_t(\x)), \quad \xi_0(\x) =\x.
    \end{equation*}
We will use $\xi_t$ to form a quotient space. Let $\sim$ be the chain equivalence relation based on $\xi_t$.

All trajectories of $\xi_t$ are periodic on $\R^2 \setminus \overV$, therefore the chain orbit of $\x \in
\R^2 \setminus \overV$ consists of all the points on the level set to which $\x$ belongs. All points of the
open set $\V$ are fixed points of the flow $\xi_t$, and every point $\x \in \V$ is its own chain orbit. All
points of $\partV$ are also fixed points of the flow, but the chain orbit of every such point is the whole
set $\partV$, since every point of $\partV$ always has a neighborhood containing points from $\R^2 \setminus
\overV.$

\begin{figure}[h]
\centerline{\psfig{file=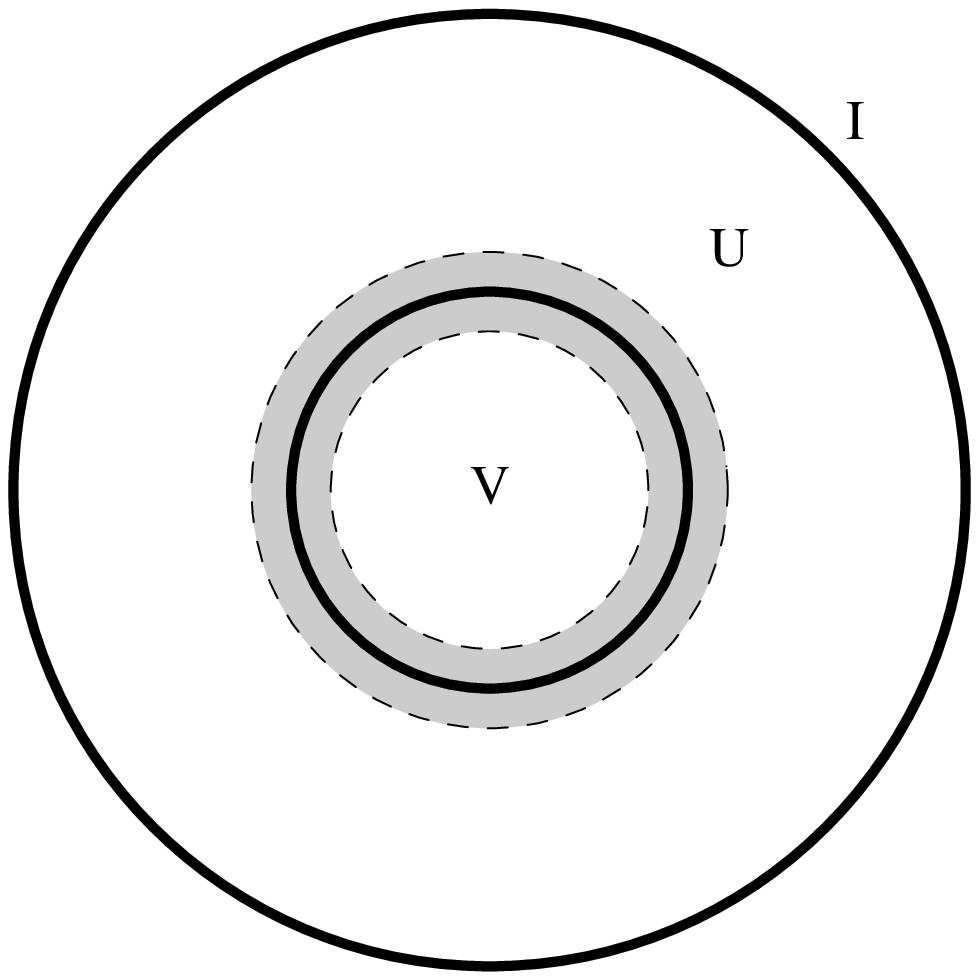,width=2in}} \vspace*{6pt} \fcaption{$\overI = \V \bigcup \partial \V
\bigcup \U \bigcup \partI$.} \label{fig:setU}
\end{figure}

Fix a level of $H$ of height $H^*>0,$ let $K^* := \left( H^*\right)^{1/n}$, and define the two open sets
    \begin{equation*}
    \U := \{\x \in \R^2: 0< H(\x) < H^*\} = \{\x \in \R^2: 0< K(\x) < K^*\}
    \end{equation*}
and
    \begin{equation*}
    \I := \{\x \in \R^2: H(\x) < H^*\} = \{\x \in \R^2: K(\x) < K^*\}.
    \end{equation*}
The sets $\V, \partial \V, \U, \partI$ are disjoint with union $\overI$ (Fig.~\ref{fig:setU}).

Define
    \begin{equation*}
     \G := \overline{\I} / \sim.
     \end{equation*}
Then $\G,$ endowed with the quotient topology defined by $\sim$, is a quotient space. For every $\x \in
\overI,$ define the natural map
    $$\pi(\x) :=[\x] = \left\{ \y \in \overI: \y \sim \x \right\}.$$
Notice that
    \begin{equation*}
    \pi(\x) =
  \begin{cases}
    H^{-1} \left(H(\x)\right)
        & \text{ if } \x \in \U \cup \partial \I, \\
    \partial \V
        & \text{ if } \x \in \partial \V, \\
    \{\x\} & \text{ if } \x \in \V.
  \end{cases}
    \end{equation*}
For $A \subset \overI$ denote $\G_{A}  = \pi(A) = \left\{ [\x] \in \G: \x \in A \right\} $ so that $ \G =
\G_{\V} \cup \G_{\partial \V} \cup \G_{\U} \cup \G_{\partial \I}. $

\begin{proposition}
There is a homeomorphism from $\G$ to the submanifold $${\sfS}^2 \cup \biggl( \{0 \} \times \{0 \} \times
\left[1,H^*+1 \right] \biggr)$$ in $\R^3$, where ${\sfS}^2$ is the $2$-sphere $\left\{\x \in \R^3: \| \x \|=1
\right\}.$
\end{proposition}
\begin{proof}[\bf{Proof}]
Denote ${\sfS}^2 \cup \biggl( \{0 \} \times \{0 \} \times \left[1,H^*+1 \right] \biggr)$ by ${\sfM}_3,$ and define the
map~$p: \G \rightarrow {\sfM}_3$ by
\begin{equation*}
p \left([\x] \right) :=
  \begin{cases}
    s \left(\x \right) & \text{ if } [\x] \in \G_{\V}, \\
    (0,0,1) & \text{ if } [\x] \in \G_{\partial \V},\\
     \left(0,0,H(\x)+1 \right)
     & \text{ if } [\x] \in \G_{\U} \cup \G_{\partial \U},
  \end{cases}
\end{equation*}
where $s: \overV \rightarrow {\sfS}^2$ is a homeomorphism such that $s \left(\partV \right) = (0,0,1).$ Notice that $p$
is the desired homeomorphism from $\G$ to ${\sfM}_3.$
\end{proof}

We see that $\G_{\V}$ is an open two-dimensional smooth manifold, $\G_{\U}$ is an open one-dimensional smooth
manifold, $\G_{\V} \cap \G_{\U} = \varnothing.$ $\G_{\partV}$ is a limit point of both $\G_{\V}$ and
$\G_{\U}$, and $\G_{\partI}$ is the other limit point of $\G_{\U}.$ Both $\G_{\partV}$ and $\G_{\partI}$ can
be treated as zero-dimensional smooth manifolds, and the Whitney regularity condition is trivially satisfied.
Therefore, $\G = \G_{\V} \cup \G_{\partial \V} \cup \G_{\U} \cup \G_{\partial \I}$ is a stratified space.

\subsection{Probabilistic formulation of the problem}
\noindent Probabilistically, we are interested in the laws of a stochastic processes. We would like to show
that a certain limiting law is Markovian in nature. The martingale problem is a standard tool for such
analysis. We will restate our problem in the language natural for the martingale approach.

For any $f \in C^{2}(\R^2),$ and for every $t \in \R,$ define the differential operator
    \begin{equation*}
    (\LL_t^\eps f)(\x)
    := \frac{1}{\eps} \left( \nabla f(\x), b(\x,t) \right)
    + (\LL_t f)(\x),
    \end{equation*}
with $\LL_t f$ defined in (\ref{eqn:L}) above. If $f \in C^{2,1}(\R^2 \times \R), $ then $(\LL_t^\eps
f)(\x,t)$ is obtained by applying $\LL_t^\eps$ to $f( \cdot, t) $ i.e., $\left(\LL_t f \right)(\x,t) =
\left(\LL_t f( \cdot,t)\right)(\x)$. We assume that the domain $\D^\eps$ of the generator $\LL_t^\eps$
contains a dense subset of functions from $C^{2}\left(\R^2 \right)$.

Let $\Om :=C([0,\infty);\R^2).$ Given $\x \in \overI,$ let $\PP^\eps \in \mathcal{P}\left(C \left([0,\infty);
\R^2\right) \right)$ be a solution of the stopped martingale problem for $ \left(\LL_t^\eps,\delta_{\x}, \I \right), t
\geq 0.$ Then for every $f \in C^{2,1} \left(\R^2 \times [0,\infty) \right)$,
\begin{equation*}
    f \left(\X_{t \wedge \tau},t \wedge \tau \right)
    - f \left(\X_{s \wedge \tau}, s \wedge \tau \right)
     - \int_{s \wedge \tau}^{t \wedge \tau} \biggl( \frac{\partial
     f}{\partial u}\left(\X_u,u\right) +
     \left(\LL_{u}^\eps f \right) \left(\X_u,u\right)\biggr) d u
    \end{equation*}
is a $\PP^\eps$-martingale.

Let $\X$ be the corresponding coordinate process defined by $ \X_t(\omega):=\omega(t)$, $t\geq~0,~\omega \in \Om.$ Set
$\F_t$ to be generated by this coordinate process up to time $t$, and define $\F:= \bigvee_t \F_t.$ Define the
$\F_t$-stopping time $\tau$ by
    \begin{equation*}
     \tau := \inf \left\{s \geq 0: \X_s \notin \I \right\}
     =  \inf \left\{s \geq 0: H(\X_s) \geq H^* \right\},
     \end{equation*}
i.e., $\tau$ is defined to be the time of the first exit from $\I.$

Our condition on $\PP^\eps$ means that the corresponding coordinate process $\X$ on $ \left(\Om,\F,\PP^\eps
\right)$ is a solution to the stopped martingale problem for $ \left(\LL_t^\eps,\delta_{\x}, \I \right), t
\geq 0.$

Define $\Y_t:=\left[\X_{t \wedge \tau} \right], \ t \geq 0.$ For each $\eps > 0,$ define the probability
measure $\PP^{\eps,*}(A) := \PP^\eps \left\{ \Y \in A \right\}, \quad A \in \mathcal{B} \left(C \left( [0,
\infty), \G \right) \right).$ Let $\Y_t(\omega)=\omega(t), t \geq 0,$ $ \omega~\!\in~\!C~\!\left( [0,
\infty), \G \right),$ be the corresponding coordinate process. For each $t \geq 0, $ define $\F^*_t:= \s
\left\{\Y_s, 0 \leq s \leq t \right\},$ and define a $\sigma$-algebra on $\Omega^*=C\left( [0, \infty), \G
\right)$ by $\F^*:= \bigvee_{t \geq 0} \F_t^*$.

With this in mind, we can formulate our task once again; this time in the language for the martingale
approach. Our goal is to show that the $\PP^\eps$-laws of $\Y$ converge to the law of a $\G$-valued Markov
process. In other words, we want to understand the asymptotics of $\PP^{\eps,*}$ as $\eps \downarrow 0$. Even
more precisely, we will prove that $\PP^{\eps,*}$ converges to the unique solution $\PP^*$ of the martingale
problem for $\left(\LL^*, \delta_{[\x]} \right)$, and we will identify the limiting generator $\LL^*$ and its
domain.

\section{Averaging}
\noindent
The main goal of our work is model reduction. To achieve this goal, we use both time averaging and
the stochastic averaging of Freidlin and Wentzell as our main techniques.

We study our system in the presence of three time scales. The periodic fluctuations of the coefficients occur
on the time scale of order $1/\eps^2$, the effect of the drift is visible in time of order $1/\eps$, and
diffusion is of order one. The hierarchy of the time scales indicate that our system will require
multiple-level averaging. First, we average the quickly-varying periodic coefficients.

\subsection{First-level averaging}

We recall that $\PP^\eps \in \mathcal{P}(C[0,\infty); \R )$ is a solution of the stopped martingale problem
for $ \left(\LL_t^\eps,\delta_{\x}, \I \right)$, $t \geq 0,$ $\x \in \overI,$ where $ \I := \{\x \in \R^2:
H(\x) < H^*\},$ and $ \tau := \inf \left\{s \geq 0: \X_s \notin \I \right\}.$ Let $\E^\eps$ be the
expectation operator associated with $\PP^\eps.$
\begin{lemma} \label{lemma:f-Mf}
If $f \in C^2 \left(\overI \times  [0, \infty); \R \right)$ is $2 \pi$-periodic in its last argument, then
there exists a positive constant $C$ such that
    \begin{equation}
    \E^\eps\left[ \left|
        \int_{0}^{t \wedge \tau}
        \left( f
        \left( \X_u, \frac{u}{\eps^2} \right) - (\M f) \left( \X_u \right)
        \right) d u \right|
            \right]
     \leq C \eps (1+t)
        \left\| f \right\|_{C^2 \left(\overI \times  [0, \infty) \right)}. \label{eqn:aver}
    \end{equation}
\end{lemma}

\begin{proof}[\bf{Proof}]
We define
    \begin{equation} \label{eqn:Phi-phi}
    \Phi_{f}(\x,t) : =  \int_0^t
        \left(f(\x,u) - (\M f)(\x) \right) d u.
    \end{equation}

For any $\x \in \overI, $ any multi-index $\alpha$ and $t \geq 0,$
    \begin{equation*}
    \begin{split}
    \partial^\alpha \Phi_f (\x,t)
    & = \int_0^t \left(
            \partial^\alpha f \left( \x,u \right)
            - \frac{1}{2 \pi} \int_0^{2 \pi}  \partial^\alpha f(\x,s) d s
                \right) d u\\
    & = \int_{2 \pi \lfloor \frac{t}{2 \pi } \rfloor}^t
            \partial^\alpha f \left( \x,u \right) d u
    - \left(\frac{t}{2 \pi } - \lfloor \frac{t}{2 \pi } \rfloor \right)
    \int_0^{2 \pi} \partial^\alpha f \left( \x,u \right) d u.
    \end{split}
    \end{equation*}
    This implies that for every $0<k \leq 2$
    \begin{equation} \label{ineqn:Phif-f}
    \left\| \Phi_f\right\|_{C^k \left(\overline{\I} \times  [0, \infty)\right)}
    \leq 4 \pi \left\| f\right\|_{C^k \left(\overline{\I} \times  [0, \infty) \right)}.
    \end{equation}
 By the stopped martingale problem for $ \left(\LL_{t/\eps^2}^\eps,\delta_{\x}, \,\I \right),$
    \begin{equation} \label{eqn:Phif}
    \begin{split}
    \Phi_f\left(\X_{t \wedge \tau},\frac{t \wedge \tau}{\eps^2} \right)
    - \Phi_f\left(\x,0\right)
    & = \frac{1}{\eps^2} \int_0^{t \wedge \tau}
    \frac{\partial \Phi_f}{\partial u}\left(\X_u, \frac{u}{\eps^2}\right) d u\\
    & \quad +  \frac{1}{\eps} \int_0^{t \wedge \tau}
        \left( \nabla_{\x} \Phi_f \left(\X_u, \frac{u}{\eps^2}\right),
            b \left(\X_u, \frac{u}{\eps^2}\right) \right)d u\\
    & \quad +  \int_0^{t \wedge \tau}
        \left(\LL_{u/\eps^2} \Phi_f \right) \left(\X_u, \frac{u}{\eps^2}\right) d u\\
    & \quad + \mathbf{M}^{\Phi_f,\eps}_{t \wedge \tau},
    \end{split}
    \end{equation}
where $\mathbf{M}^{\Phi_{f},\eps}$ is the $\PP^\eps$-martingale  with its quadratic variation
    \begin{equation*}
    \langle \mathbf{M}^{\Phi_f,\eps} \rangle_{t}
    = \int_0^{t } \langle d \Phi_f, d \Phi_f \rangle_{u/\eps^2}
            \left( \X_u, \frac{u}{\eps^2} \right)d u.
    \end{equation*}
    Then
    \begin{equation} \label{eqn:avg}
    \begin{split}
    \int_{0}^{t \wedge \tau}
        \left( f \left( \X_u, \frac{u}{\eps^2} \right) - (\M f) \left( \X_u \right)
        \right) d u
    & = \eps^2 \left( \Phi_f\left(\X_{t \wedge \tau},\frac{t \wedge \tau}{\eps^2} \right)
        - \Phi_f\left(\x,0\right) \right)\\
    & \quad - \eps \int_0^{t \wedge \tau}
        \left( \nabla_{\x} \Phi_f \left(\X_u, \frac{u}{\eps^2}\right),
            b \left(\X_u, \frac{u}{\eps^2}\right) \right)d u\\
    & \quad -  \eps^2 \int_0^{t \wedge \tau}
        \left(\LL_{u/\eps^2} \Phi_f \right) \left(\X_u, \frac{u}{\eps^2}\right) d u\\
    & \quad - \eps^2 \ \mathbf{M}^{\Phi_f,\eps}_{t \wedge \tau}.
    \end{split}
    \end{equation}
By Burkholder-Davis-Gundy inequality, there exists a universal positive constant $C_m$ such that
\begin{equation} \label{ineq:BDG}
    \E^\eps \left[\left|\mathbf{M}^{\Phi_f,\eps}_{t \wedge \tau}\right|\right]
    \leq C_m \,\E^\eps \left[ \langle \mathbf{M}^{\Phi_f,\eps} \rangle^{1/2}_{t \wedge \tau}\right]
    \end{equation}
    From (\ref{ineqn:Phif-f}), there exists a positive constant $C$ such that
\begin{equation} \label{ineq:qv}
    \langle \mathbf{M}^{\Phi_f,\eps} \rangle_{t \wedge \tau}
    \leq C t \left\|f \right\|^2_{C^1 \left(\overline{\I} \times  [0, \infty)\right)}.
    \end{equation}
Finally, from expression (\ref{eqn:avg}), applying inequalities (\ref{ineqn:Phif-f}), (\ref{ineq:BDG}), and
(\ref{ineq:qv}), there is a positive constant $C$ such that (\ref{eqn:aver}) holds.
\end{proof}

\subsection{Stochastic averaging away from the critical set}
\noindent
Let $\phi_t$ be the flow generated by $\nabla^\perp K$, i.e.,
    \begin{equation*}
     \dot{\phi}_t(\x) =\nabla^\perp K \left(\phi_t(\x) \right),\quad \phi_0(\x) = \x.
    \end{equation*}
Here we are concerned with stochastic averaging on a subset of $\overI \setminus \overV = \U \cup
\partI$. On any such subset, $\xi_t(\x) = \phi_{n \left(K(\x)\right)^{n-1} t}(\x)$ and the flow $\phi_t$ has the same orbits as the
flow $\xi_t$ generated by $\nabla^{\perp} H.$  To average over these orbits, we define
    \begin{equation*}
    \eta(K(\x)) : = \inf \left\{ t > 0: \phi_t(\x) = \x \right\},
    \end{equation*}
    and for every $f \in C\left(\overI;\R \right)$ we define
    \begin{equation} \label{def:Psi-f}
    \Psi_f(\x) := \frac{1}{\eta \left(K(\x)\right)}
        \int_0^{\eta \left(K(\x)\right)} t  f\left(\phi_t(\x) \right) d t.
    \end{equation}

Let $\D_{\A_0}$ be the set of all real-valued functions whose restriction to $\V \cup \U$ is continuous. We
define the pre-averaging operator $\A_0 : \D_{\A_0} \rightarrow C \left(\G_{\V} \cup \G_{\U}\right)$ by
    \begin{equation*}
    \left(\A_0 f\right) ([\x]):= \begin{cases}
    f(\x) & \text{ if } [\x] \in \G_{\V}\\
    \left(\int_{\y \in [\x]} \frac{1}{\|\nabla K(\y)\|} \ dl(\y) \right)^{-1}
    \int_{\y \in [\x]} \frac{f(\y)}{\|\nabla K(\y)\|} \ dl(\y) & \text{ if } [\x] \in \G_{\U},
    \end{cases}
    \end{equation*}
where $dl$ is the arc-length element.

Let $\D_{\A}$ be the subset of functions $f$ from $\D_{\A_0}$ for which
 \begin{equation*}
    \lim_{\substack{[\y] \rightarrow \G_{\partV}\\
        [\y] \in \G_\V  \cup \G_{\U}}}
        \left(\A_0 f\right) ([\y])
        \quad
        \text{and}
        \quad
        \lim_{\substack{[\y] \rightarrow \G_{\partU}\\
        [\y] \in \G_\V \cup \G_{\U} }}
        \left(\A_0 f\right) ([\y])
    \end{equation*} exist. We define the averaging operator $\A:\D_{\A} \rightarrow C(\G)$ by
    \begin{equation*}
    (\A f)([\x]) := \lim_{\substack{[\y] \rightarrow [\x]\\
        [\y] \in \G_\V \cup \G_{\U} }}
    (\A_0 f)([\y]).
    \end{equation*}
\begin{remark}
We notice that if $\x \in \U \cup \partI,$ then $[\x]=K^{-1}\left(K(\x)\right),$ i.e., $[\x]$ is totally
defined by $K(\x),$ and $ \left(\A f\right)([\x]) = \left(\A f \right) \left(K^{-1}(K(\x)) \right). $ We
define $ \left(\A_K f\right)(h)$ to be the average of $f$ on the level set $\left\{\x \in \R^2: K(\x)=h
\right\}$, i.e.,
\begin{equation*}
\left(\A_K  f\right) (K(\x))
    := \frac{1}{\eta(K(\x))}
    \int_{0}^{\eta(K(\x))}  f(\phi_t(\x)) d t.
\end{equation*} Then for any $\x \in \U \cup \partI,$
\begin{equation*}
\left(\A f\right)([\x]) = \left(\A_K f \right) \left(K(\x) \right).
\end{equation*}
\end{remark}

Let the function $\omega \in C_c^\infty(\R)$ be supported on a compact subset of $(0, \infty).$
\begin{lemma} \label{lemma:f-Af}
For any function $f \in C^3 \! \left( \overI ; \R \right)$ and for any $0< p < 1$ there exists a constant
$C_a>0$ such that
    \begin{align}
    \E^\eps\left[ \left|
        \int_{0}^{t \wedge \tau}
         \biggl( f
        \left( \X_u \right) - (\A f) \left(\left[ \X_u \right] \right)
        \biggr)
        \right. \right. & \left. \left.
        \omega \left( \frac{K \left(\X_u \right)}{\eps^p} \right) d u \right|
            \right] \notag\\
     & < C_a \ \eps^{1-p(n+2)}
    \ (1+t) \ \|f\|_{C^3 \left(\overline{\I} \right)}. \label{ineq:f-Af}
    \end{align}
\end{lemma}
\begin{proof}[\bf{Proof}] Let $\x \in \overI.$ Since $\omega$ is supported on a compact subset of $(0, \infty),$ either
\begin{equation*}
\biggl( f \left( \x \right) - (\A f) \left(\left[ \x \right] \right) \biggr) \omega \left( \frac{K \left(\x
\right)}{\eps^p} \right)=0, \end{equation*} or there exist some positive constants $C_1$ and $C_2$ such that
$C_1 \eps^p < K(\x)< C_2 \eps^p$, and therefore we may assume that $\x \in \U \cup \partI$. From the above
remark, this implies that $ \left(\A f\right)([\x]) = \left(\A_K f \right)(K(\x)).$

From (\ref{def:Psi-f}), there exist a constant $C>0$ such that
    \begin{equation} \label{ineq:Psif}
    \left\|\Psi_f\right\|_{C^k\left(\overI \right)}
    < C \left\|f\right\|_{C^k\left(\overline{\I} \right)}.
    \end{equation}
Notice that $\eta \left( K\left(\phi_t(\x)\right)\right) = \eta \left(K(\x)\right)$ and
    \begin{equation*}
    \begin{split}
    \Psi_f\left(\phi_t(\x) \right)&  = \frac{1}{\eta \left(K(\x) \right)}
        \int_0^{\eta \left(K(\x) \right)} s
         f\left(\phi_s \left( \phi_t(\x)\right) \right) d s \\
       & = \frac{1}{\eta \left(K(\x) \right)}
        \int_t^{t+\eta \left(K(\x) \right)} (s- t)
        f\left(\phi_s (\x) \right) d s.
        \end{split}
    \end{equation*}
    Differentiating the last equality, we obtain
\begin{align}
    \frac{d}{d t} \Psi_f\left(\phi_t(\x) \right)
        & = \frac{1}{\eta \left(K(\x) \right)} \frac{d}{d t}
        \left(\int_t^{t+\eta \left(K(\x) \right)} (s-t) f\left(\phi_s (\x) \right) d s \right) \notag\\
        & = f \left(\phi_t(\x)\right)
        - (\A_K f) \left( K \left(\phi_t(\x)\right)\right). \label{eqn:derPsi1}
    \end{align}
On the other hand, using the fact that the flow $\phi_t$ is generated by $\nabla^{\perp} K,$ we observe that
    \begin{equation} \label{eqn:derPsi2}
    \frac{d}{d t} \Psi_f \left(\phi_t(\x) \right)
    = \left( \nabla  \Psi_f \left(\phi_t(\x) \right),
    \nabla^{\perp} K \left(\phi_t(\x) \right) \right),
    \end{equation}
Expressions (\ref{eqn:derPsi1}) and (\ref{eqn:derPsi2}) combined yield
    \begin{equation*}
  f \left(\phi_t(\x)\right)
        - (\A_K  f) \left(K\left(\phi_t(\x)\right)\right)
        = \left( \nabla  \Psi_f \left(\phi_t(\x) \right),
    \nabla^{\perp} K \left(\phi_t(\x) \right) \right).
    \end{equation*}
This equality, being true for any $t \in \R$, implies that for any $\x \in \overI,$
    \begin{equation} \label{eqn:f-A0f=grad-grad}
    \biggl( f (\x) - (\A_K f) \left( K\left(\x\right)\right) \biggr)
    \omega \left( \frac{K \left(\x \right)}{\eps^p} \right)
    = \left( \nabla  \Psi_f (\x), \nabla^{\perp} K(\x) \right)
    \omega \left( \frac{K \left(\x \right)}{\eps^p} \right).
    \end{equation}
For any $\x \in \overI$, we define
    \begin{equation*}
    \Psi_{f}^\eps(\x): = \frac{1}{n (K(\x))^{n-1}} \Psi_f(\x) \omega \left(\frac{K(\x)}{\eps^p}\right)
    \end{equation*}
and notice that since $\left(\M b\right)(\x) = n (K(\x))^{n-1} \nabla^\perp K(\x)$, then
\begin{equation*}
     \biggl(\M \left(\nabla \Psi_f^\eps,b\right) \biggr)(\x)
     = \left( \nabla  \Psi_f (\x), \nabla^{\perp} K(\x) \right)
    \omega \left( \frac{K \left(\x \right)}{\eps^p} \right).
    \end{equation*}
Using inequality (\ref{ineq:Psif}), we see that there exists a constant $C>0,$ such that
\begin{equation} \label{ineq:Psi-epsf}
    \left\|\Psi^\eps_f\right\|_{C^k\left(\overline{\I} \right)}
    < C \eps^{-p(n-1+k)} \left\|f\right\|_{C^k \left(\overline{\I} \right)}.
    \end{equation}
This implies that there exists a positive constant $C$ such that
    \begin{equation} \label{ineq:grad-Psi-b-epsf}
    \left\| \left( \nabla \Psi^\eps_f, b \right)
    \right\|_{C^2 \left(\overline{\I} \times [0, \infty)\right)}
    < C \eps^{-p \ (n+2) } \left\|f\right\|_{C^3 \left(\overline{\I} \right)}.
    \end{equation}

By the stopped martingale problem for $ \left(\LL^\eps_{t/\eps^2},\delta_{\x}, \I \right), $
    \begin{multline} \label{eqn:Phi-f}
    \int_0^{t \wedge \tau}
        \left( \nabla \Psi^\eps_f \left(\X_u \right),
            b \left(\X_u, \frac{u}{\eps^2}\right) \right)d u
    = \eps \left(\Psi^\eps_f\left(\X_{t \wedge \tau}\right)
        - \Psi^\eps_f\left(\x \right)\right)\\
     - \eps \int_0^{t \wedge \tau}
        \left(\LL_{u/\eps^2} \Psi^\eps_f \right) \left(\X_u\right) d u
    - \eps \mathbf{M}^{\Psi^\eps_f,\eps}_{t \wedge \tau},
    \end{multline}
where $\mathbf{M}^{\Psi^\eps_{f},\eps}$ is the $\PP^\eps$-martingale with its quadratic variation
    \begin{equation*}
    \langle \mathbf{M}^{\Psi^\eps_f,\eps} \rangle_{t}
    = \int_0^{t } \langle d \Psi^\eps_f, d \Psi^\eps_f \rangle_{u/\eps^2}
            \left( \X_u \right)d u.
    \end{equation*}
From (\ref{eqn:Phi-f}), (\ref{ineq:Psi-epsf}) and Burkholder-Davis-Gundy inequality, we see that there is a
positive constant $C$ such that
    \begin{align}
    \E^\eps &
    \Biggl[ \Biggl|
        \int_0^{t \wedge \tau}
        \Biggl(
            \nabla \Psi^\eps_f \left(\X_u \right),
            b \left(\X_u, \frac{u}{\eps^2} \right) \Biggr) d u
        \Biggr|
    \Biggr]
      \leq
    \eps \E^\eps
    \left[
        \left|
        \Psi^\eps_f\left(\X_{t \wedge \tau}\right)
        - \Psi^\eps_f\left(\x \right)
        \right|
    \right]\notag\\
    & \quad + \eps \E^\eps
    \left[
        \left|
        \int_0^{t \wedge \tau}
        \left(\LL_{u/\eps^2} \Psi^\eps_f \right) \left(\X_u\right) d u
        \right|
    \right] + \eps \E^\eps
    \left[
        \left|\mathbf{M}^{\Psi^\eps_f,\eps}_{t \wedge \tau} \right|
    \right] \notag\\
    & < C \eps^{1-p(n+1)} \ (1+t) \|f\|_{C^2 \left(\overline{\I}\right)}. \label{ineq:f-Mf-1}
    \end{align}
Apply Lemma~\ref{lemma:f-Mf}, with $\left(\nabla \Psi^\eps_f,b \right)$ in place of $f$, and then recall
(\ref{ineq:grad-Psi-b-epsf}) to notice that there exists a positive constant $C$ such that
    \begin{multline} \label{ineq:f-Mf-2}
    \E^\eps
    \Biggl[ \Biggl|
    \int_{0}^{t \wedge \tau}
        \Biggl(
            \left( \nabla \Psi^\eps_f \left(\X_u \right),
            b \left(\X_u, \frac{u}{\eps^2}\right) \right)
            - \left(\M \left( \nabla \Psi^\eps_f ,
            b \right) \right) \left(\X_u \right)
        \Biggr) d u
    \Biggr|
    \Biggr]\\
    < C \eps^{1-p (n+2)}(1+t)
    \left\| f \right\|_{C^3 \left(\overline{\I} \right)}.
        \end{multline}
Finally, combine (\ref{ineq:f-Mf-1}) and (\ref{ineq:f-Mf-2}) into (\ref{ineq:f-Af}).
\end{proof}

\begin{lemma} \label{lemma:f-AMf}
For any function $f \in C^3\left(\overI \times  [0, \infty); \R \right)$ which is $2 \pi$-periodic in its
last argument and for any $0 < p < 1$ there exists a constant $C_c>0$ such that
    \begin{multline*}
    \E^\eps \left[ \left|
        \int_{0}^{t \wedge \tau}
        \biggl( f
        \left( \X_u , \frac{u}{\eps^2}\right) - (\A (\M f)) \left(\left[ \X_u \right] \right)
        \biggr)
        \omega \left( \frac{K \left(\X_u \right)}{\eps^p} \right) d u \right|
            \right]\\
    < C_c \eps^{1-p(n+2)} (1+t)
        \left\| f \right\|_{C^3 \left(\overline{\I} \times  [0, \infty) \right)}.
        \end{multline*}
\end{lemma}
\begin{proof}[\bf{Proof}] We denote
    \begin{equation*}
    q_\eps(\x,t): = f(\x,t) \omega \left( \frac{K(\x)}{\eps^p} \right), \quad \x
    \in \overI,
    \quad t \geq 0.
    \end{equation*}
Notice that $q_\eps \in C^2 \left(\overI \times  [0, \infty); \R \right)$ is $2 \pi$-periodic in its last
argument,
    \begin{equation} \label{ineq:q-eps}
    \left\|q_\eps\right\|_{C^2 \left(\overI \times  [0, \infty)\right)}
     \leq C \eps^{-2p}\,\left\|f \right\|_{C^2 \left(\overI \times  [0, \infty)\right)}
     \end{equation} for
some constant $C>0,$ and
    \begin{multline} \label{eq:f-A0f}
    \biggl(f \left( \X_u , \frac{u}{\eps^2}\right)
    - (\A (\M f)) \left(\left[ \X_u \right] \right) \biggr)
    \omega \left( \frac{K \left(\X_u \right)}{\eps^p} \right)
     = \Biggl( q_\eps \left( \X_u , \frac{u}{\eps^2}\right)
        - (\M q_\eps)(\X_u) \Biggr)\\
        + \biggl( (\M f)\left( \X_u \right)
        - (\A(\M f))\left(\left[ \X_u \right] \right) \biggr)
        \omega \left( \frac{K \left(\X_u \right)}{\eps^p} \right) .
    \end{multline}
Applying Lemma~\ref{lemma:f-Mf}, and using (\ref{ineq:q-eps}), we obtain
    \begin{equation} \label{ineq:qeps-M-qeps}
    \E^\eps\left[ \left|
        \int_{0}^{t \wedge \tau}
        \left( q_\eps
        \left( \X_u, \frac{u}{\eps^2} \right) - (\M q_\eps) \left( \X_u \right)
        \right) d u \right|
            \right]
    < C \eps^{1-2p}\, (1+t)
        \left\| f \right\|_{C^2 \left(\overline{\I} \times  [0, \infty) \right)}.
    \end{equation}
The function $\M f \in C^3 \left(\overI ; \R \right)$, and
    $
    \|\M f\|_{C^3 \left(\overI \right)} \leq \|f\|_{C^3 \left(\overI \times  [0, \infty) \right)}.
    $
Taking $\M f$ in place of $f$ in Lemma~\ref{lemma:f-Af} we have
    \begin{multline} \label{ineq:f-A0f}
    \E^\eps \left[ \left|
        \int_{0}^{t \wedge \tau}
        \biggl( (\M f)\left( \X_u \right)
        - (\A(\M f))\left(\left[ \X_u \right] \right) \biggr)
        \omega \left( \frac{K \left(\X_u \right)}{\eps^p} \right)
       d u \right|
            \right]\\
       < C_a \eps^{1-p(n+2)}(1+t)
       \|f\|_{C^3 \left(\overline{\I} \times  [0, \infty) \right)}.
       \end{multline}
Combine (\ref{eq:f-A0f}), (\ref{ineq:qeps-M-qeps}) and (\ref{ineq:f-A0f}) to complete this proof.
\end{proof}

Denote $\psi(\x,t) := \left( \nabla K(\x), b(\x,t) \right),$ then $\M \psi \equiv 0$ and by definition
(\ref{eqn:Phi-phi}),
        \begin{equation*}
        \Phi_\psi(\x,t) = \int_0^t
        \left( \nabla K(\x), b(\x,u) \right)d u.
        \end{equation*}

\begin{lemma} \label{lemma:phi+grad}
For any function $g \in C^3\left(\R; \R \right)$ and for any $0<p<1$ there exists a constant $C_b >0$ such
that
    \begin{multline*}
    \E^\eps\Biggl[ \Biggl|
        \int_{0}^{t \wedge \tau}
        \biggl\{\frac{1}{\eps}\psi \left( \X_u , \frac{u}{\eps^2}\right)
        + \left(\nabla \Phi_{\psi}, b \right)\left( \X_u , \frac{u}{\eps^2}\right)
        \biggr\} g \left( K \left(\X_u \right)\right)
        \omega \left(\frac{K(\X_u)}{\eps^p} \right)d u \Biggr|
            \Biggl]\\
    < C_b \,\eps^{1-3p} \, (1+t) \left\| g \right\|_{C^3 \left([0,K^*] \right)}.
    \end{multline*}
    \end{lemma}
    \begin{proof}[\bf{Proof}]
Let $k_\eps(h) = g(h) \,\omega \! \left(\frac{h}{\eps^p} \right), \ h \in \R.$ By the stopped martingale
problem for $ \left(\LL_{t/\eps^2}^\eps,\delta_{\x}, \I \right), t \geq 0, $
    \begin{align}
    & \eps \Phi_{\psi} \left( \X_{t \wedge \tau} , \frac{t \wedge \tau}{\eps^2}\right)
    k_\eps \left( K \left(\X_{t \wedge \tau} \right)\right)
    - \eps \Phi_{\psi}\left( \x , 0\right)
    k_\eps \left( K \left(\x \right)\right) \notag\\
    & - \frac{1}{\eps} \int_0^{t \wedge \tau} \frac{\partial \Phi_\psi}{\partial u}
    \left( \X_u , \frac{u}{\eps^2}\right) k_\eps \left( K \left(\X_u \right)\right) d u \notag\\
    & -  \int_0^{t \wedge \tau}
    \left(\nabla \Phi_{\psi}\left( \X_u , \frac{u}{\eps^2}\right),
         b \left( \X_u , \frac{u}{\eps^2}\right)\right)
         k_\eps \left( K \left(\X_u \right)\right) d u \notag\\
    & - \eps \int_0^{t \wedge \tau}
    \left(\LL_{u /\eps^2} \Phi_{\psi} \right) \left( \X_u , \frac{u}{\eps^2}\right)
         k_\eps \left( K \left(\X_u \right)\right) d u \notag\\
    & -  \int_0^{t \wedge \tau} \Phi_{\psi} \left( \X_u , \frac{u}{\eps^2}\right)
         \psi \left( \X_u , \frac{u}{\eps^2}\right)
         k_\eps'\left( K \left(\X_u \right)\right) d u \notag\\
    & -  \eps \int_0^{t \wedge \tau} \Phi_{\psi} \left( \X_u , \frac{u}{\eps^2}\right)
          (\LL_{u/\eps^2} (k_\eps \circ K))
         \left( \X_u , \frac{u}{\eps^2}\right)d u \notag\\
    & - 2 \eps \int_0^{t \wedge \tau} \langle d \Phi_\psi, d (k_\eps \circ K) \rangle_{u/\eps^2}
    \left( \X_u , \frac{u}{\eps^2}\right)
    = \eps M^\eps_{t \wedge \tau}, \label{ineq:comb}
    \end{align}
where $M^\eps$ is a $\PP^\eps$-martingale with its quadratic variation
    \begin{equation*}
    \langle M^\eps \rangle_t= \int_0^t \langle d \left( \Phi_{\psi} \left(k_\eps \circ K \right)\right),
        d \left( \Phi_{\psi} \left(k_\eps \circ K \right)\right)
        \rangle_{u/\eps^2} \left( \X_u , \frac{u}{\eps^2} \right)du.
\end{equation*}

Now notice that $$\frac{\partial \Phi_\psi}{\partial t}(\x,t) =\psi(\x,t),$$ rearrange the terms of
(\ref{ineq:comb}) and apply Burkholder-Davis-Gundy inequality to see that there exists a constant $C>0,$ such
that
    \begin{multline} \label{ineq:part1}
    \E^\eps\Biggl[ \Biggl|
        \int_{0}^{t \wedge \tau}
        \biggl\{\frac{1}{\eps}\psi \left( \X_u , \frac{u}{\eps^2}\right)
        + \left(\nabla \Phi_{\psi}\left( \X_u , \frac{u}{\eps^2}\right),
         b \left( \X_u , \frac{u}{\eps^2}\right)\right)
        \biggr\} k_\eps \left( K \left(\X_u \right)\right) d u \Biggr|
            \Biggl]\\
     \leq \E^\eps \Biggl[ \Biggl|
    \int_0^{t \wedge \tau} \Phi_{\psi} \left( \X_u , \frac{u}{\eps^2}\right)
        \psi \left( \X_u , \frac{u}{\eps^2}\right)
         k_\eps'\left( K \left(\X_u \right)\right) d u \Biggr| \Biggr]\\
          + C \eps^{1-2p} (1+t) \| g \|_{C^2 \left(\left[0,K^*\right]\right)}.
    \end{multline}
Next consider
    \begin{equation*}
    \left(\M \left( \Phi_\psi \psi \left(k_\eps' \circ K \right) \right) \right)(\x)
     =  k_\eps'( K(\x)) \left(\M \left( \Phi_\psi \psi \right) \right)(\x).
    \end{equation*}
Since $b(\x, 2 \pi) = b(\x,0)$ for any $\x,$ we have $\Phi_\psi(\x,2 \pi)= \Phi_\psi(\x,0)$ for any $\x$.
Integration by parts reveals that
    \begin{equation*}
    \begin{split}
    \left(\M \left( \Phi_\psi \psi \right) \right)(\x) = \frac{1}{2 \pi} \int_0^{2 \pi} \Phi_\psi(\x,t)
        \frac{\partial \Phi_\psi}{\partial t}(\x,t) d t \equiv 0.
    \end{split}
    \end{equation*}
Writing $r_\eps := \Phi_\psi \psi \left(k_\eps' \circ K \right)$ we have $ \M r_\eps \equiv 0.$ Apply Lemma
\ref{lemma:f-Mf} with $r_\eps$ in place of $f$ to see that there exists $C>0$ such that
\begin{multline*}
\E^\eps  \Biggl[ \Biggl|
    \int_0^{t \wedge \tau} \Phi_{\psi} \left( \X_u , \frac{u}{\eps^2}\right)
        \psi \left( \X_u , \frac{u}{\eps^2}\right)
         k_\eps'\left( K \left(\X_u \right)\right) d u \Biggr| \Biggr]\\
         = \E^\eps\Biggl[ \Biggl|\int_0^{t \wedge \tau}
         \left(r_\eps \left(\X_u , \frac{u}{\eps^2}\right)
         - (\M r_\eps)(\X_u) \right)  d u \Biggr| \Biggr]\\
        < C \eps^{1-3p} \ (1+t) \ \|g \|_{C^3 \left( [0,K^*]\right)}.
         \end{multline*}
Combine this inequality with (\ref{ineq:part1}) to complete the proof.
\end{proof}

\subsection{Time spent near the boundary of the critical set}

We will use Lemma \ref{lemma:f-Mf} to show that, under $\PP^\eps,$ the process $\X$ does not spend `too much
time' near the boundary $\partV$ of the critical set $\overV.$ This fact will be effectively used later.
\begin{lemma} \label{lemma:residence}
    Given $t  \geq 0,$
    \begin{equation*}
    \overline{\lim}_{\eps \rightarrow 0} \overline{\lim}_{\delta \rightarrow 0}
    \E^\eps \left[ \frac{1}{\delta}
    \int_0^{t \wedge \tau} \chi_{[-\delta,\delta]} \left( K(\X_u) \right) d u
    \right] < \infty.
    \end{equation*}
\end{lemma}

\begin{proof}[\bf{Proof}]
Let $v \in C^\infty(\R;[0, 1])$ be an even function, such that
    \begin{equation*}
    v(h) =
  \begin{cases}
    1 & \text{ if } |h| \leq 1, \\
    0 & \text{ if } |h| \geq 2.
  \end{cases}
 \end{equation*}
 For every $h \in \R$ and $\delta>0,$ we define a function
    \begin{equation*}
    g_\delta(h)
    := \frac{1}{\delta} \int_0^h
    \left(\int_0^r v \left( \frac{s}{\delta} \right) d s \right)d r,
    \end{equation*}
so that
    \begin{equation} \label{ineq:g-delta}
    g''_\delta(h)
    = \frac{1}{\delta} v \left( \frac{h}{\delta} \right)
    \geq  \frac{1}{\delta} \chi_{[-\delta,\delta]} \left( h \right)
    \end{equation} and $g_\delta \in \D^\eps.$
Let $0 \leq t<\tau$ and  $\x \in \overI.$ By the martingale problem  for $
\left(\LL_{t/\eps^2}^\eps,\delta_{\x}, \I \right)\!,$
    \begin{equation*}
    \begin{split}
    g_\delta \left(K \left(\X_t\right) \right) & - g_\delta \left(K \left(\x) \right) \right)
    - \frac{1}{\eps} \int_0^t g_\delta' \left(K \left(\X_u \right) \right)
        \left( \nabla K \left(\X_u \right), b \left( \X_u, \frac{u}{\eps^2}\right)\right) d u\\
    & - \int_{0}^t g_\delta' \left( K \left(\X_u\right)\right)
        \left(\LL_{u/\eps^2} K \right) \left(\X_u \right) d u\\
    & - \frac{1}{2} \int _{0}^t g_\delta'' \left( K \left(\X_u\right)\right)
    \langle d K , d K \rangle_{u/\eps^2} \left(\X_u\right) d u
    \end{split}
    \end{equation*}
is a $\PP^\eps$-martingale. Set
    \begin{equation*}
    \begin{split}
    f_1(\x,t) :& = g_\delta' \left(K\left(\x \right) \right)
        \left( \nabla K \left(\x\right), b \left( \x,t\right)\right),\\
    f_2(\x,t) :& = g_\delta' \left( K \left(\x \right)\right)
        \left(\LL_{t} K \right) \left(\x \right),\\
    f_3(\x, t) :& = g_\delta'' \left( K \left(\x\right)\right)
    \langle d K , d K \rangle_{t} \left(\x\right),
    \end{split}
    \end{equation*}
so that
    \begin{equation*}
    \begin{split}
    \left( \M f_1\right)(\x) & = g_\delta' \left(K\left(\x \right) \right)
        \left( \nabla K \left(\x\right), \left(\M b\right) \left( \x \right)\right)\\
        & = n K^{n-1}(\x) g_\delta' \left(K\left(\x \right) \right)
        \left( \nabla K \left(\x\right), \nabla^{\perp} K \left( \x \right)\right)
        \equiv 0,\\
    \left( \M f_2 \right) (\x) & = g_\delta' \left( K \left(\x \right)\right)
       \left( \M \left(\LL K \right) \right)\left(\x \right),\\
    \left( \M f_3\right)(\x ) & = g_\delta'' \left( K \left(\x\right)\right)
    \left( \M \langle d K , d K \rangle \right) \left(\x\right),
    \end{split}
    \end{equation*}
and
    \begin{equation*}
    \begin{split}
    g_\delta \left(K \left(\X_t\right) \right)
    & - g_\delta \left(K \left(\x) \right) \right)\\
    & - \frac{1}{\eps} \int_0^t \biggl( f_1 \left( \X_u, \frac{u}{\eps^2}\right)
            - \left(\M f_1 \right)\left( \X_u\right)\biggr)d u \\
    & - \int_{0}^t  \biggl( f_2 \left( \X_u, \frac{u}{\eps^2}\right)
            - \left(\M f_2 \right)\left( \X_u\right)\biggr)d u\\
    & - \frac{1}{2} \int _{0}^t \biggl( f_3 \left( \X_u, \frac{u}{\eps^2}\right)
            - \left(\M f_3 \right)\left( \X_u\right)\biggr) d u\\
    & - \int_0^t \left(\M f_2 \right)\left(\X_u\right) d u\\
    & - \frac{1}{2} \int_0^t \left(\M f_3 \right)\left(\X_u\right) d u
    \end{split}
\end{equation*} is a $\PP^\eps$-martingale.
Then for any $\eps>0,\delta>0,$ and for any $0 \leq t \leq \tau,$
    \begin{align}
    \frac{1}{2}\E^\eps \left[ \int_0^t \left(\M f_3 \right)\left(\X_u\right) d u\right]
    & \leq \E^\eps \biggl[\biggl| g_\delta \left(K \left(\X_t\right) \right)
     - g_\delta \left(K \left(\x \right) \right) \biggr|\biggr] \notag\\
     & \quad + \frac{1}{\eps}
     \E^\eps \left[\left|\int_0^t \biggl( f_1 \left( \X_u, \frac{u}{\eps^2}\right)
            - \left(\M f_1 \right)\left( \X_u\right)\biggr)d u\right|\right] \notag\\
     & \quad + \E^\eps \left[\left| \int_{0}^t  \biggl( f_2 \left( \X_u, \frac{u}{\eps^2}\right)
            - \left(\M f_2 \right)\left( \X_u\right)\biggr)d u\right|\right] \notag\\
    & \quad + \frac{1}{2} \E^\eps \left[ \left|
        \int _{0}^t \biggl( f_3 \left( \X_u, \frac{u}{\eps^2}\right)
            - \left(\M f_3 \right)\left( \X_u\right)\biggr) d u\right| \right] \notag\\
    & \quad + \E^\eps
        \left[ \left| \int_0^t \left(\M f_2 \right)\left(\X_u\right) d u\right|\right]. \label{ineq:time}
    \end{align}
We use Remark~\ref{rem:K} to make sure that
    \begin{equation*}
    C_1 := \frac{1}{2} \inf_{\y \in \overline{\EE}}
    \left(\M \langle d K, d K \rangle \right)(\y) > 0,
    \end{equation*}
which, together with inequality (\ref{ineq:g-delta}), guarantees that for all $\x \in \overI$ and for any $0<
\delta<a$
    \begin{equation*}
    \left(\M f_3 \right)(\x)
    \geq C_{1} \frac{1}{\delta} \chi_{[-\delta,\delta]} \left(K \left(\x\right) \right).
    \end{equation*}
Notice that each of the functions $f_i \in C^2\left(\overI \times  [0, \infty); \R \right), i=1,2,3,$ is $2
\pi$-periodic in its last argument, i.e., each of these three functions $f_i$ satisfy all the hypotheses on
the function $f$ in Lemma~\ref{lemma:f-Mf}. Now, inequality~(\ref{ineq:time}) and Lemma~\ref{lemma:f-Mf}
together imply that for all $\x \in \overI,$ all $0 \leq t \leq \tau,$ and for all $0< \delta< a,$
    \begin{multline} \label{ineq:time-1}
    C_1 \E^\eps  \left[ \frac{1}{\delta}\int_0^t \chi_{[-\delta,\delta]}
    \left(K \left(\X_u\right) \right) d u\right]
    \leq  \E^\eps \biggl[\biggl| g_\delta \left(K \left(\X_t\right) \right)
     - g_\delta \left(K \left(\x \right) \right) \biggr|\biggr]\\
    + C \ (1+t)
    \left( \left\| f_1 \right\|_{C^2 \left(\overI \times  [0, \infty) \right)}
        + \eps \left\| f_2 \right\|_{C^2 \left(\overI \times  [0, \infty) \right)}
        + \eps \left\| f_3 \right\|_{C^2 \left(\overI \times  [0, \infty) \right)}\right)\\
    + \E^\eps
        \left[ \left| \int_0^t \left(\M f_2 \right)\left(\X_u\right) d u\right|\right].
    \end{multline}
For all $\delta> 0,$ there exists a constant $C^*>0$ such that
    \begin{equation*}
    \sup_{-K^* \leq h \leq K^*} \left|g'_\delta(h)\right| \leq C^* \quad \text{ and } \quad
    \sup_{-K^* \leq h \leq K^*} \left|g_\delta(h)\right|
    \leq C^*.
    \end{equation*}
Then, for all $\delta > 0,$ all $0 \leq t \leq \tau,$ and $\x \in \overI,$
    \begin{equation} \label{ineq:g-g}
    \E^\eps \biggl[\biggl| g_\delta \left(K \left(\X_t\right) \right)
     - g_\delta \left(K \left(\x \right) \right) \biggr|\biggr] \leq 2 C^*
     \end{equation}
and there exists a positive constant $C^{**}$ such that
    \begin{equation} \label{ineq:Mf2}
    \E^\eps \biggl[\left| \int_0^t \left(\M f_2 \right)(\X_u) d u\right| \biggr]
    <  C^{**} t \left\| K\right\|_{C^2 \left(\overI \right)}.
    \end{equation}
Denote $C_{f_i}= \left\| f_i \right\|_{C^2 \left(\overI \times [0, \infty)\right)}, \ i=1,2,3,$  $C_K =
\left\| K\right\|_{C^2 \left(\overI \right)}.$ Plug inequalities (\ref{ineq:g-g}) and (\ref{ineq:Mf2}) into
(\ref{ineq:time-1}) to obtain that for all $t >0,$
    \begin{equation*}
    \begin{split}
     \overline{\lim}_{\delta \rightarrow 0}
     \E^\eps & \left[ \frac{1}{\delta}\int_0^{t \wedge \tau} \chi_{[-\delta,\delta]}
    \left(K \left(\X_u\right) \right) d u\right]\\
    & \leq \frac{1}{C_1} \left( 2C^*
    + C(1+t)
    \left( C_{f_1} + \eps C_{f_2} + \eps C_{f_3} \right) + C^{**} C_K t \right).
    \end{split}
    \end{equation*}
    Finally, let $\eps$ tend to $0$ to see that the last inequality implies that
    \begin{equation*}
    \begin{split}
     \overline{\lim}_{\eps \rightarrow 0} \overline{\lim}_{\delta \rightarrow 0}
     \E^\eps & \left[ \frac{1}{\delta}\int_0^{t \wedge \tau} \chi_{[-\delta,\delta]}
    \left(K \left(\X_u\right) \right) d u\right]
    \leq \overline{C} + \overline{\overline{C}} t,
    \end{split}
    \end{equation*}
where $\overline{C} = \left(2C^*+C C_{f_1}\right)/C_1 <\infty,$ and $\overline{\overline{C}} = \left(C^{**}
C_K+ C C_{f_1}\right)/C_1 < \infty .$
\end{proof}

\begin{remark}
For the purpose of our proof, the lemma above is not stated in the most general form. The following more
general result holds.

Given $t \geq 0,$ and a function $v \in L_1\left( \R \right)$,
    \begin{equation} \label{eqn:exp-residence}
    \overline{\lim}_{\eps \rightarrow 0}
    \overline{\lim}_{\delta \rightarrow 0}
    \E^\eps \left[ \frac{1}{\delta}
    \int_0^{t \wedge \tau} v
    \left(\frac{ K(\X_u) }{\delta} \right)d u
    \right] < \infty.
    \end{equation}
We can prove this statement analogously, setting $g''_\delta(h)
    = \frac{1}{\delta} v \left( \frac{h}{\delta} \right)$. Then $g'_\delta(h)
    =  \int_0^{h/\delta} v ( s ) ds,$ and $g_\delta(h)
    =  \delta \int_0^{h/\delta}
    \left( \int_0^s v \left( u \right)du \right)
    ds.$

For all $\delta> 0,$ there exists a constant $C^*>0$ such that
    \begin{equation*}
    \sup_{-K^* \leq h \leq K^*} \left|g'_\delta(h)\right|
    \leq \int_0^\infty v \left( s \right)ds
    \leq  C^*
    \end{equation*} and
    \begin{equation*}
    \sup_{-K^* \leq h \leq K^*} \left|g_\delta(h)\right|
    \leq \sup_{-K^* \leq h \leq K^*} \left|\delta \int_0^{h/\delta}
    \left(\int_0^\infty v \left( u \right)du \right) ds
    \right|
    \leq C^*.
    \end{equation*}
Now, to verify inequality (\ref{eqn:exp-residence}), we can repeat all the arguments of the proof of
Lemma~\ref{lemma:residence}.
\end{remark}

\subsection{Stochastic averaging around the boundary of the critical set}
\noindent
 In order to perform analysis close to $\partV,$ and to achieve sufficient accuracy, we will
introduce new coordinates defined on $\EE$. We will use the flow $\zeta_t$ defined on $\EE$ by
    \begin{equation}
    \dot{\zeta}_t(\x) = \frac{\nabla K}{\left\| \nabla K \right\|^2} \left(\zeta_t(\x)
    \right), \,
    \zeta_0(\x)  = \x. \label{eqn:flow-on-E}
    \end{equation} Notice that for any $\x \in \EE$
\begin{equation*}
    \begin{split}
    \frac{d K}{d t} \left(\zeta_t(\x)\right)
    = \left(\nabla K \left(\zeta_t(\x)\right),\dot{\zeta}_t(\x) \right) \equiv 1,
    \end{split}
    \end{equation*}
and if $\x \in \partV$, then $K(\zeta_t(\x)) = t$ for all $-a < t < a.$ Fix a point $\x^* \in
\partV,$ and define the curve
    \begin{equation*}
    \CC := \left\{\zeta_{-s}(\x^*): -a<s<a \right\}.
    \end{equation*}
We reverse time to emphasize that a trajectory of the constructed flow moves a point $\x \in \EE \cup \CC$
towards (not away from) $\x^* \in \partV.$

Notice that
     $
    \EE = \left\{\phi_t(\x): \x \in \CC,
    0 \leq t < \eta \left(K(\x) \right) \right\},
     $
and in particular the map
     $
    (t,s) \mapsto \phi_t \left(\zeta_s \left( \x^* \right)\right)
     $
is a homeomorphism from
    \begin{equation*}
    \left\{(t,s) \in \R^2: -a < s < a, 0 \leq t <   \eta \left(K \left(\zeta_s \left(\x^* \right) \right)\right)
    \right\}
     \end{equation*} to $\EE$. This is the foundation of our new coordinate system for $\EE$. It is
     advantageous to normalize the first coordinate by setting
    \begin{equation*}
  \beta(s) := \int_{\y \in K^{-1}(s)} \frac{\left( \M \langle d K, d K \rangle
    \right)\left( \y \right) }{\left\| \nabla K(\y)\right\|} \,d l(y), \quad \Theta (\x ) := \int_0^t \left( \M \langle d K, d K \rangle
    \right)\left( \phi_s (\x)\right) d s,
    \end{equation*}
where $t$ satisfies $\phi_t \left( \zeta_s \left(\x^*\right)\right)=\x$. The map $\rho: \x \mapsto \left(\Theta(\x),
K(\x) \right)$ is a homeomorphism from the set $\EE$ onto $ \left\{(\theta,s) \in \R^2:  -a<s<a , 0 \leq \theta <
\beta(s)\right\}. $ For future reference, notice that for any $\x \in \EE \setminus \CC$
    \begin{equation} \label{eqn:Theta-Kperp}
    \left( \nabla \Theta (\x), \nabla^\perp K(\x)\right)
    = \left( \M \langle d K, d K \rangle \right)(\x).
    \end{equation}
This equation reflects the essential property of our new coordinates, that equal increments of $\Theta$
correspond to equal amounts of diffusion across $\partV$. Recall the definition (\ref{def:Psi-f}) of
$\Psi_f$, and define
    \begin{equation*}
    \breve{\Psi}_f:= \Psi_f \circ \rho^{-1}.
    \end{equation*}
That is, for any $\x \in \EE, \, \Psi_f(\x) = \breve{\Psi}_f (\theta,s)$ with $ \theta= \Theta(\x)$ and $s =
K(\x).$  We can specify $f$ so that $ \breve{\Psi}_f (\theta,0) \in C^5 \left(\R \right),$ and $
\breve{\Psi}_f \left( \theta + \beta(s) ,s \right) = \breve{\Psi}_f (\theta,s), \theta \in \R, -a < s < a.$
Using property (\ref{eqn:Theta-Kperp}) of the new coordinates, for all $\x \in \EE \setminus \CC,$
    \begin{multline} \label{eqn:new-coord-prop}
    \left( \nabla \Psi_f(\x), \nabla^\perp K(\x) \right)
    = \biggl( \frac{\partial \breve{\Psi}_f}{\partial \theta}
         \nabla \Theta(\x)
        + \frac{\partial \breve{\Psi}_f}{\partial s}
           \nabla K(\x),
         \nabla^\perp K(\x) \biggr)\\
    =  \frac{\partial \breve{\Psi}_f}{\partial \theta}
        \biggl(\nabla \Theta(\x),\nabla^\perp K(\x) \biggr)
    = \frac{\partial \breve{\Psi}_f}{\partial \theta}
    \left( \M \langle d K, d K \rangle \right)(\x),
    \end{multline}
with all the derivatives of $\breve{\Psi}_f$ computed at $\left(\Theta(\x), K(\x) \right).$

As before, we assume that function $\omega \in C_c^\infty(\R)$ is supported on a compact subset of $(0,
\infty).$ We will borrow the following proposition from \cite{S.1}, only with minor changes.

\begin{proposition} \label{prop:PDE}
Let $p = 1/(n+1)$ and
    \begin{equation*}
    \oomega(h) := \omega(h^{2p})\,h^{2p}, \, h \in \R.
    \end{equation*} There exists a function $B \in C^2\left(\R \times [0, \infty)
\right)$ such that
    \begin{align}
    & n \frac{\partial B}{\partial \theta}
    + \frac{(n+1)^2}{8} \frac{\partial^2 B}{\partial s^2}
    + \frac{n^2-1}{8 s} \frac{\partial B}{\partial s}
    =\frac{\partial  \breve{\Psi}_f\left(\theta,0\right)}{\partial \theta}
        \oomega(s), \notag\\
    & B(\theta +\beta(0), s) = B(\theta,s), \notag\\
    & \left| \partial^\alpha B(\theta,s)\right|
    \leq C \left\| \breve{\Psi}_f \left( \cdot,0 \right)\right\|_{C^5 \left(\R\right)}
        e^{-s/C} \label{eqn:PDE-1}
    \end{align}
for all $(\theta, s) \in \R \times (0, \infty),$ for some constant $C>0,$ and for any multi-index $\alpha,
|\alpha| \leq 2$, and such that
    \begin{equation} \label{eqn:initial-cond}
    \begin{split}
    & \frac{\partial B(\theta,0)}{\partial s} = 0
    \end{split}
    \end{equation} for all $\theta \in \R.$
\end{proposition}

It was shown in \cite{S.1} that solving (\ref{eqn:PDE-1}) involves Fourier analysis and leads to an
inhomogeneous Bessel's ordinary differential equation. The unique explicit solution is expressed via Bessel
functions of small order and purely imaginary argument.

Define the functions $\dd_n: \R \rightarrow \R$ by $\dd_n (h) := \biggl(\max \left\{0, h \right\}\biggr)^n.$
This sequence satisfies the recursive relation $\dd_n(h) = h \, \dd_{n-1}(h).$ For some $n>2, \, H(\x) =
\dd_n\left(K(\x) \right)$ and $ \nabla^\perp H(\x) = n \dd_{n-1}\left(K(\x) \right) \nabla^\perp K(\x)$ for
all $\x \in \R^2$.

\begin{proposition} \label{prop:averaging-partV}
There exists a constant $C>0$ such that for any $f \in C^5 \left(\R^2; \R \right)$ and for any $0<p<1$
    \begin{multline} \label{eqn:eps-G-f}
    \left|\E^\eps
    \left[ \int_{s \wedge \tau}^{t \wedge \tau} \biggl(f \left(\X_u \right)
        - \left(\A f \right)\left([\X_u] \right)\biggr)
        \dd_n\left(K \left( \X_u
        \right)\right)
        \omega \left(\frac{K \left( \X_u\right)}{\eps^p} \right)
         d u \prod_{k=1}^{m} h_k \left( \X_{t_k} \right)
    \right] \right| \\
    < C \eps^{1+p} \,(1+ t) \left\| f \right\|_{C^5 (\overI)}.
    \end{multline}
whenever $0 \leq t_1<t_2< \dots<t_m \leq s <t,$ and $h_1, \ldots,h_m \in C_b(\R^2).$
\end{proposition}

\begin{remark}
Since $\omega$ has compact support we can find a constant $C^*>0$ such that for any $\x \in \overI$
    \begin{equation*}
    \left|\dd_n\left(K(\x)\right)\omega \left(\frac{K \left( \x \right)}{\eps^p} \right)\right|
    < C^* \eps^{n p}\chi_{\left[-\eps^p, \eps^p\right]}\left(K(\x)\right).
    \end{equation*}
 Then, by Lemma~\ref{lemma:residence}, there exist constants $C^{**}, C >0,$ such that
    \begin{equation*}
    \begin{split}
    \left|\E^\eps
    \left[ \int_{s \wedge \tau}^{t \wedge \tau} \biggl(f \left(\X_u \right)
    \biggr. \right. \right.& \left. \left. \biggl.
        - \left(\A f \right)\left([\X_u] \right)\biggr)\dd_n\left(K \left( \X_u
        \right)\right)
        \omega \left(\frac{K \left( \X_u\right)}{\eps^p} \right)d u
    \right] \right| \\
    & < C^{**} \eps^{n p}\left\|f \right\|_{C(\overI)}
        \E^\eps \int_{s \wedge \tau}^{t \wedge \tau}
        \chi_{\left[-\eps^p,\eps^p\right]} \left(K \left( \X_u\right)\right)d u\\
    & < C \eps^{p(n+1)} \left\|f \right\|_{C(\overI)}.
        \end{split}
    \end{equation*}
For values of $p< 1/n$, Proposition \ref{prop:averaging-partV} provides a better bound.
\end{remark}

\begin{proof}[\bf{Proof}]
For every $\x \in \overI$ and for all $f \in C^4 \left(\R^2;\R\right)$ we define
    \begin{equation*}
    \GG_f^\eps(\x):= \frac{1}{\eps}
     \biggl(f \left(\x \right)- \left(\A_0 f \right)([\x]) \biggr)
        \dd_n \left(K(\x)\right)\omega \left(\frac{K \left( \x \right)}{\eps^p} \right).
    \end{equation*}

We will perform stochastic averaging using the martingale problem. This approach was developed by
Papanicolaou and Kohler \cite{P.-K.}. We will look for a function $\Psi^\eps \in C^{2} \left(\R^2\right)$
that satisfies
    \begin{equation*}
  \lim_{\eps \rightarrow 0}
    \eps^{-p} \left\|\Psi^\eps \right\|_{C \left(\overI\right)}
    < \infty,
    \end{equation*}
and which, for every $\eps>0,$ $t \geq 0,$ solves the approximate partial differential equation
    \begin{equation} \label{eqn:Leps-approx}
    \left( \M \left(\LL^\eps \Psi^\eps\right)\right)(\x)
    \thickapprox  \GG_f^\eps(\x), \quad \x \in \overI.
    \end{equation}
From the approximate equation (\ref{eqn:Leps-approx}), by Lemma \ref{lemma:f-Mf}, and by the stopped
martingale problem for $\left(\LL^\eps_{t/\eps^2}, \delta_{\x}, \I\right),$
    \begin{multline*}
    \int_0^{t \wedge \tau} \GG_f^\eps \left(\X_u\right)d u
    \thickapprox
    \int_0^{t \wedge \tau} \left(\M \left(\LL^\eps \Psi^\eps\right)\right)
    \left(\X_u\right)d u
    \thickapprox
    \int_0^{t \wedge \tau} \left( \LL^\eps_{u/\eps^2} \Psi^\eps\right)
    \left(\X_u\right)d u \\
    = \Psi^\eps \left(\X_{t \wedge \tau} \right)
    - \Psi^\eps \left(\X_0 \right) + M^\eps_{t \wedge \tau},
    \end{multline*}
where $M^\eps_{t \wedge\tau}$ is a $\PP^\eps$-martingale. Then, we will finish the proof by optional sampling
for this martingale.

The operator
    \begin{equation*}
    \eps \, \M \left(\LL^\eps \Psi^\eps\right)
    = \left( \nabla \Psi^\eps(\x), \nabla^\perp H (\x)\right)
    + \eps \, \M \left(\LL \Psi^\eps\right)
    \end{equation*}
has drift of order $1$ and small diffusion, which leads to a singular perturbation problem.

Notice from (\ref{eqn:f-A0f=grad-grad}) and from property (\ref{eqn:new-coord-prop}) of the new coordinates
that for every $\x \in \EE \setminus \CC$
    \begin{equation} \label{eqn:Gf}
    \GG_f^\eps(\x):= \eps^{-1}
    \frac{\partial \breve{\Psi}_f\left(\Theta(\x), K(\x)\right)}{\partial \theta}
        \dd_n\left(K(\x)\right)\omega \left(\frac{K \left( \x \right)}{\eps^p}
        \right) \
        \left( \M \langle d K, d K \rangle \right)(\x).
    \end{equation}
Recall that $ \oomega(h) = \omega(h^{2p})\,h^{2p}, h \in \R,$ and write $\dd_{n}(h)$ as $\dd_{n-1}(h) h.$
Then, from (\ref{eqn:Gf}),
    \begin{equation*}
    \GG_f^\eps(\x):= \eps^{p-1}
    \frac{\partial \breve{\Psi}_f\left(\Theta(\x), K(\x)\right)}{\partial \theta}
        \dd_{n-1}\left(K(\x)\right)
        \oomega
        \left(\frac{\dd_{1/2p} \left(K \left( \x \right)\right)}{\sqrt{\eps}}
        \right) \
        \left( \M \langle d K, d K \rangle \right)(\x).
    \end{equation*}
For all $\x \in \EE$ we define
    \begin{equation*} \Theta^\circ(\x):= \frac{\beta(0)}{\beta(K(\x))} \,\Theta(\x).
    \end{equation*}
To continue with this proof, we state and prove the following lemma.
\begin{lemma} \label{lemma:Beps}
Let $f \in C^5 \left(\R^2;\R \right)$. Let $B$ satisfy the hypotheses of Proposition~\ref{prop:PDE}, and for
every $\eps>0$ define
\begin{equation*}
B^\eps(\x) := \eps^p \,B \left(\Theta^\circ(\x),
    \frac{\dd_{1/2p}\left(K(\x)\right)}{\sqrt{\eps}}\right)
\end{equation*} for all $\x \in \EE$. Then
\begin{equation*}
B^\eps \in C^1(\EE) \bigcap C^2\left(\EE \setminus (\partV \cup \mathcal{C})\right),
\end{equation*} and
\begin{equation*}
\left( \M \left(\LL^\eps B^\eps\right) \right)(\x):= \lim_{\substack{\y \rightarrow \x\\\y \in \EE \setminus
(\partV \cup \mathcal{C})}} \left( \M \left( \LL^\eps B^\eps\right) \right)(\y)
\end{equation*}
exists for all $\x \in \partV \cup \mathcal{C}.$ There exist positive constants $C_1,C_2>0$ such that
\begin{equation*}
\left|\left( \M \left(\LL^\eps B^\eps \right) \right)(\x) - \GG_f^\eps(\x) \right| \leq C_1 \left\|f
\right\|_{C^5(\overI)}\exp \left\{-\frac{|K(\x)|^{1/2p}}{C_2 \sqrt{\eps}}\right\},
\end{equation*}
\begin{equation*}
\sqrt{ \left(\M \left\langle d B^\eps, d B^\eps \right\rangle \right)(\x)} \leq C_1 \left\|f \right\|_{C^5(\overI)}\exp
\left\{-\frac{|K(\x)|^{1/2p}}{C_2 \sqrt{\eps}}\right\},
\end{equation*} and
\begin{equation*}
\left|B^\eps (\x)\right| \leq C_1 \eps^p \left\|f \right\|_{C^5(\overI)}
\end{equation*} for all $\x \in \EE \setminus
(\partV \cup \mathcal{C})$ and  $\eps>0.$
\end{lemma}

\begin{proof}[\bf{Proof}]
First of all, recall that the small set $\EE = \left\{\y \in \R^2: \left|K(\y) \right|< a \right\} \supset
\partV $ is homeomorphic to the set $\left\{ (\theta,s) \in \R^2: -a<s<a , 0 \leq \theta < \beta(s) \right\},$
and that $\mathcal{C}$ is a fixed curve in $\EE,$ which crosses all $K^{-1}(s)$ transversally.

Notice that the definition of $B^\eps$ can be rewritten more explicitly as
\begin{equation*}
B^\eps(\x): = \begin{cases}
    \eps^p B \left( \,\Theta^\circ(\x), 0 \,\right)
        & \text{ if }  -a < K(\x) < 0,\\
    \eps^p B \left( \Theta^\circ(\x),
    \sqrt{\frac{K^{n+1}(\x)}{\eps}} \,\right)
        & \text{ if }\, 0 \leq K(\x)<a.
    \end{cases}
\end{equation*}

We know that $B \in C^2 \left(\R \times [0,\infty)\right)$ is periodic in its first argument, which, together
with the properties of $\Theta^\circ$ and $\dd_{1/2p} \circ K$ away from $\partV \cup \mathcal{C},$
guarantees that $B^\eps \in C^2 \left(\EE \setminus \left(
\partV \cup \mathcal{C}\right)\right).$

The function $K$ is smooth inside $\EE,$ and $1/2p = (n-1)/2>0$ implies differentiability of $\dd_{1/2p}$ on
$\left(-a,a\right).$ Therefore, we are guaranteed that $\dd_{1/2p} \circ K \in C^1 (\EE)$. By the definition
of $\Theta^\circ,$ it is guaranteed to be differentiable on $\EE.$ Notice that given any $\x \in \EE,$ we
have $\text{span}\left\{\nabla K(\x), \nabla^\perp K(\x)\right\} = \R^2$. From the periodicity in the first
argument and from the boundary conditions on $B$, we also see that for every $\x \in \EE$ both $\lim_{\y
\rightarrow \x} \left(\nabla K, \nabla B^\eps \right)(\y) $ and $\lim_{\y \rightarrow \x}\left( \nabla^\perp
K, \nabla B^\eps\right)(\y)$ exist. Therefore, $B^\eps \in C^1\left(\EE\right).$ Thus, $B^\eps \in C^1
\left(\EE\right) \bigcap C^2 \left(\EE \setminus \left( \partV \cup \mathcal{C}\right) \right).$

Let $\x \in \EE \setminus \left( \partV \cup \mathcal{C} \right),$ and consider
    \begin{align}
    \LL^\eps_t B \left(\Theta^\circ(\x),
    \frac{\dd_{1/2p}\left(K(\x)\right)}{\sqrt{\eps}}\right)
    & = \frac{1}{2} \frac{\partial^2 B}{\partial \theta^2}
        \langle d \Theta^\circ d \Theta^\circ \rangle_t(\x)
        + \frac{\partial B}{\partial \theta}
        \left(\LL_t
        \Theta^\circ \right)(\x) \notag\\
    & \quad + \frac{1}{\sqrt{\eps}} \dd_{1/2p}'\left( K(\x)\right)
        \frac{\partial^2 B}{\partial s \partial \theta}
        \langle d K,d \Theta^\circ \rangle_t (\x) \notag\\
    & \quad + \frac{1}{2\eps} \frac{\partial^2 B}{\partial s^2}
        \left(\dd_{1/2p}'\left( K(\x)\right)\right)^2
        \langle d K,d K \rangle_t (\x) \notag\\
    & \quad + \frac{1}{2\sqrt{\eps}} \frac{\partial B}{\partial s}
        \dd_{1/2p}''\left( K(\x)\right)
        \langle d K,d K \rangle_t (\x) \notag\\
    & \quad + \frac{1}{\sqrt{\eps}} \frac{\partial B}{\partial s}
        \dd_{1/2p}'\left( K(\x)\right)
        \left(\LL_t K\right) (\x), \label{eqn:Lepst-Psi}
    \end{align}
    and
    \begin{align}
    \left\langle dB, dB \right\rangle_t(\x)
    & = \left(\frac{\partial B}{\partial \theta}\right)^2
        \left\langle d \Theta^\circ, d \Theta^\circ
        \right\rangle_t(\x) \notag\\
    & \quad + \frac{1}{\eps}
         \left(\frac{\partial B}{\partial s}\right)^2
        \left(\dd'_{1/2p}\left(K(\x)\right)\right)^2
        \left\langle d K, d K \right\rangle_t(\x) \notag\\
    & \quad + \frac{2}{\sqrt{\eps}} \frac{\partial B}{\partial \theta}
        \frac{\partial B}{\partial s}
        \dd'_{1/2p}\left(K(\x)\right)
        \left\langle d \Theta^\circ, d K \right\rangle_t(\x). \label{eqn:dB-dB}
    \end{align}
Here and below, $B$ and its partial derivatives are evaluated at $\left(\Theta^\circ(\x),
\frac{\dd_{1/2p}\left(K(\x)\right)}{\sqrt{\eps}} \right),$ unless specified otherwise. From
(\ref{eqn:Lepst-Psi}) and (\ref{eqn:dB-dB}), using $$\dd_{1/2p}'(h) = \frac{1}{2p} \dd_{1/2p-1}(h)$$ and
$$\dd_{1/2p}''(h) = \frac{1}{2p} \left(\frac{1}{2p}-1 \right)\dd_{1/2p-2}(h),$$ we can see that
    \begin{align}
    \left(\LL^\eps_t B^{\eps}\right)(\x)
    & = \frac{1}{\eps} \left(\nabla B^\eps(\x),
        b\left(\x,t\right) \right)
        + \left(\LL_t B^{\eps}\right)(\x) = \eps^{p-1} \frac{\partial B}{\partial \theta}
        \left( \nabla \Theta^\circ(\x), b \left(\x,t\right)\right) \notag\\
    & \quad + \eps^{p-3/2} \frac{\partial B}{\partial s}m
        \dd_{1/2p-1}\left(K(\x)\right) \left( \nabla K(\x),
        b\left(\x,t\right)\right) \notag\\
    & \quad + \eps^p \left( \frac{1}{2} \frac{\partial^2 B}{\partial \theta^2}
        \langle d \Theta^\circ d \Theta^\circ \rangle_t(\x)
        + \frac{\partial B}{\partial \theta}
        \left(\LL_t
        \Theta^\circ \right)(\x) \right. \notag\\
    & \quad + \frac{1}{\sqrt{\eps}} \frac{1}{2p} \dd_{1/2p-1}\left( K(\x)\right)
        \frac{\partial^2 B}{\partial s \partial \theta}
        \langle d K,d \Theta^\circ \rangle_t(\x) \notag\\
    & \quad + \frac{1}{2\eps} \frac{\partial^2 B}{\partial s^2}
       \left( \frac{1}{2p}\right)^2 \dd_{1/2p-1}^2\left( K(\x)\right)
        \langle d K,d K \rangle_t (\x) \notag\\
    & \quad + \frac{1}{2\sqrt{\eps}} \frac{\partial B}{\partial s}
        \frac{1}{2p}\left(\frac{1}{2p}-1\right)\dd_{1/2p-2}\left( K(\x)\right)
        \langle d K,d K \rangle_t (\x) \notag\\
    & \quad+ \left.\frac{1}{\sqrt{\eps}} \frac{\partial B}{\partial s}
        \frac{1}{2p} \dd_{1/2p-1}\left( K(\x)\right)
        \left(\LL_t K\right) (\x)\right) \label{eqn:Lepst-Psieps}
    \end{align}
    and
    \begin{align}
    \left\langle dB^\eps, dB^\eps \right\rangle_t(\x)
    & = \eps^{2p}\left(\frac{\partial B}{\partial \theta}\right)^2
        \left\langle d \Theta^\circ, d \Theta^\circ
        \right\rangle_t(\x) \notag\\
    & \quad + \eps^{2p-1} \left(\frac{1}{2p}\right)^2
        \left(\frac{\partial B}{\partial s}\right)^2
        \dd^2_{1/2p-1}\left(K(\x)\right)
        \left\langle d K, d K \right\rangle_t(\x) \notag\\
    & \quad + 2 \eps^{2p-1/2} \frac{\partial B}{\partial \theta}
        \frac{\partial B}{\partial s}
        \frac{1}{2p}\dd_{1/2p-1}\left(K(\x)\right)
        \left\langle d K, d \Theta^\circ \right\rangle_t(\x). \label{eqn:dBeps-dBeps}
    \end{align}
Recall that $\left(\M \left(\nabla K, b\right)\right)(\x)= 0.$ Then, from (\ref{eqn:Lepst-Psieps}),
    \begin{equation*}
    \begin{split}
    \left( \M \left(\LL^\eps B^\eps \right)\right)(\x)
    & = \frac{1-2p}{8p^2} \eps^{p-1/2} \dd_{1/2p-2}\left(K(\x)\right)
        \frac{\partial B}{\partial s}
        \left(\M \langle dK,dK\rangle\right)(\x)\\
    & \quad + \frac{1}{8 p^2} \eps^{p-1} \dd^2_{1/2p-1}\left(K(\x)\right)
        \frac{\partial^2 B}{\partial s^2}
        \left(\M \langle dK,dK\rangle\right)(\x)\\
    & \quad + \frac{1}{2p} \eps^{p-1/2} \dd_{1/2p-1}\left(K(\x)\right)
        \frac{\partial B}{\partial s}
        \left(\M \left(\LL K\right)\right)(\x)\\
    & \quad + \frac{1}{2p} \eps^{p-1/2} \dd_{1/2p-1}\left(K(\x)\right)
        \frac{\partial^2 B}{\partial s \partial \theta }
        \left(\M \langle d K,d \Theta^\circ \rangle \right)(\x)\\
    & \quad + \eps^p \left[ \frac{\partial B}{\partial \theta}
        \left(\M \left(\LL \Theta^\circ\right)\right)(\x)
        + \frac{1}{2} \frac{\partial^2 B}{\partial \theta^2 }
        \left(\M \langle d \Theta^\circ,d \Theta^\circ \rangle\right)
        (\x)\right]\\
    & \quad + n \eps^{p-1} \dd_{n-1}\left(K(\x)\right)
        \frac{\partial B}{\partial \theta}
        \left( \nabla \Theta^\circ(\x), \nabla^\perp K(\x)\right),
    \end{split}
    \end{equation*}
and from (\ref{eqn:dBeps-dBeps}),
    \begin{equation*}
    \begin{split}
    \left( \M \left\langle dB^\eps, dB^\eps \right\rangle \right)(\x)
    & = \eps^{2p}\left(\frac{\partial B}{\partial \theta}\right)^2
        \left( \M \left\langle d \Theta^\circ, d \Theta^\circ
        \right\rangle \right)(\x)\\
    & \quad + \frac{1}{4p^2}\eps^{2p-1}
        \left(\frac{\partial B}{\partial s}\right)^2
        \dd^2_{1/2p-1}\left(K(\x)\right)
        \left( \M \left\langle d K, d K \right\rangle \right)(\x)\\
    & \quad +  \frac{1}{p}\eps^{2p-1/2}
        \frac{\partial B}{\partial \theta}
        \frac{\partial B}{\partial s}
        \dd_{1/2p-1}\left(K(\x)\right)
        \left( \M  \left\langle d K, d \Theta^\circ \right\rangle \right)(\x).
    \end{split}
    \end{equation*}
From (\ref{eqn:Theta-Kperp}),
    \begin{equation*}
    \left( \nabla \Theta^\circ(\x),\nabla^\perp K(\x)\right)
    = \frac{\beta(0)}{\beta\left(K(\x)\right)}
    \left(\M \langle dK,dK\rangle \right)(\x).
    \end{equation*}
Take $p=\frac{1}{n+1}$. Then, $\frac{1-2p}{8p^2} = \frac{n^2-1}{8},$ $\frac{1}{8p^2} = \frac{(n+1)^2}{8},$
and
    \begin{align}
    \left( \M \left(\LL^\eps B^\eps\right) \right)(\x)
    & = \eps^{p-1} \dd_{n-1}\left(K(\x)\right)
        \left(\M \langle dK,dK\rangle \right)(\x) \notag\\
    & \quad \times
        \left(n \frac{\partial B}{\partial \theta}
        + \frac{(n+1)^2}{8} \frac{\partial^2 B}{\partial s^2}
        + \frac{n^2-1}{8}
        \left( \frac{\dd_{1/2p}\left(K(\x)\right)}{\sqrt{\eps} }\right)^{-1}
        \frac{\partial B}{\partial s}
        \right) \notag \\
    & \quad + I_1^\eps(\x)+I_2^\eps(\x) +I_3^\eps(\x)
    = \GG_f^\eps(\x) + \sum_{i=1}^4 I^\eps_i(\x), \label{eqn:M-L-Beps}
         \end{align}
    where
    \begin{equation*}
    I_1^\eps(\x) = \frac{1}{2p} \eps^{p-1/2} \dd_{1/2p-1}\left(K(\x)\right)
        \left[ \frac{\partial B}{\partial s}
        \left(\M \left(\LL K \right)\right)(\x)
        + \frac{\partial^2 B}{\partial s \partial \theta }
        \left(\M \langle d K,d \Theta^\circ \rangle \right)
        (\x)\right],
    \end{equation*}
     \begin{equation*}
    I_2^\eps(\x) = \eps^p \left[ \frac{\partial B}{\partial \theta}
        \left(\M \left(\LL \Theta^\circ\right) \right)(\x)
        + \frac{1}{2} \frac{\partial^2 B}{\partial \theta^2 }
        \left(\M \langle d \Theta^\circ,d \Theta^\circ \rangle \right)
        (\x)\right],
    \end{equation*}
    \begin{equation*}
    I_3^\eps(\x) = n \eps^{p-1} \dd_{n-1}\left(K(\x)\right)
        \left(\M \langle d K,d K \rangle \right)(\x)
       \frac{\beta(0)-\beta
       \left(K(\x)\right)}{\beta\left(K(\x)\right)}
       \frac{\partial B}{\partial \theta},
    \end{equation*}
    \begin{equation*}
    I_4^\eps(\x) = \eps^{p-1} \dd_{n-1}\left(K(\x)\right)
    \left(\frac{\partial \breve{\Psi}_f}{\partial \theta}
        \left(\Theta(\x),K(\x)\right)
    - \frac{\partial \breve{\Psi}_f}{\partial \theta}(\Theta^\circ,0)\right)
        \left(\M \langle dK,dK\rangle \right)(\x) .
    \end{equation*}
    From this we see that
    \begin{equation*}
    \lim_{\substack{\y \rightarrow \x\\\y \in \EE \setminus (\partV
    \cup \mathcal{C})}} \left( \M \left( \LL^\eps B^\eps\right)
    \right)(\y)
    \end{equation*}
exists for all $\x \in \partV \cup \mathcal{C}$. Also,
    $
    \left( \M \left\langle dB^\eps, dB^\eps \right\rangle \right)(\x)
     = \sum_{i=5}^7 I^\eps_i(\x),
    $ where
    \begin{equation*}
    I_5^\eps(\x)
        = \eps^{2p}\left(\frac{\partial B}{\partial \theta}\right)^2
        \left( \M \left\langle d \Theta^\circ, d \Theta^\circ
        \right\rangle \right)(\x),
    \end{equation*}
    \begin{equation*}
    I_6^\eps(\x) = \frac{1}{4p^2}\eps^{2p-1}
        \left(\frac{\partial B}{\partial s}\right)^2
        \dd_{n-1}\left(K(\x)\right)
        \left( \M \left\langle d K, d K \right\rangle \right)(\x),
    \end{equation*} and
    \begin{equation*}
    I_7^\eps(\x) = \frac{1}{p}\eps^{2p-1/2}
        \frac{\partial B}{\partial \theta}
        \frac{\partial B}{\partial s}
        \dd_{1/2p-1}\left(K(\x)\right)
        \left( \M  \left\langle d K, d \Theta^\circ \right\rangle
        \right)(\x).
    \end{equation*}

    From Proposition \ref{prop:PDE}, we know that
$B \in C^2\left(\R \times [0, \infty) \right)$ and
    \begin{equation} \label{ineq:B(theta,s)}
    \left| \partial^\alpha B(\theta,s)\right|
    \leq C \left\| \breve{\Psi}_f (\cdot,0) \right\|_{C^5 \left(\R  \right)}
        e^{-s/C}
    \end{equation}
for all $(\theta, s) \in \R \times (0, \infty),$ for some constant $C>0,$ and for any multi-index $\alpha,
|\alpha| \leq 2.$ Also, there exists a positive constant which does not depend of $f$ such that
    \begin{equation} \label{ineq:Psi-breve}
    \left\| \breve{\Psi}_f (\cdot,0) \right\|_{C^5 \left(\R \right)}
    \leq C \left\| f \right\|_{C^5 \left(\overI \right)}.
    \end{equation}
Combining inequalities (\ref{ineq:B(theta,s)}) and (\ref{ineq:Psi-breve}), we see that there exist constants
$C_1,C_2>0$ such that
\begin{equation} \label{ineq:B(theta,s)+Psi-breve}
    \left| \partial^\alpha B(\theta,s)\right|
    \leq C_1 \left\| f \right\|_{C^5 \left(\overI \right)}
        e^{-s/C_2}
    \end{equation}
for all $(\theta, s) \in \R \times (0, \infty),$ and for any multi-index $\alpha, |\alpha| \leq 2.$ Now, for
every $\x \in \EE \setminus \left(\partV \cup \mathcal{C} \right),$ there exist positive constants $C_1,C_2$
such that
    \begin{equation*}
    \left| B^\eps(\x)\right|
    \leq \eps^p \left| \partial^0 B(\theta,0)\right|
    \leq C_1 \eps^p \left\| f \right\|_{C^5 \left(\overI \right)},
    \end{equation*}
and positive constants $C_{ij}$ such that
    \begin{equation*}
    \begin{split}
    \left|I_1^\eps(\x) \right|
    & \leq C_{11} \eps^{p-1/2} \eps^{p \left((1/2p)-1 \right)}
        \left\| f \right\|_{C^5 \left(\overI \right)}
        \exp \left\{-\frac{|K(\x)|^{1/2p}}{C \sqrt{\eps}}\right\}\\
        & \quad = C_{12} \left\| f \right\|_{C^5 \left(\overI \right)}
        \exp \left\{-\frac{|K(\x)|^{1/2p}}{C \sqrt{\eps}}\right\},
        \end{split}
    \end{equation*}
    \begin{equation*}
    \begin{split}
    \left|I_2^\eps(\x) \right|
    & \leq C_{21} \eps^{p}
        \left\| f \right\|_{C^5 \left(\overI \right)}
        \exp \left\{-\frac{|K(\x)|^{1/2p}}{C \sqrt{\eps}}\right\},
        \end{split}
    \end{equation*}
    \begin{equation*}
    \begin{split}
    \left|I_3^\eps(\x) \right|
    & \leq C_{31} n \eps^{p-1}
        \left|\frac{\beta'(K(\x))}{\beta(K(\x))} \right|
        \dd_{n} \left(K(\x)\right)
        \left\| f \right\|_{C^5 \left(\overI \right)}
        \exp \left\{-\frac{|K(\x)|^{1/2p}}{C \sqrt{\eps}}\right\}\\
    & \quad \leq C_{32} n  \eps^{p-1+pn}
        \left\| f \right\|_{C^5 \left(\overI \right)}
        \exp \left\{-\frac{|K(\x)|^{1/2p}}{C \sqrt{\eps}}\right\}\\
    & \quad \leq C_{33}
        \left\| f \right\|_{C^5 \left(\overI \right)}
        \exp \left\{-\frac{|K(\x)|^{1/2p}}{C \sqrt{\eps}}\right\},
        \end{split}
    \end{equation*}
    \begin{equation*}
    \begin{split}
    \left|I_4^\eps(\x) \right|
    & \leq C_{41} \eps^{p}
        \left\| f \right\|_{C^5 \left(\overI \right)}
        \exp \left\{-\frac{|K(\x)|^{1/2p}}{C \sqrt{\eps}}\right\},
        \end{split}
    \end{equation*}
    \begin{equation*}
    \begin{split}
    \left|I_5^\eps(\x) \right|
    & \leq C_{51}^2 \eps^{2p}
        \left\| f \right\|^2_{C^5 \left(\overI \right)}
        \exp \left\{-\frac{2|K(\x)|^{1/2p}}{C \sqrt{\eps}}\right\},
        \end{split}
    \end{equation*}
    \begin{equation*}
    \begin{split}
    \left|I_6^\eps(\x) \right|
    & \leq C_{61}^2 \eps^{2p-1+p(n-1)}
        \left\| f \right\|^2_{C^5 \left(\overI \right)}
        \exp \left\{-\frac{2|K(\x)|^{1/2p}}{C \sqrt{\eps}}\right\}\\
        & \quad = C_{62}^2
        \left\| f \right\|^2_{C^5 \left(\overI \right)}
        \exp \left\{-\frac{2|K(\x)|^{1/2p}}{C \sqrt{\eps}}\right\},
        \end{split}
    \end{equation*} and
    \begin{equation*}
    \begin{split}
    \left|I_7^\eps(\x) \right|
    & \leq C_{71}^2 \eps^{p}
        \left\| f \right\|^2_{C^5 \left(\overI \right)}
        \exp \left\{-\frac{2|K(\x)|^{1/2p}}{C \sqrt{\eps}}\right\}.
        \end{split}
    \end{equation*}
Then, there exist positive constants $C_1$ and $C_2$ such that for every $\x \in \EE \setminus \left(\partV
\cup \mathcal{C}\right)$
    \begin{equation*}
    \begin{split}
    \left|\left( \M \left(\LL^\eps B^\eps \right) \right)(\x) -
        \GG_f^\eps(\x) \right|
    \leq \sum_{i=1}^5 \left|I_i^\eps(\x) \right|
     \leq C_1 \left\|f \right\|_{C^5 \left(\overI \right)}
        \exp \left\{-\frac{|K(\x)|^{1/2p}}{C_2 \sqrt{\eps}}\right\}
    \end{split}
    \end{equation*} and
    \begin{equation*}
    \begin{split}
    \sqrt{ \left( \M \left\langle d B^\eps,
        d B^\eps \right\rangle \right)(\x)}
     \leq \sum_{i=5}^7 \sqrt{ \left| I_i^\eps(\x) \right|}
     \leq C_1 \left\|f \right\|_{C^5 \left(\overI \right)}
        \exp \left\{-\frac{|K(\x)|^{1/2p}}{C_2 \sqrt{\eps}}\right\}.
    \end{split}
    \end{equation*} This completes our proof of Lemma~\ref{lemma:Beps}.
    \end{proof}

Examining expressions (\ref{eqn:Leps-approx}) and (\ref{eqn:M-L-Beps}), we see that we would like to apply
the martingale problem to a smoothed version of $B^\eps$ which we call $\Psi^\eps$. From Proposition
\ref{lemma:Beps}, $B^\eps$ is defined on $\EE$ only, and its second derivatives are not guaranteed to be
continuous on $\partV \cup \mathcal{C} \subset \EE$. On the other hand, we would like to have $\Psi^\eps \in
C^2 \left(\R^2\right)$.

First, we define $v_a \in C_c^\infty \left(\R;[0,1] \right)$ by $ v_a(h) :=v \left( 3 h/a\right),$ then
    \begin{equation*}
    v_a(h) = \begin{cases} 1 & \text{ if } |h| \leq a/3,\\0 & \text{ if } |h| \geq 2a/3.
    \end{cases}
    \end{equation*}
Using this function we cut off $B^\eps$ near the edges of $\EE$ so that
    \begin{equation*}
    \Psi^\eps(\x)
    = v_a\left( K(\x) \right) B^\eps(\x)
    \end{equation*}
for every $\x \in \R^2$. Now, notice that
    \begin{multline*}
    \left( \M \left(\LL^\eps \Psi^\eps\right)\right)(\x)
        - \GG_f^\eps(\x)
    = v_a \left(K(\x)\right)
    \left( \left( \M \left(\LL^\eps B^\eps\right)\right)(\x)
        - \GG_f^\eps(\x)\right)\\
    + B^\eps(\x) \left(v'_a\left(K(\x) \right)
        \left(\M \left(\LL K\right)\right)(\x)
        + \frac{1}{2} v_a''\left(K(\x) \right)
        \left(\M \langle d K, dK \rangle \right)(\x)\right)\\
    + v'_a \left( K(\x)\right)
            \left(\M \langle d B^\eps, dK \rangle \right)(\x)
            \end{multline*}
for all $\x \in \R^2 \setminus \left( \partV \cup \mathcal{C}\right).$ Using Lemma~\ref{lemma:Beps}, there
exists a positive constant $C$ such that
\begin{multline*}
    \left|\left(
        \M \left(\LL^\eps \Psi^\eps\right)
        \right)(\x) - \GG_f^\eps(\x)
    \right|
    \leq
    \left|\left( \left( \M \left(\LL^\eps B^\eps\right)\right)(\x)
        - \GG_f^\eps(\x)\right)\right|\\
        + \left|B^\eps(\x)\right| \cdot \left|\beta'\left(K(\x) \right)
        \left(\M \left(\LL K\right)\right)(\x)
        + \frac{1}{2} \beta''\left(K(\x) \right)
        \left(\M \langle d K, dK \rangle \right)(\x)\right|\\
    + \left|\beta' \left( K(\x)\right)\right| \cdot
            \left|\left(\M \langle d B^\eps, dK \rangle
            \right)(\x)\right|\\
     \leq C \left\|f \right\|_{C^5 \left(\overI \right)}
        \left( \eps^p + \exp \left\{-\frac{|K(\x)|^{1/2p}}{C_2 \sqrt{\eps}}\right\}
        \right).
        \end{multline*}

Notice that $\Psi^\eps \in C^1\left( \R^2 \right) \cup C^2 \left( \R^2 \setminus \left( \partV \cup
\mathcal{C} \right)\right).$ Also, $$\lim_{\substack{\y \rightarrow \x\\\y \in \EE \setminus (\partV \cup
\mathcal{C})}} \left( \LL^\eps_t \Psi^\eps\right)(\y) $$ exists for all $t \geq 0$ and for all $ \x \in
\partV \cup \mathcal{C} $ since
    \begin{equation*}
        \lim_{\substack{\y \rightarrow \x\\\y \in \EE \setminus (\partV \cup
    \mathcal{C})}} \left( \M \left( \LL^\eps \Psi^\eps\right) \right)(\y)
    \end{equation*}
exists for all $\x \in \partV \cup \mathcal{C} $.

We will use an even function $u_2 \in C_c^\infty \left(\R^2 \right)$ with $\int_{\R^2}u_2(\x) d\x  = 1.$ For
$\delta>0,$ we set the mollifier
    \begin{equation*}
    \Psi^\eps_\delta (\y)
    := \delta^{-2} \int_{\R^2}
        u_2\left( \frac{\y -\z}{\delta }\right)
        \Psi^\eps(\z) d \z, \ \y \in \R^2,
    \end{equation*}
and observe that $\Psi^\eps_\delta \in C^\infty \left(\R^2\right) \subset \D^\eps.$ Now, if $0 \leq t_1<t_2<
\dots<t_m \leq s < t,$ and $h_1, \ldots,h_m \in C_b(\R^2),$ then
\begin{equation*}
    \begin{split}
    \E^\eps &
        \Biggl[
        \Biggl\{
            \left(\Psi^\eps_\delta  \right)
                \left( \X_{t \wedge \tau }\right)
            - \left(\Psi^\eps_\delta  \right)
            \left( \X_{s \wedge \tau}\right)
        \Biggr.\Biggr.\Biggl.\Biggl.
        - \int_{s \wedge \tau}^{t \wedge \tau}
              \left( \LL^\eps_{u/\eps^2} \Psi^\eps_\delta
             \right)
                \left( \X_u \right) d u
        \Biggr\}
            \prod_{k=1}^m h_k \left( \X_{t_k} \right)
        \Biggr]=0.
     \end{split}
\end{equation*}
Using properties of mollifiers, $ \Psi^\eps_\delta \rightarrow \Psi^\eps$ a.e. as $\delta \rightarrow 0$, and
    \begin{equation*}
    \LL^\eps_t \Psi^\eps_\delta (\y) = \delta^{-2}\int_{\R^2}
        u_2\left( \frac{\y -\z}{\delta }\right)
        \LL^\eps_t \Psi^\eps(\z)
        d \z.
        \end{equation*}
Changing the order of integration,
\begin{align}
    \E^\eps &
        \Biggl[
        \Biggl\{
            \Psi^\eps \left( \X_{t \wedge \tau }\right)
            - \Psi^\eps \left( \X_{s \wedge \tau}\right)
        \Biggr.\Biggr.\Biggl.\Biggl. -
        \int_{s \wedge \tau}^{t \wedge \tau}
              \left( \LL^\eps_{u/\eps^2} \Psi^\eps
             \right)
                \left( \X_u \right) d u
        \Biggr\}
            \prod_{k=1}^m h_k \left( \X_{t_k} \right)
        \Biggr] \notag\\
        & =
    \lim_{\delta \rightarrow 0 }\E^\eps
        \Biggl[
        \Biggl\{
            \Psi^\eps_\delta \left( \X_{t \wedge \tau }\right)
            - \Psi^\eps_\delta \left( \X_{s \wedge \tau}\right)
        \Biggr.\Biggr.\Biggl.\Biggl.
        - \int_{s \wedge \tau}^{t \wedge \tau}
              \left( \LL^\eps_{u/\eps^2} \Psi^\eps_\delta
             \right)
                \left( \X_u \right) d u
        \Biggr\}
            \prod_{k=1}^m h_k \left( \X_{t_k} \right)
        \Biggr] \notag\\
    & = 0. \label{eqn:mtg-prob-Psi-eps}
     \end{align}

Applying Lemma~\ref{lemma:f-Mf}, there exist positive constants $C_1$ and $C_2$ such that
\begin{multline*}
    \left|\E^\eps
        \Biggl[ \int_{s \wedge \tau}^{t \wedge \tau}
              \left\{ \left( \LL^\eps_{u/\eps^2} \Psi^\eps
             \right)
                \left( \X_u \right)
                - \left(\M \left(\LL^\eps \Psi^\eps\right)\right)
                \left( \X_u \right)
        \right\} du
            \prod_{k=1}^m h_k \left( \X_{t_k} \right)
        \Biggr]\right|\\
        < C_1\, \eps \,(1+t) \left\|\Psi^\eps \right\|_{C^4
        \left(\overI\right)}< C_2 \,\eps^{p+2} \left\|f \right\|_{C^5
        \left(\overI\right)}.
     \end{multline*}
Next, using the ideas described at the beginning of the proof of Proposition \ref{prop:averaging-partV}, and
Eq.~(\ref{eqn:mtg-prob-Psi-eps}),
    \begin{equation*}
    \begin{split}
    & \biggl|\E^\eps
    \biggl[ \int_{s \wedge \tau}^{t \wedge \tau}  \biggl(f \left(\X_u \right)
        - \left(\A_0 f \right)\left([\X_u] \right)\biggr)
            \dd_n\left(K \left( \X_u \right)\right)
        \omega \left(\frac{K \left( \X_u\right)}{\eps^p} \right)
         d u \prod_{k=1}^m h_k \left( \X_{t_k} \right)
    \biggr] \biggr| \\
    & \quad  = \left|\E^\eps
    \left[ \int_{s \wedge \tau}^{t \wedge \tau}
    \eps \GG_f^\eps \left(\X_u \right) d u
    \prod_{k=1}^m h_k \left( \X_{t_k} \right)
    \right] \right|\\
    & \quad \leq
     \eps \left|\E^\eps
    \left[ \int_{s \wedge \tau}^{t \wedge \tau}
    \biggl( \GG_f^\eps \left(\X_u \right)
    - \left(\M \left(\LL^\eps \Psi^\eps\right) \right)\left(\X_u \right) \biggr) d u
    \prod_{k=1}^m h_k \left( \X_{t_k} \right)
    \right] \right|\\
     & \quad + \eps \left|\E^\eps
    \left[ \int_{s \wedge \tau}^{t \wedge \tau}
   \biggl(\left(\M \left(\LL^\eps \Psi^\eps\right)\right) \left(\X_u\right)
    - \left(\LL^\eps_{u/\eps^2} \Psi^\eps \right)\left(\X_u\right)
    \biggr)d u
    \prod_{k=1}^m h_k \left( \X_{t_k} \right)
    \right] \right|\\
    & \quad + \eps \left|\E^\eps \left[
        \biggl(
            \Psi^\eps \left( \X_{t \wedge \tau }\right)
            - \Psi^\eps \left( \X_{s \wedge \tau}\right)
        \biggr)
            \prod_{k=1}^m h_k \left( \X_{t_k} \right) \right]
            \right|.
    \end{split}
    \end{equation*}
Therefore, there exists a positive constant $C$ such that
    \begin{multline*}
    \left|\E^\eps
    \left[ \int_{s \wedge \tau}^{t \wedge \tau}
        \!\biggl(f \left(\X_u \right)
        - \left(\A_0 f \right)\left([\X_u] \right)
        \biggr)
        \dd_n\left(K \left( \X_u \right)\right)
        \omega \left(\frac{K \left( \X_u\right)}{\eps^p} \right)
         d u \prod_{k=1}^m h_k \left( \X_{t_k} \right)
    \right] \right| \\
    \quad \leq C \eps \left\| f \right\|_{C^5 \left(\overI \right)}
        \left( \eps^p + \eps^{p+2}+ \E^\eps
        \left[ \int_{s \wedge \tau}^{t \wedge \tau}
        \left(\eps^p
            + \exp \left\{-\frac{|K \left(\X_u \right)|^{1/2p}}{C_2 \sqrt{\eps}}\right\}
        \right)du \right]
        \right)
    \end{multline*}
whenever $0 \leq t_1<t_2< \dots<t_m \leq s < t,$ and $h_1, \ldots,h_m \in C_b(\R^2).$ Finally, combining this
inequality and inequality (\ref{eqn:exp-residence}) with $v(h) = \exp \left\{ - \frac{1}{C_2} \left| h
\right|^\alpha  \right\}, \alpha = 1/2p$ and $ \delta = \eps^p$, we obtain (\ref{eqn:eps-G-f}). This
completes the proof of Proposition \ref{prop:averaging-partV}.
\end{proof}

\section{The Martingale Problem}
\noindent
We want to prove that a family of probability measures converges to a unique solution $\PP^*$ of
the martingale problem for $\left(\LL^*, \delta_{[\x]} \right)$, and to identify the limiting generator
$\LL^*$ and its domain.

\subsection{The limiting generator}
\noindent Having the homeomorphism $p$ into $\R^3$, we can define the metric $\rho$ on $\G$ by $$ \rho \left(
[\x],[\y]\right) = \left\|p \left([\x]\right) - p \left([\y]\right) \right\|_{\R^3} $$ for any $[\x], [\y]
\in \G$. The constructed metric space $\left( \G, \rho \right)$ is Polish. Using Prohorov's theorem, we can
show that the family of probability measures $\PP^{\eps,*}  \in \mathcal{P} \left(C \left( [0, \infty), \G
\right) \right)$ is tight in the Prohorov topology. We choose to omit the proof, since it is analogous to the
proof of tightness in \cite{S.1}. By Prohorov's theorem, tightness implies that $\left\{\PP^{\eps,*}
\right\}$ has a cluster point. We will show that any such cluster point satisfies a certain martingale
problem. Our next step is to identify the generator of the limiting martingale problem.

By properties of $K$, $K^{-1} \left( K(\x)\right) = H^{-1} \left(H(\x) \right)$ for any $\x \notin \overV$,
and there exists a smooth extension of $K(\x)= \left(H(\x)\right)^{1/n}$ inside $\overV$.

Given a function $F$ and a set $[A] \in \G$, we denote the restriction of $F$ to $[A]$ as $F_A$. In
particular, if $F \in C( \G ; \R),$ then $F \circ \pi \in C(\overI,\R)$ and
    \begin{itemize}
    \item $F_{\V} \circ \pi(\x) = F_{\V}(\x)$ for $\x \in \V;$
    \item $F_{\overU} \circ \pi (\x) = \left( F_{\overU} \circ K^{-1}\right) \circ K$ for $\x \in
    \overU;$
    \item $F_K:= F_{\overU}\circ K^{-1}$ maps $ \left[ 0,K^* \right]$ into $\R;$
    \item $F \circ \pi(\x)= F_\V (\x) {1\!\!1}_{\V}(\x) + F_K(K(\x)) {\II}_{\overI \setminus \V}(\x)$ for $ \x \in \overI.$
    \end{itemize}

If $F \in C^2 \left(\G_\V \cup \G_{\U}; \R \right)$ then $F \circ \pi \in C^2 \left(\I \setminus \partV; \R
\right)$, and for every $t \in \R$ the operator $ \LL_t(F \circ \pi) $ is well-defined on $\I \setminus
\partV .$ If $\M \left(\LL(F \circ \pi) \right) \in \D_\A,$ then for $[\x] \in \G_\V \cup \G_{\U},$ we set
    \begin{equation*}
    \left(\A^* F \right)([\x])
    = \begin{cases}
    \bb \left(K(\x)\right)F_K'(K(\x))
    + \frac{1}{2}\sg^2 \left(K(\x)\right)
        F_K''(K(\x)) & \text{ if } [\x] \in \G_{\U}, \\
    \biggl(\M \left(\LL F_\V\right) \biggr)(\x) & \text{ if } [\x] \in \G_\V.
  \end{cases}
    \end{equation*}
    where
    \begin{equation*}
    \bb
    := \A_K \Biggl(\M \biggl( \left(\LL K \right)
            - \left(\nabla \Phi_{\psi},b \right)
                        \biggr)
                        \Biggr)
    \end{equation*} and
    \begin{equation*}
    \sg^2
    := \A_K \biggl(\M \langle d K, d K \rangle \biggr).
    \end{equation*}

Also, for future reference, define the second-order differential operator $\M \LL^K$ on $C^2 \left(\left(0,
K^*\right) \right)$ by
    \begin{equation*}
    \left(\M \left( \LL^K f\right)\right)(h) = \bb(h)f'(h)
    + \frac{1}{2}\sg^2 \left(h\right) f''(h).
    \end{equation*}
For every $\x \in \U \setminus \overV$ and for every $F \in C^2 \left(\G_\V \cup \G_{\U} \right) $ such that $\M
\left(\LL \left( F \circ \pi \right)\right) \in \D_{\A}$, we have $\left(\A^* F \right)([\x]) = \left( \M \left( \LL^K
F_K\right)\right)\left(K(\x)\right).$

We introduce two gluing operators. For $F $ with $F_{\V}\in C^1 \left(\G_\V\right), $ we define an 'inner'
gluing operator by
\begin{equation*}
\begin{split}
\underline{G} F & := \lim_{\substack{[\x] \rightarrow \G_{\partV} \\
[\x] \in \G_{\V}}} \biggl(  \A_0  \left( \M \langle dF_{\V}, d K \rangle\right) \biggr)([\x])\\
& = \lim_{\substack{[\x] \rightarrow \G_{\partV} \\
        [\x] \in\G_{\V}} }
        \left(\int_{\y \in [\x]} \frac{1}{\|\nabla K(\y)\|}\ dl(\y)
        \right)^{-1}
        \int_{\y \in [\x]}
        \frac{\left(\M\langle d F_{\V}, d K \rangle\right)(\y)} {\|\nabla K(\y)\|} \ dl(\y),
    \end{split}
\end{equation*} if this limit exists.
For $F$ with $F_{\U} \in C^1 \left(\G_{\U}\right), $ we define an 'outer' gluing operator
    \begin{equation*}
    \begin{split}
    \overline{G} F
    &:= \lim_{\substack{[\x] \rightarrow \G_{\partV} \\
        [\x]\in \G_{\U}}}
    \biggl( \A_0
        \left( \M \langle dF_K, dK \rangle\right)
    \biggr)([\x])\\
    & = \lim_{\substack{[\x] \rightarrow \G_{\partV} \\
        [\x] \in\G_{\U}} }
        \left(\int_{\y \in [\x]} \frac{1}{\|\nabla K(\y)\|}\ dl(\y)
        \right)^{-1}
        \int_{\y \in [\x]}
        \frac{\left(\M \langle d F_K, d K \rangle\right)(\y)} {\|\nabla K(\y)\|} \ dl(\y),
    \end{split}
\end{equation*} if this limit exists.

    \begin{definition} \label{def:D*}
    Define the domain
    \begin{equation*}
    \begin{split}
    \D^*:= & \biggl\{ F \in C(\G):
        F|_{\G_\V \cup \G_{\U} } \in  C^2 \left(\G_\V \cup \G_{\U} \right),
        \M \left(\LL(F \circ \pi) \right) \in \D_\A,
        \biggr. \\
    & \quad \biggl.
        \underline{G}F = \overline{G}F,
    \lim_{[\y] \rightarrow \G_{\partI}} (\A^* F) ([\y]) = 0
     \biggr\}.
     \end{split}
     \end{equation*}
For any $F \in \D^*,$ we define
     \begin{equation*}
     (\LL^*F)([\x]) := \lim_{\substack{[\y] \rightarrow [\x] \\ [\y] \in \G_{\U} \cup \G_{\V} }}(\A^* F)([\y]), \quad [\x] \in \G.
    \end{equation*}
    \end{definition}

\subsection{Approximate test functions}
\noindent
If $F \in \D^*,$ then $ F \in C(\G),$ and for all $\x \in
\overI,$
\begin{equation*}
\begin{split}
    F \circ \pi(\x)
    & = F_\V (\x) \II_{\V}(\x) + F_K(K(\x)) \II_{\overI \setminus \V}(\x).
    \end{split}
    \end{equation*}
Unfortunately, given $F \in \D^*,$ we are not guaranteed that $F \circ \pi  \in \D^\eps.$ According to the
definition of $\D^*,$ $F|_{\G_\V \cup \G_{\U} } \in  C^2 \left(\G_\V \cup \G_{\U} \right),$ that is, function
$F$ has enough smoothness on both $\G_{\V}$ and $\G_{\U}.$ But this does not immediately guarantee enough
smoothness at $\G_{\partV}.$

\begin{proposition} \label{lemma:FV-FK}
    Let $F \in \D^*,$ then $ F_{\V} \in \HH^2(\V)$ and $ F_K \in C^2
    \left(\left[0,K^*\right]\right).$
 \end{proposition}

 \begin{proof}[\bf{Proof}] To study the regularity
properties of $F$ near $\G_{\partV},$ we will use the continuity of $\LL^* F$ on $\G.$

\begin{lemma} \label{lemma:F-V}
Let $g \in C\left( \overV\right)$ and $r \in \R.$ There exists a unique solution $u_{g,r} \in \HH^2 \left(\V
\right) \cap C \left(\overV \right)$ of the PDE
\begin{equation} \label{eqn:inhomPDE}
\begin{split}
\M \left(\LL  u_{g,r} \right)  = g \quad \text{in} \quad \V, \quad u_{g,r}  = r \quad \text{on} \quad \partV.
\end{split}
\end{equation}
If $g \in C^k \left(\overV \right)$ and $k \geq 0,$ then $u_{g,r} \in \HH^{k+2} \left(\V \right) \cap C^k
\left(\overV\right)$ and there exists a constant $C_k$ such that
\begin{equation} \label{ineq:regularity}
\begin{split}
\left\|u_{g,r}\right\|_{\HH^{k+2}\left(\V \right)} & \leq C_k
    \left\{ \left\|g\right\|_{C^k \left(\overV \right)}+|r|\right\}\\
\left\|u_{g,r}\right\|_{C^k \left(\overV\right)} & \leq C_k
    \left\{ \left\|g\right\|_{C^k \left(\overV \right)}+|r|\right\}.
\end{split}
\end{equation}
\end{lemma}
\begin{proof}[\bf{Proof}]
Let $v \in \HH_0^1 \left( \V \right)$ be a weak solution of the PDE
\begin{equation*}
\begin{split}
\M \left(\LL v \right) & = 0 \quad \text{in} \quad \V,\\
v & = 0 \quad \text{on} \quad \partV.
\end{split}
\end{equation*}
We know that $\partV \in C^\infty$ and all the coefficients of $\M \LL$ are smooth. Then, by the infinite
differentiability up to the boundary \cite[Theorem 6, pg. 326]{E.}, $v \in C^\infty \left(\overV \right).$ By
the weak maximum principle \cite[Theorem 1, pg. 327]{E.}, $v \equiv 0$ on $\overV.$ The Fredholm alternative
implies that a unique weak solution $u \in \HH^1_0 \left(\V \right)$ of the corresponding inhomogeneous PDE
(\ref{eqn:inhomPDE}) exists. By boundary $\HH^2$-regularity \cite[Theorem 4, pg. 317]{E.}, $u \in \HH^2
\left(\V \right)$. By standard Sobolev theory, $u \in C \left(\overV \right).$

Now, we define $u_{g,r} := u + r.$ Notice that this $u_{g,r}$ satisfies (\ref{eqn:inhomPDE}). Let $u'_{g,r}$
be another solution of (\ref{eqn:inhomPDE}), and take $u'= u_{g,r}- u'_{g,r}$ on $\overV.$ Then $\M \left(\LL
u'\right)  = 0$ on $\V,$ and by the infinite differentiability in the interior \cite[Theorem 3, pg. 316]{E.},
$u' \in C^\infty \left(\V \right).$ On the other hand, $u'_{g,r} \in C \left(\overV \right)$ implies that $u'
\in C \left(\overV \right)$, and by the maximum principle, $u'\equiv 0$ on $\overV.$

If $g \in C^k \left( \overV \right),$ then $g \in \HH^k \left( \V \right),$ and by higher boundary regularity
\cite[Theorem 5, pg.~323]{E.}, $u \in \HH^{k+2} \left( \V \right),$ and there exists a constant $C_k$ such
that inequalities (\ref{ineq:regularity}) hold.
\end{proof}

\begin{lemma} \label{lemma:F-K}
Let $g \in C^1\left( [0,K^*]\right)$ and $r_1,r_2 \in \R.$

There exists a unique solution $u_{g,r_1,r_2} \in C^2 \left([0,K^*] \right)$ of the PDE
\begin{equation} \label{eqn:inhomPDE-exp}
\begin{split}
\M \left(\LL^K  u_{g,r_1,r_2} \right) & = g \text{ on } \left(0,K^*\right),\, u_{g,r_1,r_2} (0)  = r_1,
\,u_{g,r_1,r_2} \left(K^* \right)  = r_2.
\end{split}
\end{equation}
If $g \in C^k \left(\left[0,K^* \right] \right)$ and $k \geq 0,$ then $u_{g,r_1,r_2} \in C^{k+1}
\left(\left[0,K^* \right] \right)$ and there exists a constant $C_k$ such that
\begin{equation} \label{ineq:bound-u}
\begin{split}
\left\|u_{g,r_1,r_2}\right\|_{C^{k+1}\left(\left[0,K^*\right] \right)} & \leq C_k
    \left\{ \left\|g\right\|_{C^k \left(\left[0,K^* \right] \right)}
    + \left|r_1 \right|+ \left|r_2 \right|\right\}.
\end{split}
\end{equation}
\end{lemma}

\begin{proof}[\bf{Proof}]
Let
    \begin{equation*}
    I \left(h \right) = \exp \left\{2 \int_{0}^h \frac{\bb(s)}{\sg^2(s)}
    ds\right\}, \quad h \in \left[0,K^* \right].
    \end{equation*}
Notice that
    \begin{multline*}
    u_{g,r_1,r_2}(h)
        = \frac{1}{\int_0^{K^*} \frac{ds}{I(s)}}
            \left(r_1 \int_h^{K^*} \frac{ds}{I(s)} + r_2 \int_0^{h}
            \frac{ds}{I(s)}\right)\\
        + 2 \int_0^h \frac{1}{I(q)} \left(\int_0^q
            \frac{g(s)}{\sg^2(s)}I(s) ds \right)dq \\
        - 2 \frac{ \int_0^h \frac{ds}{I(s)}}{\int_0^{K^*}
        \frac{ds}{I(s)}} \int_0^{K^*} \frac{1}{I(q)} \left(\int_0^q
    \frac{g(s)}{\sg^2(s)}I(s) ds \right) dq,
    \end{multline*}
solves the PDE~(\ref{eqn:inhomPDE-exp}) uniquely.

If $g \in C^k \left(\left[0,K^*\right] \right),$ then $u_{g,r_1,r_2}\in C^{k+1} \left(\left[0,K^*\right]
\right)$ and there exists a constant $C_k$ such that (\ref{ineq:bound-u}) holds.
\end{proof}

Now, apply Lemma~\ref{lemma:F-V} with $
    g  = \M \left( \LL F_{\V} \right)$ and $
    r  = F \left( \partV \right)
    $
to see that $u_{g,r} = F_{\V} \in \HH^2\left(\V \right).$ Lemma~\ref{lemma:F-K} with
    $
    g  = \M \left( \LL^K F_{K} \right)$,
    $r_1  = F_K (0),$ and $
    r_2 = F_K \left( K^* \right)
    $
implies that $u_{g,r_1,r_2} = F_K \in C^2 \left(\left[0,K^*\right] \right).$ Therefore, the proof of
Proposition~\ref{lemma:FV-FK} is complete.
\end{proof}

Although $F \circ \pi$ itself is not guaranteed to belong to $\D^\eps,$ it can be approximated by a collection of
functions $\left\{F^\eps, \eps >0\right\}$ from $ C^2 \left(\R^2\right) \subset \D^\eps,$ so that
\begin{equation*}
    \lim_{\eps \rightarrow 0}
    \sup_{\x \in \overI}\left| F^\eps (\x)
        - F \circ\pi(\x)\right|
    =0.
\end{equation*} We will define $F^\eps$ similarly to the approximate test functions defined in \cite{S.1}.
The special construction of these functions is sketched below and will be used later.

First, we define the cutoff function $v_+ \in C^\infty(\R)$ by
\begin{equation*}
    v_+(h) :=\begin{cases}
    0 & \text{ if } h \leq 1, \\
    1 & \text{ if } h \geq 2.
  \end{cases}
 \end{equation*}

Next, we build suitable extensions of the functions $F_{\V}$ and $F_K$.
 \begin{lemma} \label{lemma:extension-FV}
Let $F \in \D^*.$ Suppose that $\partV \subset Z_1 \subset \R^2.$ There exists a family of functions
$\left\{F^\eps_{\V}, \,\eps > 0\right\}$, such that every $F^\eps_{\V} \in C^8_c\left(Z_1\right),$ and
    \begin{equation*}
    \begin{split}
    & F^\eps_{\V}(\x) =  F(\G_{\partV})  \text{ for all } \x \in
    \partV, \eps>0,\\
    & \lim_{\eps \rightarrow 0} \left\| F^\eps_{\V} -
        F_{\V}\right\|_{\HH^2(\V)} = 0,\\
    & \lim_{\eps \rightarrow 0} \left\|
    \M \left( \LL F^\eps_{\V}\right) -
        \LL^* F_{\V}\right\|_{C(\overV)} = 0,\\
    & \overline{\lim}_{\eps \rightarrow 0}\eps^{kq}
        \left\| F^\eps_{\V}\right\|_{C^k\left(Z_1\right)} < \infty, k \leq 8.
    \end{split}
    \end{equation*}
\end{lemma}

 \begin{lemma} \label{lemma:extension-FK}
Let $F \in \D^*$. Suppose that $\left[0,K^*\right] \subset Z_2 \subset \R.$ There exists a sequence of
functions $\left\{F^\eps_{K}, \eps > 0\right\}$ such that every $F^\eps_{K} \in C^8_c\left(Z_2\right)$ and
    \begin{equation*}
    \begin{split}
    & F^\eps_{K}(0) = F(\G_{\partV}) \text{ for all } \eps>0,\\
    & \lim_{\eps \rightarrow 0}
        \left\| F^\eps_{K} - F_{K}
        \right\|_{C^2\left( \left[0,K^*\right]\right)} = 0,\\
    & \lim_{\eps \rightarrow 0}
        \left\|\M \left( \LL^K F^\eps_{K} \right)-
        \left(\LL^* F_K\right)^{K}
        \right\|_{C\left( \left[0,K^* \right] \right)} = 0,\\
    & \overline{\lim}_{\eps \rightarrow 0}\eps^{kq}
        \left\| F^\eps_{K}\right\|_{C^k\left(Z_2\right)} < \infty, k \leq 8.
    \end{split}
    \end{equation*}
 \end{lemma}

Lemma~\ref{lemma:extension-FV} and  Lemma~\ref{lemma:extension-FK} can be proved in the same fashion as Lemma 4.5 and
Lemma 4.6 from \cite{S.1}, but instead of the properties of the generator $\mathcal{L}$ from \cite{S.1}, we have to pay
attention to the continuity properties of the time average of $\LL^\eps_{t/\eps^2}$.

Finally, given any $\x \in \overI, \eps>0,$ we let
\begin{equation} \label{eqn:Feps-def}
F^\eps(\x) :=
    F^\eps_{\V}(\x) \left(1- v_+ \left(\frac{K(\x)}{\eps^p}\right) \right)
    + F^\eps_K \left(K(\x) \right) v_+ \left(\frac{K(\x)}{\eps^p}
    \right),
\end{equation}
where $F^\eps_\V$ is constructed in Lemma~\ref{lemma:extension-FV}, and $F_K^\eps$ is constructed in
Lemma~\ref{lemma:extension-FK}.

\subsection{Gluing} \noindent
The hypothesis that $\underline{G}F = \overline{G}F$  will be called the `gluing
requirement'. It can be rewritten in the form
    \begin{equation*}
    \begin{split}
    & \lim_{\substack{[\x] \rightarrow \G_{\partV} \\ [\x] \in \G_\V}}
   \biggl( \A_0 \left( \M \langle dF_{\V}, dK \rangle\right)
        \biggr)([\x])
    = \lim_{\substack{[\x] \rightarrow \G_{\partV} \\ [\x] \in \G_{\U} }}
    \biggl(\A_0 \left( \M \langle dF_K, dK \rangle \right)
        \biggr)([\x]).
    \end{split}
    \end{equation*}
Notice that if $ [\x] \in \G_{\V},$ then
    \begin{equation*}
    \begin{split}
    \biggl( \A_0 \left( \M \langle dF_{\V}, dK \rangle \right)
        \biggr)([\x])
    = \biggl( \A_K \left(\M \langle dF_{\V}, d K \rangle \right) \biggr) \left(K(\x) \right).
    \end{split}
    \end{equation*}
If $ [\x] \in \G_{\U},$ then
    \begin{multline*}
    \biggl( \A_0 \left(\M \langle d F_K, dK \rangle
        \right)\biggr)([\x])
        = F'_K(K(x))
            \biggl( \A_0 \left(\M \langle d K, d K \rangle \right)\biggr)([\x])\\
            = \sg^2(K(\x)) F'_K(K(\x)).
    \end{multline*}
    So, formally, the gluing requirement demands that
\begin{equation*}
    \begin{split}
   \lim_{h \uparrow 0}
     \biggl(\A_K \left(\M \langle dF_{\V}, d K \rangle \right)\biggr) (h)=
      \lim_{h \downarrow 0} \sg^2(h) F'_K(h).
    \end{split}
    \end{equation*}

\begin{proposition}
Let $ F \in C(\G)$ be such that $ F|_{\G_\V \cup \G_{\U} } \in C^2 \left(\G_\V \cup \G_{\U} \right)$ and $ \M
\left(\LL(F \circ \pi) \right) \in \D_\A$, then both $\underline{G}F$ and $\overline{G}F$ exist.
\end{proposition}

\begin{proof}[\bf{Proof}]
From Lemma~\ref{lemma:FV-FK}, $ \lim_{h \downarrow 0} \left(F_K \right)'(h)$ exists; we denote this limit by
$\left(F_K \right)'(0).$  Notice that from our assumptions on function $K$ the limit $\lim_{h \downarrow 0}
\sg^2(h) = \sg^2(0)$ exists and $\sg^2(0)>0$. Then the quantity $\overline{G}F$ is well-defined and
\begin{equation*} \overline{G}F = \sg^2(0) \left(F_K \right)'(0).
\end{equation*}
For every $h \in (-a,a),$ we introduce a bounded linear operator $\T_h: \HH^1(\V) \rightarrow
L^2\left(K^{-1}(h)\right)$ such that if $f \in \HH^1(\V) \cap C(\overV)$ then  $\T_h f = f|_{\partV}$ and
    \begin{equation} \label{ineq:trace}
    \left\| \T_h f\right\|_{L^2(\partV)} \leq C \|f\|_{\HH^1(\V)}.
    \end{equation}
The Trace theorem~\cite[Theorem 1, pg. 258]{E.} guarantees the existence of such operators $\T_h$. For any $f
\in \HH^1\left(\V\right)$ and for $h \in (-a,0]$ we define
    \begin{equation*}
    \left(\A^{\tt Tr}f\right)(h)
    = \left(\int_{\y \in K^{-1}(h)} \frac{1}{\|\nabla K(\y)\|}\ dl(\y)
        \right)^{-1}
        \int_{\y \in K^{-1}(h)}
        \frac{\left(\T_h f\right)(y)} {\|\nabla K(\y)\|} \ dl(\y).
    \end{equation*}

If $f \in C\left( \overV\right),$ then $\A^{\tt Tr} f = A_K f$ on $(-a,0]$. By Lemma \ref{lemma:FV-FK}, $F_\V
\in \HH^2 \left(\V \right)$, which implies that $\biggl(\A^{\tt Tr} \left( \M \langle dF_{\V}, dK \rangle
\right)\biggr)(0)$ is well-defined and that
    \begin{equation} \label{eqn:G-inner}
    \underline{G}F
    = \lim_{h \uparrow 0} \biggl(\A^{\tt Tr}
    \left( \M \langle dF_{\V}, dK \rangle \right)\biggr)(h),
    \end{equation}
if this limit exists. We consider the sequence of functions $F_\V^\eps$ from Lemma \ref{lemma:extension-FV}.
By property (\ref{ineq:trace}) of the trace operator, there exists a constant $C>0$ such that
    \begin{equation*}
    \left|\biggl(\A^{\tt Tr}
    \left( \M \langle dF_{\V}, dK \rangle \right)\biggr)(h)
    - \biggl( \A^{\tt Tr}
    \left( \M \langle dF^\eps_{\V}, dK \rangle \right)\biggr)(h) \right|
    \leq C\left\| F_{\V} - F^\eps_{\V}\right\|_{\HH^2(\V)}
    \end{equation*}
for all $h \in (-a,0].$

The following result is borrowed from \cite{S.1} (Lemma 9.6).
\begin{lemma} \label{lemma:h1-h2-trace}
There exists a constant $C>0$ such that if $f \in \HH^2(\R^2),$ then
\begin{equation} \label{eqn:h1-h2-trace}
    \left|\biggl(\A^{\tt Tr}
    f\biggr)(h')
    - \biggl( \A^{\tt Tr}
    f\biggr)(h'')
    \right| \leq C \|f\|_{\HH^1(\R^2)} \sqrt{\left|h'-h''\right|}
    \end{equation} for all $h',h'' \in (-a,a).$
\end{lemma}


From (\ref{eqn:h1-h2-trace}), there exists a constant $C>0$ such that
    \begin{equation*}
    \begin{split}
    & \left|  \biggl(\A^{\tt Tr}
    \left( \M \langle dF_{\V}, dK \rangle \right)\biggr)(h)
    - \biggl( \A^{\tt Tr}
    \left( \M \langle dF_{\V}, dK \rangle \right)\biggr)(0)
    \right|\\
     & \leq \left|\biggl(\A^{\tt Tr}
    \left( \M \langle dF_{\V}, dK \rangle \right)\biggr)(h)
    - \biggl( \A^{\tt Tr}
    \left( \M \langle dF^\eps_{\V}, dK \rangle \right)\biggr)(h)
    \right|\\
    & \quad + \left|\biggl(\A^{\tt Tr}
    \left( \M \langle dF^\eps_{\V}, dK \rangle \right)\biggr)(h)
    - \biggl( \A^{\tt Tr}
    \left( \M \langle dF^\eps_{\V}, dK \rangle \right)\biggr)(0)
    \right|\\
    & \quad + \left|\biggl(\A^{\tt Tr}
    \left( \M \langle dF^\eps_{\V}, dK \rangle \right)\biggr)(0)
    - \biggl( \A^{\tt Tr}
    \left( \M \langle dF_{\V}, dK \rangle \right)\biggr)(0)
    \right| \\
    & \leq 2 C \left\| F_{\V} - F^\eps_{\V}\right\|_{\HH^2(\V)}
    + C \sqrt{|h|} \left\| F^\eps_{\V}\right\|_{\HH^2(\V)}.
    \end{split}
    \end{equation*}
Therefore the limit in (\ref{eqn:G-inner}) exists and $ \underline{G}F = \biggl(\A^{\tt Tr} \left( \M \langle
dF_{\V}, dK \rangle \right)\biggr)(0) $ is well-defined.
\end{proof}

\subsection{The limiting martingale problem}
\noindent
To prove that any cluster point of the $\PP^{\eps,*}$s satisfies the martingale problem for
$\left(\LL^*,\delta_{[\x]}\right)$, we will need to connect elements of the limiting domain $\D^*$ back to
the functions from $\D^\eps$. The propositions and lemmas below allow us to establish this connections.

Recall $F^\eps,$ constructed above to satisfy (\ref{eqn:Feps-def}).
\begin{proposition}
Given $\x \in \overI,$ for every $t \geq 0$ and for every $\eps>0$
\begin{equation} \label{eqn:Leps-Feps}
\begin{split}
    \left(\LL^\eps_{t/\eps^2}F^\eps\right)(\x)
     = L^\eps \left( \x, \frac{t}{\eps^2} \right)
        + \sum_{k=0}^3 G^\eps_i \left( \x, \frac{t}{\eps^2} \right)
        + R^\eps \left( \x, \frac{t}{\eps^2} \right),
    \end{split}
\end{equation} where
    \begin{equation*}
    \begin{split}
    L^\eps \left( \x, t \right)
    & = \left( 1 - v_+  \left(\frac{K(\x)}{\eps^p}\right)\right)
        \left(\LL_t F_{\V}^\eps\right)(\x)
        + v_+  \left(\frac{K(\x)}{\eps^p} \right)
        \LL_t F^\eps_K \left(K(\x) \right)\\
    & \quad + \frac{1}{\eps}
        v_+  \left(\frac{K(\x)}{\eps^p}\right)
        \left(F_K^\eps \right)'(K(\x))
        \left(\nabla K(\x),b(\x,t) \right),\\
    G^\eps_0(\x,t) & = \frac{1}{\eps^{p+1}} v'_+  \left(\frac{K(\x)}{\eps^p}\right)
        \left( \nabla K(\x) ,b\left(\x,t\right)\right)
        \biggl( F^\eps_K \left(K(\x) \right) - F_{\V}^\eps(\x)\biggr),\\
    G^\eps_1(\x,t) & = \frac{1}{\eps^p}v'_+  \left(\frac{K(\x)}{\eps^p}\right)
        \biggl(\left\langle d F_{\V}^\eps,
        d K \right\rangle_t \left(\x\right)
        - \left\langle d \left(F_K^\eps \circ K\right),
        d K\right\rangle_t (\x)\biggr),\\
    G^\eps_2(\x,t) & = \frac{1}{2 \eps^{2p}}
        v''_+  \left(\frac{K(\x)}{\eps^p}\right)
        \biggl( F^\eps_K \left(K(\x) \right) - F_{\V}^\eps(\x)\biggr)
        \langle d K,d K \rangle_t (\x),\\
    G^\eps_3(\x,t) & = \frac{1}{\eps}
        \left( 1 - v_+  \left(\frac{K(\x)}{\eps^p}\right)\right)
        \left(\nabla F_{\V}^\eps(\x),b(\x,t) \right),\\
    R^\eps(\x,t)& = \frac{1}{ \eps^{p}}
        v'_+  \left(\frac{K(\x)}{\eps^p}\right)
        \biggl( F^\eps_K \left(K(\x) \right) - F_{\V}^\eps(\x)\biggr)
        \left(\LL_t K\right)(\x).
    \end{split}
    \end{equation*}
\end{proposition}
\begin{proof}[\bf{Proof}]
Below, for any appropriate $f$ and $g$ involved, we understand $\LL^\eps_{t/\eps^2} f(g(\x))$ and $df(g(\x))$
as $\left(\LL^\eps_{t/\eps^2} (f \circ g)\right)(\x)$ and $d(f \circ g)(\x)$ correspondingly. For example,
    \begin{equation*}
    \begin{split}
    \LL^\eps_{t/\eps^2} v_+\left(\frac{K(\x)}{\eps^p} \right)
    & =  \frac{1}{\eps}
        \left( \nabla v_+\left(\frac{K(\x)}{\eps^p} \right),
        b\left(\x,\frac{t}{\eps^2}\right) \right) + \LL_{t/\eps^2} v_+\left(\frac{K(\x)}{\eps^p} \right)\\
    & = \frac{1}{\eps^{p+1}}
        v'_+  \left(\frac{K(\x)}{\eps^p}\right)
        \left( \nabla K(\x) ,b\left(\x,\frac{t}{\eps^2}\right)
        \right)\\
    & \quad + \frac{1}{2}
        \sum_{i,j=1}^2 \gamma_{ij}\left(\x,\frac{t}{\eps^2}\right)
        \frac{\partial^2 v_+\left(\frac{K(\x)}{\eps^p} \right)}
        {\partial x_i \partial x_j}(\x).
        \end{split}
    \end{equation*} This implies that
    \begin{align}
    \LL^\eps_{t/\eps^2} v_+\left(\frac{K(\x)}{\eps^p} \right)
    & = \frac{1}{\eps^{p+1}} v'_+  \left(\frac{K(\x)}{\eps^p}\right)
        \left( \nabla K(\x) ,b\left(\x,\frac{t}{\eps^2}\right)
        \right) \notag\\
    & \quad + \frac{1}{ \eps^{p}}
        v'_+  \left(\frac{K(\x)}{\eps^p}\right)
        \left(\LL_{t/\eps^2} K\right)(\x) \notag\\
    & \quad + \frac{1}{2 \eps^{2p}}
        v''_+  \left(\frac{K(\x)}{\eps^p}\right)
        \langle d K,d K \rangle_{t/\eps^2} (\x) \label{eqn:Leps-w}
        \end{align}
    and
    \begin{equation} \label{eqn:dw}
    d v_+\left(\frac{K(\x)}{\eps^p}\right)
    = \frac{1}{\eps^p} v'_+  \left(\frac{K(\x)}{\eps^p}\right) dK(\x).
    \end{equation}
We know that $F^\eps \in C^2\left(\R^2\right) \subset \D^\eps.$ Applying the generator $\LL^\eps_{t/\eps^2}$
to the function $F^\eps $ defined in (\ref{eqn:Feps-def}), we obtain
    \begin{equation*}
    \begin{split}
    \left( \LL^\eps_{t/\eps^2}  F^\eps \right)(\x)
    & = \left( 1 - v_+  \left(\frac{K(\x)}{\eps^p}\right)\right)
        \left( \frac{1}{\eps}
        \left(\nabla F_{\V}^\eps(\x),b\left(\x,\frac{t}{\eps^2} \right)
        \right) +
        \left(\LL_{t/\eps^2}F_{\V}^\eps\right)(\x) \right)\\
    & \quad + v_+  \left(\frac{K(\x)}{\eps^p} \right)
        \LL_{t/\eps^2}F^\eps_K \left(K(\x) \right)\\
    & \quad + \frac{1}{\eps}
        v_+  \left(\frac{K(\x)}{\eps^p}\right)
        \left(F_K^\eps \right)'(K(\x))
        \left(\nabla K(\x),b \left(\x,\frac{t}{\eps^2}\right) \right)\\
    & \quad + \biggl( F^\eps_K \left(K(\x) \right) - F_{\V}^\eps(\x)\biggr)
         \LL_{t/\eps^2} v_+  \left(\frac{K(\x)}{\eps^p} \right)\\
    & \quad - \left\langle d F_{\V}^\eps\left(\x\right),
        d v_+\left(\frac{K(\x)}{\eps^p} \right)\right\rangle_{t/\eps^2}\\
    & \quad + \left\langle d F_K^\eps\left(K(\x)\right),
        d v_+\left(\frac{K(\x)}{\eps^p}
        \right)\right\rangle_{t/\eps^2}.
    \end{split}
    \end{equation*}
We unravel this expression further, using (\ref{eqn:Leps-w}) and (\ref{eqn:dw}), to see that
    \begin{equation*}
    \begin{split}
    \left( \LL^\eps_{t/\eps^2} F^\eps \right) &(\x)
    = \left( 1 - v_+  \left(\frac{K(\x)}{\eps^p}\right)\right)
        \left( \frac{1}{\eps}
        \left(\nabla F_{\V}^\eps(\x),b\left(\x,\frac{t}{\eps^2} \right)
        \right) +
        \left(\LL_{t/\eps^2}F_{\V}^\eps\right)(\x) \right)\\
    & + v_+  \left(\frac{K(\x)}{\eps^p} \right)
        \LL_{t/\eps^2}F^\eps_K \left(K(\x) \right)\\
    & + \frac{1}{\eps}
        v_+  \left(\frac{K(\x)}{\eps^p}\right)
        \left(F_K^\eps \right)'(K(\x))
        \left(\nabla K(\x),b \left(\x,\frac{t}{\eps^2}\right)
        \right)\\
    &  + \frac{1}{\eps^{p+1}} v'_+  \left(\frac{K(\x)}{\eps^p}\right)
        \left( \nabla K(\x) ,b\left(\x,\frac{t}{\eps^2}\right)\right)
        \biggl( F^\eps_K \left(K(\x) \right) - F_{\V}^\eps(\x)\biggr)\\
    &  + \frac{1}{ \eps^{p}}
        v'_+  \left(\frac{K(\x)}{\eps^p}\right)
        \biggl( F^\eps_K \left(K(\x) \right) - F_{\V}^\eps(\x)\biggr)
        \left(\LL_{t/\eps^2} K\right)(\x)\\
    &  + \frac{1}{2 \eps^{2p}}
        v''_+  \left(\frac{K(\x)}{\eps^p}\right)
        \biggl( F^\eps_K \left(K(\x) \right) - F_{\V}^\eps(\x)\biggr)
        \langle d K,d K \rangle_{t/\eps^2} (\x)\\
    & - \frac{1}{\eps^p}v'_+  \left(\frac{K(\x)}{\eps^p}\right)
        \left(\left\langle d F_{\V}^\eps,
        d K \right\rangle_{t/\eps^2}\left(\x\right)
        - \left\langle d \left(F_K^\eps \circ K\right),
        d K\right\rangle_{t/\eps^2}(\x)\right).
    \end{split}
    \end{equation*}
    This immediately implies (\ref{eqn:Leps-Feps}).
    \end{proof}

\begin{proposition} \label{prop:Leps-L*f}
Fix $F \in \D^*$ and $t >0.$ Then
    \begin{equation} \label{eqn:L-ML-together}
    \overline{\lim_{\eps \rightarrow 0} }\E^\eps
    \left[\left| \int_0^{t \wedge \tau}
        \biggl( \left(\M L^\eps \right)
        \left(\X_s \right) - \left(\LL^* F\right)
        \left(\left[\X_s \right]\right)\biggr) ds\right| \right]
    = 0.
    \end{equation}
\end{proposition}

\begin{proof}[\bf{Proof}] Consider
    \begin{equation*}
    \begin{split}
    \left(\M L^\eps \right)\left( \x \right)
    & = \left( 1 - v_+  \left(\frac{K(\x)}{\eps^p}\right)\right)
       \left(\M \left(\LL F_{\V}^\eps\right)\right)(\x)\\
    & + v_+  \left(\frac{K(\x)}{\eps^p} \right)
        \left(\M \left(\LL \left(F^\eps_K \circ K \right)\right)\right)(\x)\\
    & = \left( 1 - v_+  \left(\frac{K(\x)}{\eps^p}\right)\right)
        \left(\M \left(\LL F_{\V}^\eps\right)\right)(\x)\\
    & + v_+  \left(\frac{K(\x)}{\eps^p} \right)
        \left(F^\eps_K\right)' \left(K(\x) \right)
        \left( \M \left(\LL K \right) \right)(\x)\\
    & + \frac{1}{2}v_+  \left(\frac{K(\x)}{\eps^p} \right)
        \left(F^\eps_K\right)'' \left(K(\x) \right)
        \left( \M \left\langle dK, dK \right\rangle \right)(\x).
        \end{split}
    \end{equation*}
Recall that $v, v_+ \in C^\infty(\R;[0,1])$ were defined to satisfy
\begin{equation*}
v (h)=
    \begin{cases}
    0 & \text{ if } |h| \leq 1,\\
    1 & \text{ if } |h| \geq 2
    \end{cases} \text{ and }
v_+(h)=
    \begin{cases}
    0 & \text{ if } h \leq 1,\\
    1 & \text{ if } h \geq 2.
    \end{cases}
\end{equation*}
Define $v_- \in C^\infty(\R;[0,1])$ by $v_-(h):= v_+(-h),$ then
\begin{equation*}
v_- (h)=
    \begin{cases}
    1 & \text{ if } h \leq -2,\\
    0 & \text{ if } |h| \geq -1.
    \end{cases}
\end{equation*}
Without loss of generality assume that these three functions $v,v_-$ and $v_+$ form a smooth partition of
unity. Notice that $$\supp \left(v_-\right) \subset \left(- \infty,-1\right), \supp \left(v\right) \subset
\left(-2,2 \right), \supp \left(v_+ \right) \subset \left(1,\infty \right).$$ Fix $r$ so that $\eps^{2q } \gg
\eps^r \gg \eps^p$ as $\eps \rightarrow 0$. Then we have $\supp \left(v_- \left(\frac{K(\,\cdot\,)}{\eps^r}
\right) \right) \subseteq \V$ and $ \supp \left(v_+ \left(\frac{K(\,\cdot\,)}{\eps^r} \right) \right)
\subseteq \U.$ Notice that
    \begin{equation*}
    \begin{split}
    & v_- \left(\frac{K(\x)}{\eps^r}
    \right) \left(1 -v_+  \left(\frac{K(\x)}{\eps^p}\right)
    \right) = v_- \left(\frac{K(\x)}{\eps^r}
    \right),
    \end{split}
    \end{equation*}
which implies that
    \begin{equation} \label{eqn:average-L-vV}
    \begin{split}
    \left(\M L^\eps \right)\left( \x \right)
    v_- \left(\frac{K(\x)}{\eps^r}
    \right) = \left(\M \left(\LL F_{\V}^\eps\right)\right)(\x)
        v_- \left(\frac{K(\x)}{\eps^r} \right).
        \end{split}
    \end{equation}
Now, using Lemma \ref{lemma:extension-FV} and expression (\ref{eqn:average-L-vV})
    \begin{multline} \label{eqn:L-ML-V}
    \overline{\lim_{\eps \rightarrow 0}} \
    \E^\eps
        \left[ \left|
        \int_{0}^{t \wedge \tau}
        \biggl( \left( \M L^\eps \right) \left(\X_u \right)
            - \left(\LL^* F\right)\left(\Y_u \right)
        \biggr) v_- \left(\frac{K \left(\X_u \right)}{\eps^r}
    \right) du
         \right|\right]\\
    \leq \overline{\lim_{\eps \rightarrow 0}} \
    \biggl\| \M \left(\LL F^\eps_{\V}\right)
    - \LL^* F \biggr\|_{C\left( \overV\right)} = 0.
    \end{multline}
Next, by Lemma \ref{lemma:extension-FV}, Lemma \ref{lemma:extension-FK}, and using the fact that  $\LL^* F$
is bounded, we see that
    \begin{equation*}
    \left|\biggl( \left( \M L^\eps \right) \left(\x \right)
            - \left(\LL^* F\right)\left([\x] \right)
        \biggr)
    v \left(\frac{K(\x)}{\eps^r}\right)
    \right| < C \eps^{-2 q}
    \chi_{ \left[-2 \eps^r,2 \eps^r \right]} \left(K(\x) \right).
    \end{equation*}
By Lemma \ref{lemma:residence}, there exists a positive constant $C$ such that
    \begin{multline} \label{eqn:L-ML-partV}
    \overline{\lim_{\eps \rightarrow 0}} \
    \E^\eps
        \left[ \left|
        \int_{0}^{t \wedge \tau}
        \biggl( \left( \M L^\eps \right) \left(\X_u \right)
            - \left(\LL^* F\right)\left(\Y_u \right)
        \biggr) v \left(\frac{K \left(\X_u \right)}{\eps^r}
    \right) du
         \right|\right]\\
    \leq C \eps^{r -2q} \overline{\lim_{\eps \rightarrow 0}} \
    \E^\eps
        \left[ \frac{1}{\eps^r}
        \int_{0}^{t \wedge \tau}
        \chi_{\left[ -2\eps^r, 2\eps^r\right]}
        \left(K \left(\X_u \right)\right) du \right] = 0.
    \end{multline}
Notice that
    \begin{align}
     \E^\eps
     & \left[ \left|
        \int_{0}^{t \wedge \tau}
        \biggl( \left( \M L^\eps \right) \left(\X_u \right)
            - \left(\LL^* F\right)\left(\Y_u \right)
        \biggr)
        v_+ \left(\frac{K \left(\X_u \right)}{\eps^r}\right) du
         \right|\right] \notag\\
     &  \leq \E^\eps
        \left[ \left|
        \int_{0}^{t \wedge \tau}
        \biggl( \left( \M L^\eps \right) \left(\X_u \right)
            - \left( A_0\left( \M L^\eps \right) \right) \left(\Y_u \right)
        \biggr) v_+ \left(\frac{K \left(\X_u \right)}{\eps^r}
    \right) du
         \right|\right] \notag\\
         &\quad  + \E^\eps
        \left[ \left|
        \int_{0}^{t \wedge \tau}
        \biggl( \left(A_0 \left( \M L^\eps \right)\right) \left(\Y_u \right)
            - L^\eps \left(\X_u, \frac{u}{\eps^2}\right)
        \biggr) v_+ \left(\frac{K \left(\X_u \right)}{\eps^r}
    \right) du
         \right|\right] \notag\\
    & \quad + \E^\eps
        \left[ \left|
        \int_{0}^{t \wedge \tau}
        \biggl( L^\eps \left(\X_u, \frac{u}{\eps^2}\right)
            - \left(\LL^* F\right)\left(\Y_u \right)
        \biggr) v_+ \left(\frac{K \left(\X_u \right)}{\eps^r}
    \right) du
         \right|\right], \label{eqn:MLeps-L*}
        \end{align}
where
    \begin{align}
    \biggl( \biggr. & \biggl.
        L^\eps \left(\x, \frac{t}{\eps^2}\right)
        - \left(\LL^* F \right)
        \left( [\x] \right)
    \biggr)
            v_+ \left( \frac{K \left(\x \right)}{\eps^r}\right) \notag\\
    & = \biggl(
        \LL_{t/\eps^2} F^{\eps}_K \left( K(\x)\right)
            + \frac{1}{\eps}
            \left(F^{\eps}_K \right)' \left(K(\x)\right)
            \psi \left(\x, \frac{t}{\eps^2} \right)
    \biggr)
        v_+ \left( \frac{K \left(\x \right)}{\eps^r} \right) \notag\\
    & \quad - \A_0 \biggl(
        \M \left( (\LL K)
            - \left( \nabla \Phi_{\psi}, b \right)
            \right)
            \biggr)([\x])
            F_K' \left(K(\x)\right)
            v_+ \left( \frac{K \left(\x \right)}{\eps^r} \right) \notag\\
    & \quad - \frac{1}{2}
        \A_0 \biggl( \M \left( \langle dK,dK \rangle \right)
            \biggr)([\x])
        F_K'' \left(K(\x)\right)
        v_+ \left( \frac{K \left(\x \right)}{\eps^r}\right) \notag\\
    & = v_+ \left( \frac{K \left(\x \right)}{\eps^r}\right)
    \sum_{i=1}^5 K_i^\eps \left(\x,\frac{t}{\eps^2} \right) \label{eqn:sum-Ki}
    \end{align}
and
    \begin{equation*}
    \begin{split}
    K_1^\eps \left(\x,\frac{t}{\eps^2} \right)
    & = \biggl(
        \LL_{t/\eps^2} F^{\eps}_K \left( K(\x)\right)
        - \left(\A_0
            \left( \M
                \left( \LL\left(F^\eps_K \circ K\right)
                \right)
            \right)
        \right)
        \left( [\x] \right)
    \biggr),\\
     K_2^\eps \left(\x,\frac{t}{\eps^2} \right)
     & = \biggl( \left(\A_0
            \left( \M
                \left( \LL\left(F^\eps_K \circ K\right)
                \right)
            \right)
        \right)
        \left( [\x] \right)
     - \left(\A_0
            \left( \M
                \left( \LL\left(F_K \circ K\right)
                \right)
            \right)
        \right)
        \left( [\x] \right)
    \biggr), \\
    K_3^\eps \left(\x,\frac{t}{\eps^2} \right)
    & = \biggl(F_K' \left(K(\x)\right)
        - \left(F^\eps_K \right)' \left(K(\x)\right)
        \biggr)
    \A_0 \left(
        \M \left( \nabla \Phi_{\psi}, b \right)
            \right)([\x]), \\
    K_4^\eps \left(\x,\frac{t}{\eps^2} \right)
    & = \biggl(\A_0 \left(
        \M \left( \nabla \Phi_{\psi}, b \right)
            \right)([\x])
             - \left( \nabla \Phi_{\psi} ,
                b  \right)\left(\x, \frac{t}{\eps^2}
             \right) \biggr)
            \left(F^\eps_K \right)' \left(K(\x)\right),\\
    K_5^\eps \left(\x,\frac{t}{\eps^2} \right)
    & =  \biggl( \frac{1}{\eps}
            \psi \left(\x, \frac{t}{\eps^2} \right)
             + \left( \nabla \Phi_{\psi} ,
                b  \right)\left(\x, \frac{t}{\eps^2}
             \right) \biggr)
            \left(F^\eps_K \right)' \left(K(\x)\right).
    \end{split}
    \end{equation*}
By Lemma \ref{lemma:f-AMf}, there exists a constant $C_c>0$ such that
\begin{multline} \label{eqn:K1}
\E^\eps
     \left[ \left|
        \int_{0}^{t \wedge \tau}
     K_1^\eps \left(\X_u,\frac{t}{\eps^2} \right)
        v_+ \left(\frac{K \left(\X_u \right)}{\eps^r}\right) du
         \right|\right] \\
         \leq C_c \eps^{1-r(n+1)}(1+t)
         \left\|F^\eps_K
         \right\|_{C^5\left(\left[0,K^* \right]\right)}.
    \end{multline}
    By Lemma \ref{lemma:extension-FK},
\begin{equation} \label{eqn:K2-K3}
\overline{\lim_{\eps \rightarrow 0}}\E^\eps
     \left[ \left|
        \int_{0}^{t \wedge \tau}
     K_i^\eps \left(\X_u,\frac{t}{\eps^2} \right)
        v_+ \left(\frac{K \left(\X_u \right)}{\eps^r}\right) du
         \right|\right] =0, \quad i=2,3.
    \end{equation}
By Lemma \ref{lemma:f-AMf} and Lemma \ref{lemma:extension-FK}, there exists a constant $C>0$ such that
\begin{multline}
\E^\eps
     \left[ \left|
        \int_{0}^{t \wedge \tau}
     K_4^\eps \left(\X_u,\frac{t}{\eps^2} \right)
        v_+ \left(\frac{K \left(\X_u \right)}{\eps^r}\right) du
         \right|\right] \\
        < C \eps^{1-r(n+1)}(1+t) \left\|F^\eps_K
            \right\|_{C^1\left(\left[0,K^*\right]\right)}. \label{eqn:K4}
    \end{multline}
By Lemma \ref{lemma:phi+grad}, there exists a constant $C_b>0$ such that
    \begin{multline}
\E^\eps
     \left[ \left|
        \int_{0}^{t \wedge \tau}
     K_5^\eps \left(\X_u,\frac{t}{\eps^2} \right)
        v_+ \left(\frac{K \left(\X_u \right)}{\eps^r}\right) du
         \right|\right] \\
        < C_b \eps^{1-2r}(1+t) \left\|F^\eps_K
            \right\|_{C^3\left(\left[0,K^*\right]\right)}. \label{eqn:K5}
    \end{multline}
Combining equality (\ref{eqn:sum-Ki}) with inequalities (\ref{eqn:K1}), (\ref{eqn:K2-K3}), (\ref{eqn:K4}),
and (\ref{eqn:K5}), we can see that $$ r < \frac{1-4q}{n+1}$$ guarantees that
\begin{equation} \label{eqn:ML-Leps}
\overline{\lim_{\eps \rightarrow 0}}\E^\eps
     \left[ \left|
    \int_{0}^{t \wedge \tau}
        \biggl( L^\eps \left(\X_u, \frac{u}{\eps^2}\right)
            - \left(\LL^* F\right)\left(\left[\X_u\right] \right)
        \biggr)
        v_+ \left(\frac{K \left(\X_u \right)}{\eps^r}\right) du
         \right|\right] =0.
    \end{equation}

By Lemma \ref{lemma:f-Af}, Lemma \ref{lemma:f-AMf}, and Lemma \ref{lemma:extension-FK}, there exists a
positive constant $C$ such that
 \begin{align}
     \E^\eps
        & \left[ \left|
        \int_{0}^{t \wedge \tau}
        \biggl( \left( \M L^\eps \right) \left(\X_u \right)
            - \left( A_0\left( \M L^\eps \right) \right) \left(\left[\X_u\right] \right)
        \biggr)  v_+ \left(\frac{K \left(\X_u \right)}{\eps^r}
    \right) du
         \right|\right] \notag \\
        & + \E^\eps
        \left[ \left|
        \int_{0}^{t \wedge \tau}
        \biggl( \left(A_0 \left( \M L^\eps \right)\right) \left(\left[\X_u\right] \right)
            - L^\eps \left(\X_u, \frac{u}{\eps^2}\right)
        \biggr) v_+ \left(\frac{K \left(\X_u \right)}{\eps^r}
    \right) du
         \right|\right] \notag \\
         & < C \eps^{1-r(n+1)-5q}(1+t) \label{eqn:ML-Leps-1-2}
        \end{align}
Therefore, combining (\ref{eqn:MLeps-L*}), (\ref{eqn:ML-Leps}), and (\ref{eqn:ML-Leps-1-2}), we conclude that
\begin{equation} \label{eqn:L-ML-U0}
    \begin{split}
    \overline{\lim_{\eps \rightarrow 0}} \
    & \E^\eps
        \left[ \left|
        \int_{0}^{t \wedge \tau}
        \biggl( \left( \M L^\eps \right) \left(\X_u \right)
            - \left(\LL^* F\right)\left(\left[\X_u\right] \right)
        \biggr) v_+ \left(\frac{K \left(\X_u \right)}{\eps^r}
    \right) du
         \right|\right]= 0.
    \end{split}
    \end{equation}
Finally, combine (\ref{eqn:L-ML-V}), (\ref{eqn:L-ML-partV}), and (\ref{eqn:L-ML-U0}) to obtain
(\ref{eqn:L-ML-together}), completing the proof of Proposition \ref{prop:Leps-L*f}.
\end{proof}

\begin{remark}
Notice that we require $$2 q < r < \min\left\{p, \frac{1-5q}{n+1} \right\}.$$ Therefore, we need to make sure
that $q < \min \left\{\frac{p}{2},\frac{1}{2n+7} \right\}$. Recall that $p= \frac{1}{n+1},$ which implies
that we need $q < \frac{1}{2n+7}.$
\end{remark}

    \begin{proposition} \label{prop:Geps-negl}
Fix $0 \leq t_1<t_2< \cdots<t_m \leq s < t$ and $h_1, \cdots , h_m \in C_b\left( \R^2 \right).$ If $F \in
\D^*,$ then
    \begin{equation} \label{eqn: MG's}
    \lim_{\eps \rightarrow 0} \E^\eps
    \left[ \int_{s \wedge \tau}^{t \wedge \tau}
        \left(\M G_i^\eps \right)
        \left(\X_s \right) ds \prod_{k=1}^m h_k \left(\X_{t_k} \right)
        \right]
    = 0, \quad i = 0,1,2.
    \end{equation}
    \end{proposition}

    \begin{proof}[\bf{Proof}]
    Notice that
    \begin{align}
    \left( \M G^\eps_0\right)(\x)
        & = \frac{1}{\eps^{p+1}} v'_+  \notag \\
        & \quad \times \left(\frac{K(\x)}{\eps^p}\right)
        \biggl( F^\eps_K \left(K(\x) \right) - F_{\V}^\eps(\x)\biggr)
        \left( \M \left( \nabla K ,b\right)\right)(\x) = 0, \notag\\
    \left( \M G^\eps_1\right)(\x)
         & = \frac{1}{\eps^p}v'_+  \left(\frac{K(\x)}{\eps^p}\right) \notag\\
        & \quad \times \biggl(\left( \M \left\langle d F_{\V}^\eps,
        d K \right\rangle \right) \left(\x\right)
        - \left( F_K^\eps\right)'(K(\x))
        \left( \M \left\langle d K, d K\right\rangle \right) (\x)
            \biggr), \notag\\
    \left( \M G^\eps_2\right)(\x)
        & = \frac{1}{2 \eps^{2p}}
        v''_+  \left(\frac{K(\x)}{\eps^p}\right)
        \biggl( F^\eps_K \left(K(\x) \right) - F_{\V}^\eps(\x)\biggr)
        \left(\M \langle d K,d K \rangle \right) (\x). \label{eq:EMG-i-0-2}
    \end{align}
For $\eps>0$ we define two functions $f_1^\eps, f_2^\eps \in C^4 \left(\overI \cup \left\{ \x: \eps^p < K(\x)
< 2 \eps^{p}\right\} \right)$ by
    \begin{equation*}
    f_1^\eps(\x) :=
        \left(
        \M \left\langle d F_{\V}^\eps, d K \right\rangle
        \right) \left(\x\right)
        - \left( F_K^\eps\right)'(K(\x))
        \left( \M \left\langle d K,  d K\right\rangle \right)(\x)
    \end{equation*} and
    \begin{equation*}
    f_2^\eps(\x)
    : = \frac{F^\eps_K \left(K(\x) \right) - F_{\V}^\eps(\x)}{K(\x)}
    \left( \M \left\langle d K,  d K\right\rangle \right)(\x).
    \end{equation*}
We also define
    \begin{equation*}
    \omega_1(h) := \frac{v_+'(h)}{h^n}, \quad
    \omega_2(h) := \frac{v_+''(h)}{2 h^{n-1}}, \quad h \in \R.
    \end{equation*} Then, for all $\x \in \overI,$
    \begin{equation*}
    \begin{split}
    \frac{1}{\eps^p} v'_+  \left( \frac{K(\x)}{\eps^p}\right)
        & = \eps^{-p(n+1)}
        \omega_1 \left( \frac{K(\x)}{\eps^p}\right)
        \left(K(\x)\right)^n\\
        & = \eps^{-1} \omega_1
        \left( \frac{K(\x)}{\eps^p}\right) \dd_n(K(\x))
    \end{split}
    \end{equation*} and
    \begin{equation*}
    \begin{split}
    \frac{1}{2\eps^{2p}} K(x) v''_+  \left( \frac{K(\x)}{\eps^p}\right)
        & = \eps^{-2p} K(\x)
        \omega_2 \left( \frac{K(\x)}{\eps^p}\right)
        \left(\frac{K(\x)}{\eps^p}\right)^{n-1}\\
        & = \eps^{-1} \omega_2
        \left( \frac{K(\x)}{\eps^p}\right) \dd_n(K(\x)).
        \end{split}
    \end{equation*}
    Therefore,
    \begin{equation*}
    \left(\M G_1^\eps\right)(\x)
        =\eps^{-1} \omega_1
        \left( \frac{K(\x)}{\eps^p}\right) \dd_n(K(\x))
        f_1^\eps(\x)
    \end{equation*}
and
    \begin{equation*}
    \left(\M G_2^\eps\right)(\x)
    = \eps^{-1} \omega_2
        \left( \frac{K(\x)}{\eps^p}\right) \dd_n(K(\x))
        f_2^\eps(\x).
    \end{equation*}
Notice that
    \begin{multline} \label{eqn:MGi-1-2}
    \left( \M G_i^\eps\right)(\x) =
    \eps^{-1}
     \biggl(f_i^\eps \left(\x \right)- \left(\A_0 f_i^\eps
        \right)([\x]) \biggr)
        \dd_n \left(K(\x)\right)
        \omega_i \left(\frac{K \left( \x \right)}{\eps^p}
        \right)\\
      + \eps^{-1}
        \left(\A_0 f_i^\eps \right)([\x])
        \dd_n \left(K(\x)\right)
        \omega_i \left(\frac{K \left( \x \right)}{\eps^p}
        \right), \quad i =1,2.
    \end{multline} Define
    \begin{equation*}
    \begin{split}
    \gamma_i(\eps):= \sup_{\eps^p<h<2 \eps^{p}}
    \left| \biggl(\A_0 f_i^\eps \biggr)
    \left( \left[K^{-1}(h) \right] \right)\right|
    = \sup_{\eps^p<h<2 \eps^{p}}
    \left| \biggl(\A_K f_i^\eps \biggr)
    \left( h\right)\right| , \quad i=1,2.
    \end{split}
    \end{equation*}
Notice that there exists a positive constant $C^*$ such that
    \begin{equation} \label{ineq:E-A0}
    \begin{split}
    \left|\E^\eps
    \left[ \int_{s \wedge \tau}^{t \wedge \tau}  \left(\A_0 f^\eps_i \right)\left([\X_u] \right)
     \dd_n\left(K \left( \X_u
        \right)\right)
        \omega_i \left(\frac{K \left( \X_u\right)}{\eps^p} \right) d u
    \right] \right| < C^* \eps \gamma_i(\eps).
    \end{split}
    \end{equation}
Now combine (\ref{eqn:MGi-1-2}), (\ref{ineq:E-A0}) and Proposition \ref{prop:averaging-partV} with $\omega_i$
for $\omega$ and with $ f_i^\eps, i = 1,2,$ for $f$ to see that
    \begin{multline} \label{ineq:E-MGi-1-2}
    \left|\E^\eps
    \left[ \int_{s \wedge \tau}^{t \wedge \tau}
        \left(\M G_i^\eps \right)
        \left(\X_s \right) ds \prod_{k=1}^m  h_k \left(\X_{t_k} \right)
        \right]\right|\\
     < C(1+t) \eps^p \left\|f^\eps_i
        \right\|_{C^5\left(\EE\right)}
    + C^* \gamma_i(\eps).
    \end{multline}

    \begin{lemma} \label{lemma:gamma-1}
    \begin{equation*}
    \sup_{\eps >0} \left(\eps^{6q}
     \left\| f_1^\eps \right\|_{C^5 \left(\EE \right)} \right)<
     \infty \quad
    \text{ and } \quad
    \lim_{\eps \rightarrow 0} \gamma_1(\eps) = 0.
    \end{equation*}
    \end{lemma}
    \begin{proof}[\bf{Proof}]
Recall that
    \begin{equation*}
    f_1^\eps(\x): =
        \left(
        \M \left\langle d F_{\V}^\eps, d K \right\rangle
        \right) \left(\x\right)
        - \left( F_K^\eps\right)'(K(\x))
        \left( \M \left\langle d K,  d K\right\rangle \right)(\x)
    \end{equation*} and
    \begin{equation*}
    \begin{split}
         \biggl(A_K f_1^\eps\biggr)
            \left( h\right)= \biggl(A_K\left(
        \M \left\langle d F_{\V}^\eps, d K \right\rangle
        \right) \biggr)\left( h\right)
        - \left( F_K^\eps\right)'(h)
        \biggl(A_K \left( \M \left\langle d K,
        d K\right\rangle \right)\biggr)
        \left(h \right).
        \end{split}
    \end{equation*}
Notice that all derivatives of $f_1^\eps$ of order $ \leq 5$ involve derivatives of $F_\V^\eps$ and
$F_K^\eps$ of order $ \leq 6$. Therefore, by Lemma \ref{lemma:extension-FV} and Lemma
\ref{lemma:extension-FK}, there exists a constant $C>0$ such that
    \begin{equation*}
    \left\|f_1^\eps \right\|_{C^5 \left(\EE \right)}
    < C \eps^{- 6q}.
    \end{equation*}

Using the gluing requirement $\underline{G}F = \overline{G}F$ on the function $F \in \D^*,$ where
$$\underline{G}F = \biggl(A^{\tt Tr}\left( \M \left\langle d F_{\V}, d K \right\rangle \right)
\biggr)\left(0\right)$$ and $$ \overline{G}F = F_K'(0)\biggl(A_K\left( \M \left\langle d K, d K \right\rangle
\right) \biggr)\left(0\right),$$ we see that
    \begin{equation*}
    \begin{split}
        &  \left|\biggl(A_K f_1^\eps\biggr)
            \left( h\right)\right|
        \leq
        \left|\biggl(A^{\tt Tr}\left(
        \M \left\langle d F_{\V}^\eps, d K \right\rangle
        \right) \biggr)\left( h\right)
        -  \biggl(A^{\tt Tr}\left(
        \M \left\langle d F_{\V}^\eps, d K \right\rangle
        \right) \biggr)\left(0 \right)\right|\\
    & \quad + \left|\biggl(A^{\tt Tr}\left(
        \M \left\langle d F_{\V}^\eps, d K \right\rangle
        \right) \biggr)\left(0 \right)
        - \biggl(A^{\tt Tr}\left(
        \M \left\langle d F_{\V}, d K \right\rangle
        \right) \biggr)\left(0\right) \right|\\
    & \quad + \left|F_K'(0)\biggl(A_K\left(
        \M \left\langle d K, d K \right\rangle
        \right) \biggr)\left(0\right) -  \left( F_K^\eps\right)'(0)
        \biggl(A_K \left( \M \left\langle d K,
        d K\right\rangle \right)\biggr)
        \left(0\right)\right|\\
    & \quad + \left|\left( F_K^\eps\right)'(0)
        \biggl(A_K \left( \M \left\langle d K,
        d K\right\rangle \right)\biggr)
        \left(0 \right)
        - \left( F_K^\eps\right)'(h)
        \biggl(A_K \left( \M \left\langle d K,
        d K\right\rangle \right)\biggr)
        \left(h\right)\right|.
        \end{split}
        \end{equation*}
for all $h \in (-a,a)$. Using properties of the trace operators, Lemma \ref{lemma:h1-h2-trace}, and the
approximation results of Lemma \ref{lemma:extension-FV} and Lemma \ref{lemma:extension-FK}, we can show that
there exists a constant $C>0$ such that
\begin{equation*}
\left| \gamma_1(\eps)\right| \leq C \left( \left\|F^\eps_\V \right\|_{\HH^2(\V)} \sqrt{|h|}+
\left\|F_{\V}^\eps - F_{\V} \right\|_{\HH^2(\V)}+ \left\|F_{K}^\eps - F_{K}
\right\|_{C^1\left(\left(0,K^*\right]\right)}+ \eps^{p-2q}\right)
\end{equation*} for all $h \in (-a,a).$
With $q<p/2,$ the inequality above implies that $ \lim_{\eps\rightarrow 0} \gamma_1(\eps) =0.$
\end{proof}

Before finishing the proof of Proposition~\ref{prop:Geps-negl}, we will need several auxiliary results.
\begin{lemma} \label{lemma:f-RfU}
Let $f \in C^n \left(\R\right)$ for some integer $n \geq 2.$

There exists a function $R_{f,\U} \in C^{n-2}\left(\R \right)$ such that
\begin{equation} \label{eqn:Taylor}
f(h) = f(0) + h f'(h) + h^2 R_{f,\U}(h),
\end{equation}
For each integer $m, \ 0 < m \leq n,$ there exists a constant $C_m>0$ such that
\begin{equation} \label{ineq:RfU}
\left\| R_{f,\U}\right\|_{C^{m-2} \left( [-a,a] \right)} \leq C_m \left\|f \right\|_{C^{m} \left( [-a,a]
\right)}.
\end{equation}
\end{lemma}
A very simple proof of Lemma~\ref{lemma:f-RfU} can be found in \cite{S.1} (Lemma 9.4).

\begin{lemma} \label{lemma:f-RfV}
Let $f \in C^n \left(\R^2\right)$ for some integer $n \geq 2,$ and suppose that $f(\x) \equiv f_0 $ for $\x
\in \partV.$ Then
    \begin{equation*}
    \left(\M \langle df, dK \rangle \right)(\x)
    = \frac{ \left( \nabla f, \nabla K \right) (\x)}
    {\left\|\nabla K (\x)\right\|^2 }
    \left( \M \langle dK, dK\rangle\right) (\x)
    \end{equation*} for all $\x \in \partV.$
There exists a function $R_{f,\V} \in C^{n-2}\left(\R^2 \right)$ such that for every $\x \in \EE$,
    \begin{equation} \label{eqn:f-RV-expand}
    f(\x) = f_0 + K(\x)
        \frac{\left( \M \langle df, dK \rangle\right)(\x)}
        {\left(\M \langle dK, dK \rangle \right)(\x)}
        + K^2(\x) R_{f,\V}(\x).
\end{equation}
For each integer $m, 0 < m \leq n,$ there exists a constant $C_m>0$ such that
\begin{equation} \label{ineq:RfV}
\left\| R_{f,\V}\right\|_{C^{m-2} \left( \overline{\EE} \right)} \leq C_m \left\|f \right\|_{C^{m} \left(
\overline{\EE} \right)}.
\end{equation}
\end{lemma}

\begin{proof}[\bf{Proof}]
Let $\x \in\partV,$ then $\left(\nabla f, \nabla^\perp K \right)(\x) = 0$. By projecting $\nabla f$ onto two
orthogonal directions of $\nabla K$ and $\nabla^\perp K,$ we see that
    \begin{equation*}
    \begin{split}
    \nabla f (\x)
    & = \frac{\left(\nabla f, \nabla K \right)(\x)}
    {\left\|\nabla K (\x)\right\|^2} \nabla K(\x)
    + \frac{\left(\nabla f, \nabla^\perp K \right)(\x)}
    {\left\|\nabla^\perp K (\x)\right\|^2} \nabla^\perp K(\x)\\
    &  = \frac{\left(\nabla f, \nabla K \right)(\x)}
    {\left\|\nabla K (\x)\right\|^2} \nabla K(\x).
    \end{split}
    \end{equation*}
Therefore, using the definition of $\langle df,dg \rangle_t$, we see that
    \begin{equation*}
    \left( \M \langle df, dK \rangle \right)(\x)
    = \frac{ \left( \nabla f, \nabla K \right) (\x)}
    {\left\|\nabla K (\x)\right\|^2 } \
     \left( \M\langle dK, dK\rangle \right)(\x)
    \end{equation*} for all $\x \in \partV.$

Using the flow $\zeta_t$ on $\EE$ introduced in (\ref{eqn:flow-on-E}) we consider the function $F:t \mapsto f
\left(\zeta_{-t K(\x)}(\x) \right)$. Notice that
    \begin{equation*}
    \dot{\zeta}_{-tK(\x)}(\x)
    = - K(\x)\frac{\nabla K}{\|\nabla K\|^2} \left(\zeta_{-tK(\x)}(\x)
    \right),
    \end{equation*}
$F(0)=f\left(\zeta_{0}(\x)\right) = f(\x),$ and $F(1)=f\left(\zeta_{-K(\x)}(\x)\right) = f_0.$ Using the
identity
    \begin{equation*}
    F(t) = F \left(s\right)
        + F'(t) \left(t-s \right)
        + \int_s^t\left(s-r\right)F''(r)dr
    \end{equation*}
    with $t=0$ and $s=1,$ we obtain
    \begin{equation*}
    \begin{split}
    f(\x)
    & = f_0 + K(\x) \frac{ \left( \nabla f, \nabla K \right) (\x)}
    {\left\|\nabla K (\x)\right\|^2 } + K^2(\x) R_{f,\V}(\x)\\
    & = f_0 + K(\x) \frac{ \left( \M \langle df, dK \rangle \right)(\x)}
    {\left( \M\langle dK, dK\rangle \right)(\x) } + K^2(\x) R_{f,\V}(\x),
    \end{split}
    \end{equation*}for all $\x \in \EE$, where
    \begin{equation*}
    R_{f,\V}(\x)
    = \int_0^1 (t-1) \left( \nabla
        \left( \frac{ \left( \nabla f, \nabla K \right)}
    {\left\|\nabla K \right\|^2 }\right),
    \frac{\nabla K}{\|\nabla K\|^2}\right)
    \left(\zeta_{-tK(\x)}(\x)\right)dt.
    \end{equation*}
Let $0<m \leq n$, then from this expression for $R_{f,\V}(\x)$ we see that there exists a constant $C_m>0$
such that the bound (\ref{ineq:RfV}) holds.
\end{proof}
\begin{lemma} \label{lemma:f3eps}
There exists a function $f_3^\eps \in C^5 \left( \EE \right)$ satisfying $\A_0f_3^\eps \equiv 0,$ and
    \begin{equation} \label{eqn:bound-f3}
    \sup_{\eps>0}\eps^{8q} \left\| f_3^\eps
    \right\|_{C^5\left(\EE\right)}< \infty,
    \end{equation}such that for all $\x \in
\EE,$
    \begin{equation} \label{eqn:Kperp-F-V}
    \left(\nabla F_{\V}^\eps(\x), \nabla^\perp K(\x)\right)
    = K(\x) f_3^\eps(\x).
    \end{equation}
\end{lemma}

\begin{proof}[\bf{Proof}]
By Lemma \ref{lemma:f-RfV},
\begin{equation*}
    F^\eps_{\V}(\x) = F(\partV) + K(\x)
        \frac{\left( \M \langle d F^\eps_{\V}, dK \rangle\right)(\x)}
        {\left(\M \langle dK, dK \rangle \right)(\x)}
        + K^2(\x) R_{F^\eps_{\V},\V}(\x)
\end{equation*}
for all $\x \in \EE$, which implies that
\begin{multline*}
    \left(\nabla F^\eps_{\V}(\x),\nabla^\perp K(\x)\right)
     = K(\x)\left(
    \nabla \left( \frac{\left(\M \langle dF^\eps_{\V},dK \rangle\right)(\x)}
        {\left(\M \langle dK,dK
        \rangle\right)(\x)}\right), \nabla^\perp K(\x)\right)\\
         + K^2(\x) \left(  \nabla
        R_{F^\eps_{\V},\V}(\x), \nabla^\perp K(\x)\right).
\end{multline*}
 Take
    \begin{equation*}
    f_3^\eps(\x) = \left( \nabla^\perp K(\x),
    \nabla \left( \frac{\left(\M \langle dF^\eps_{\V},dK \rangle\right)(\x)}
        {\left(\M \langle dK,dK \rangle\right)(\x)}\right)\right)
        + K(\x) \left( \nabla^\perp K(\x), \nabla
        R_{F^\eps_{\V},\V}(\x)\right).
    \end{equation*}
Then (\ref{eqn:Kperp-F-V}) holds. From this identity and from (\ref{eqn:grad-f3}) we see that for any $\x \in
\EE$
\begin{equation*}
\left(\A_0 \left(K(\x) f_3^\eps \right)\right)([\x])= K(\x) \left(\A_0 f^\eps_3 \right)([\x]) \equiv 0.
\end{equation*} Therefore, $\A_0 f^\eps_3 \equiv
0.$

To verify that (\ref{eqn:bound-f3}) holds, we use inequality (\ref{ineq:RfV}) and Lemma
\ref{lemma:extension-FV}. This finishes the proof of Lemma \ref{lemma:f3eps}.
\end{proof}

\begin{lemma} \label{lemma:gamma-2}
    \begin{equation*}
    \sup_{\eps >0} \left(\eps^{7q}
     \left\| f_2^\eps \right\|_{C^5 \left(\EE \right)} \right)<
     \infty \quad
    \text{ and } \quad \lim_{\eps \rightarrow 0} \gamma_2(\eps) = 0.
    \end{equation*}
    \end{lemma}
    \begin{proof}[\bf{Proof}] Recall that
    \begin{equation*}
    f_2^\eps(\x)
    : = \frac{F^\eps_K \left(K(\x) \right) - F_{\V}^\eps(\x)}{K(\x)}
    \left( \M \left\langle d K,  d K\right\rangle \right)(\x),
    \end{equation*}
$F_K^\eps(0)= F \left( \G_{\partV} \right), $ and $ F_{\V}^\eps (\x) = F \left( \G_{\partV} \right)$ for $ \x
\in \partV.$

By Lemma \ref{lemma:f-RfU} and Lemma \ref{lemma:f-RfV}, if $\x \in \EE$ then
    \begin{equation*}
    \begin{split}
    \frac{F^\eps_K \left(K(\x) \right) - F_{\V}^\eps(\x)}{K(\x)}
     & = \frac{F^\eps_K \left(K(\x) \right)
     - F \left( \G_{\partV} \right)}{K(\x)}
     - \frac{F_{\V}^\eps(\x)
     - F \left( \G_{\partV} \right)}{K(\x)}\\
     & =  \left(F_K^\eps\right)' \left(K(\x) \right)
        + K(\x) R_{F^\eps_K,\U}\left(K(\x) \right)\\
     & \quad - \frac{\left( \M \langle dF^\eps_{\V}, dK \rangle\right)(\x)}
        {\left(\M \langle dK, dK \rangle\right)(\x)} - K(\x) R_{F^\eps_{\V},\V}(\x)
     \end{split}
    \end{equation*}
Now, by Lemma \ref{lemma:extension-FK}, Lemma \ref{lemma:extension-FV}, and applying inequalities
(\ref{ineq:RfU}) and (\ref{ineq:RfV}), there exists a constant $C>0$ such that
    \begin{equation*}
    \left\|f_2^\eps \right\|_{C^3 \left(\EE \right)}
    < C \eps^{- 7q}.
    \end{equation*}

Similarly, notice that
\begin{multline*}
    f_2^\eps(\x) - f_1^\eps(\x)
    = \left( \frac{F^\eps_K \left(K(\x) \right) - F_{\V}^\eps(\x)}{K(\x)}
    - \left( F_K^\eps\right)'(K(\x))\right)
    \left( \M \left\langle d K,  d K\right\rangle \right)(\x)\\
    - \left(
        \M \left\langle d F_{\V}^\eps, d K \right\rangle
        \right) \left(\x\right),
        \end{multline*} and so there exists a constant $C>0$ such that
    \begin{equation*}
    \left| f_2^\eps (\x) - f_1^\eps(\x) \right| \leq C \eps^{p-2q}
    \end{equation*}
for all $\x \in \overI \cap \left\{\y: \eps^p < K(\y) < 2 \eps^p \right\}.$ Now,
    \begin{equation*}
    \gamma_2(\eps) \leq \gamma_1(\eps) + C \eps^{p-2q}.
    \end{equation*}
Therefore, by Lemma \ref{lemma:gamma-1}, and since $q<p/2$, $ \lim_{\eps \rightarrow 0}\gamma_2(\eps) = 0.$
    \end{proof}

From Lemma \ref{lemma:gamma-1} and Lemma \ref{lemma:gamma-2}, $\lim_{\eps \rightarrow 0} \gamma_1(\eps) =
\lim_{\eps \rightarrow 0} \gamma_2(\eps)  = 0.$ Therefore, in the case of $i=1,2$, (\ref{eqn: MG's}) follows
from (\ref{ineq:E-MGi-1-2}). The case $i=0$ is trivial due to the fact that $\left(\M G_0^\eps\right)(\x)
\equiv 0.$ This completes the proof of Proposition~\ref{prop:Geps-negl}.
\end{proof}

 \begin{proposition} \label{prop:Geps3-negl}
Fix $0 \leq t_1<t_2< \cdots< t_m \leq s < t$ and $h_1, \cdots , h_m \in C_b\left( \R^2 \right).$ If $F \in
\D^*$ then
    \begin{equation*}
    \lim_{\eps \rightarrow 0} \E^\eps
    \left[ \int_{s \wedge \tau}^{t \wedge \tau}
        \left(\M G_3^\eps \right)
        \left(\X_s \right) ds \prod_{k=1}^m h_k \left(\X_{t_k} \right)
            \right]
    = 0.
    \end{equation*}
    \end{proposition}
    \begin{proof}[\bf{Proof}]
Notice that
\begin{equation*}
\begin{split}
\left( \M G^\eps_3\right)(\x) & = \frac{1}{\eps}
        \left( 1 - v_+  \left(\frac{K(\x)}{\eps^p}\right)\right)
        \left(\nabla F_{\V}^\eps(\x),\nabla^\perp H(\x) \right)\\
        & = \frac{1}{\eps}
        \left( 1 - v_+  \left(\frac{K(\x)}{\eps^p}\right)\right)
        n \dd_{n-1}(K(\x))
        \left(\nabla F_{\V}^\eps(\x),\nabla^\perp K(\x) \right)
        \end{split}
\end{equation*} and
    \begin{align}
    \left(\A_0 \left(\nabla F_{\V}^\eps, \nabla^\perp K\right)
        \right)([\x])
    & = \frac{1}{\eta \left(K(\x)\right)} \int_{0}^{\eta \left(K(\x)\right)}
    \left(\nabla F_{\V}^\eps, \nabla^\perp K\right)(\phi_t(\x)) d t \notag\\
    &  = \frac{1}{\eta \left(K(\x)\right)} \int_{0}^{\eta \left(K(\x)\right)}
    \frac{d}{dt} F_{\V}^\eps (\phi_t(\x)) dt \notag\\
    & = F_{\V}^\eps \left(\phi_{\eta \left(K(\x)\right)}\left(\x\right)\right)
        - F_{\V}^\eps (\x)=0. \label{eqn:grad-f3}
    \end{align}

Let $f_3^\eps$ be a function from Lemma \ref{lemma:f3eps}, and define
\begin{equation*}
\omega_3(h):=n \left(1 - v_+(h)\right), \quad h \in \R.
\end{equation*} Then,
\begin{equation} \label{eqn:MGi-3}
    \left(\M G_3^\eps\right)(\x)
        =\eps^{-1} \omega_3
        \left( \frac{K(\x)}{\eps^p}\right) \dd_n(K(\x))
        f_3^\eps(\x).
    \end{equation}
Now combine (\ref{eqn:MGi-3}), and Proposition \ref{prop:averaging-partV} with $\omega_3$ for $\omega$ and
with $ f_3^\eps$ for $f$ to see that
    \begin{equation} \label{ineq:E-MGi-3}
    \left|\E^\eps
    \left[ \int_{s \wedge \tau}^{t \wedge \tau}
        \left(\M G_3^\eps \right)
        \left(\X_s \right) ds  \prod_{k=1}^m h_k \left(\X_{t_k} \right)
        \right]\right| < C\, \eps^p \,(1+t)  \left\|f^\eps_3
    \right\|_{C^5\left(\EE\right)}.
    \end{equation} Thus, Proposition \ref{prop:Geps3-negl} holds.
    \end{proof}

    \begin{proposition} \label{prop:Keps-negl}
    Fix $F \in \D^*$ and $t >0.$ Then
    \begin{equation*}
    \lim_{\eps \rightarrow 0} \E^\eps
    \left[ \left|\int_0^{t \wedge \tau}
        \left(\M R^\eps \right)
        \left(\X_s \right) ds \right|\right]
    = 0.
    \end{equation*}
    \end{proposition}

    \begin{proof}[\bf{Proof}]
Consider
\begin{equation*}
\left( \M R^\eps\right) (\x) = \frac{1}{ \eps^{p}}
        v'_+  \left(\frac{K(\x)}{\eps^p}\right)
        \biggl( F^\eps_K \left(K(\x) \right) - F_{\V}^\eps(\x)\biggr)
        \left( \M \left(\LL K\right)\right)(\x)
\end{equation*}
Define \begin{equation*}
    \beta_K(\eps): = \sup_{\substack{\x \in \overI\\\eps^p \leq
    K(\x) \leq 2 \eps^{p}}} \left|F^\eps_K \left(K(\x) \right) - F_{\V}^\eps(\x)
     \right|,
\end{equation*}
which satisfies
    \begin{equation} \label{ineq:betaK}
    \beta_K(\eps) \leq
    \sup_{\substack{\x \in \overI\\
        \eps^p \leq K(\x) \leq 2\eps^{p}}}
    \biggl(
        \left|F^\eps_K \left(K(\x) \right) - F \left(\G_{\partV}\right) \right|
        + \left| F_{\V}^\eps(\x)- F \left(\G_{\partV}\right) \right|
    \biggr).
\end{equation}
Notice that
    \begin{equation*}
    F^\eps_K(K(\x)) - F\left(\G_{\partV}\right)
    = F^\eps_K(K(\x)) - F^\eps_K(0)
    = \left(F^\eps_K \right)'(0) K(\x) + o\left( |K(\x)|\right).
    \end{equation*}
From this expression and Lemma \ref{lemma:extension-FK}, there exist constants $C_1,C_{11}>0$ such that
    \begin{equation} \label{ineq:FK-F}
 \left|F^\eps_K(K(\x)) - F \left(\G_{\partV}\right) \right|< C_1 \eps^p \left\|F^\eps_{K}
 \right\|_{C^1 \left(\left[0,K^*\right] \right)}
        < C_{11} \eps^{p-q}.
    \end{equation}
By the definition of $K$, if $\x^* \in \partV,$ then $K(\x^*) = 0$, and $\nabla K(\x^*) \neq 0.$ Therefore,
there exists $\x^* \in \partV,$ such that
    \begin{equation*}
    \begin{split}
    F^\eps_{\V}(\x) - F\left(\G_{\partV}\right)
    & = F^\eps_{\V}(\x) - F^\eps_{\V}\left(\x^*\right)\\
    & = \left(\nabla F^\eps_{\V}\left(\x^*\right), \x-\x^*\right)
    + o \left( \left\|\x - \x^*\right\|\right)\\
    & = \left\|\nabla F^\eps_{\V}\left(\x^*\right)\right\|
    \left\|\x-\x^* \right\|
    + o \left( \left\|\x - \x^*\right\|\right)\\
    & = \frac{\left\|\nabla F^\eps_{\V}\left(\x^*\right)\right\|}
    {\left\| \nabla K\left(\x^*\right) \right\|}K(\x)
    + o\left( |K(\x)|\right).
    \end{split}
    \end{equation*}
From this expression and Lemma \ref{lemma:extension-FV}, there exist constants $C_2,C_{22}>0$ such that
    \begin{equation} \label{ineq:FV-F}
 \left|F^\eps_{\V}(\x) - F\left(\G_{\partV}\right) \right|< C_2 \eps^p \left\|F^\eps_{\V}
 \right\|_{C^1 \left(\overI \right)}
        < C_{22} \eps^{p-q}.
    \end{equation}
Now, combine inequalities (\ref{ineq:betaK}), (\ref{ineq:FK-F}), and (\ref{ineq:FV-F}) to conclude that
$\beta_K(\eps) \leq \left(C_{11}+C_{22}\right) \eps^{p - q}$ for all $\eps>0.$ Therefore, since
$$\left|v'_+  \left(\frac{K(\x)}{\eps^p} \right) \right|\leq C_{w}\chi_{[\eps^p,2 \eps^p]} \left( K(\x) \right)$$
for some $C_{w}>0$ and $$\sup_{\x \in \overI} \left| \left(\M \left(\LL K \right)\right)(\x)\right| \leq
C_{K}$$ for some $C_{K}>0,$ there exists a positive constant $C := C_{w} C_{K} \left(C_{11}+C_{22}\right)$
such that
    \begin{equation*}
    \E^\eps \left[\left|
        \int_0^{t \wedge \tau}
        \left(\M R^\eps\right) \left(\X_u \right)d u
            \right| \right]
             \leq C \eps^{p - q}
             \E^\eps \left[ \frac{1}{\eps^p}
    \int_0^{t \wedge \tau} \chi_{[\eps^p,2 \eps^p]} \left( K(\X_u) \right) d u
    \right].
    \end{equation*}
This inequality and Lemma \ref{lemma:residence} give the desired result.
    \end{proof}

    \begin{proposition}
    Fix $0 \leq t_1<t_2< \cdots<t_m \leq s < t$ and $h_1, \cdots , h_m \in C_b\left( \R^2 \right).$ If $F \in
\D^*,$ then
    \begin{equation} \label{eqn:Eeps(Leps-L)}
    \lim_{\eps \rightarrow 0} \E^\eps
    \left[ \int_0^{t \wedge \tau}
        \biggl( \left(\M \left(\LL^\eps F^\eps \right)\right)
        \left(\X_s \right) - \left(\LL^* F\right)
        \left(\left[\X_s \right]\right)\biggr) ds \prod_{k=1}^m h_k \left(\X_{t_k} \right) \right]
    = 0.
    \end{equation}
    \end{proposition}

    \begin{proof}[\bf{Proof}]
Given $\x \in \overI,$ consider
    \begin{multline*}
    \left(\M \left(\LL^\eps F^\eps \right)\right)
        \left(\x \right) - \left(\LL^* F\right)
        \left([\x] \right)
     = \biggl( \left(\M L^\eps \right)
        \left(\x \right) - \left(\LL^* F\right)
        \left([\x] \right)\biggr)\\
     +  \left( \M K^\eps\right) \left( \x \right)
        + \sum_{k=0}^3 \left( \M G^\eps_i\right) \left( \x\right)
        + \left( \M R^\eps \right) \left( \x\right).
        \end{multline*}
Now, combine Proposition \ref{prop:Leps-L*f}, Proposition \ref{prop:Geps-negl}, Proposition
\ref{prop:Geps3-negl}, and Proposition \ref{prop:Keps-negl} to verify that (\ref{eqn:Eeps(Leps-L)}) holds.
    \end{proof}

Finally, we can prove the desired statement about the limiting martingale problem.

\begin{theorem} \label{thm:Pstar}
Fix $\x \in \overline{\I}$. Let $\PP^* \in \mathcal{P} \left( C \left( [0,\infty); \G \right) \right)$ be a
cluster point of the $\PP^{\eps,*}$s. Then $\PP^*$ satisfies the martingale problem for
$\left(\LL^*,\delta_{[\x]}\right)$.
\end{theorem}

\begin{proof}[\bf{Proof}]
Let $\E^{\eps}$ be the expectation operator associated with the law $\PP^{\eps}$ and let $\E^*$ be the
expectation operator associated with the law $\PP^*. $

We know that every $\PP^\eps \in \mathcal{P}(C[0,\infty); \R^2 )$ is a solution of the stopped martingale
problem for $ \left(\LL_{t/\eps^2}^\eps,\delta_{\x}, \I \right),$ and that $F^\eps \in C^2\left(\R^2 \right)
\subset \D^\eps.$ Therefore, we know that for any $0 \leq s < t, $ and $\eps>0, $
    \begin{equation} \label{eqn:mpX}
    \begin{split}
    F^\eps\left( \X_{t \wedge \tau }\right)
            - F^\eps \left( \X_{s  \wedge \tau} \right)
        - \int_{s \wedge \tau}^{t \wedge \tau}
            \left( \LL_{u/\eps^2}^\eps F^\eps \right)
                \left( \X_u\right) d u
     \end{split}
\end{equation} is a $\PP^\eps$-martingale. We also know that
$\PP^\eps \left\{ \X_0 = \x  \right\}=1.$

Our theorem claims that $\PP^* \left\{ \Y_0 = [\x] \right\} = 1$ and if $F \in \D^*,$ $0 \leq t_1<t_2<
\dots<t_m \leq s < t,$ and $h^*_1, \ldots,h^*_m \in C_b(\G)$ then
    \begin{equation} \label{eqn:mpY}
    \E^*
        \Biggl[
        \Biggl\{
            F\left(\Y_{t}\right)
            - F \left(\Y_{s} \right)
        - \int_{s }^{t}
            \left( \LL^* F \right)
                \left( \Y_u \right) d u
        \Biggr\}
            \prod_{k=1}^m h^*_k \left( \Y_{t_k} \right)
        \Biggr]=0.
    \end{equation}
Notice that if $B \subset \G$ is any set which does not contain $[\x]$, then $\pi^{-1}(B) \subset \overI$
does not contain $\x$, and
\begin{equation*}
\PP^* \left\{ \Y_0 \in B \right\} \leq \overline{\lim_{\eps \rightarrow 0}} \PP^\eps \left\{ \Y_0 \in B
\right\} = \overline{\lim_{\eps \rightarrow 0}} \PP^\eps \left\{ \X_0 \in \pi^{-1}(B)\right\} = 0.
\end{equation*}
Therefore $ \PP^* \left\{ \Y_0 = [ \x ] \right\} = 1.$

 By the definition of $\D^*$ we know that $\left(\LL^*
F \right) \left( \G_{\partI}\right)=0,$ which implies that for any $0 \leq s < t,$
\begin{equation*}
\int_{s \wedge \tau}^{t \wedge \tau} \left( \LL^* F \right)(\Y_s) ds = \int_s^t \left( \LL^* F \right)(\Y_s)
ds.
\end{equation*}

Now,
    \begin{align}
    \E^*
        & \Biggl[ \Biggl\{
            F\left(\Y_{t}\right) - F \left(\Y_{s} \right)
                - \int_{s }^{t} \left( \LL^* F \right) \left( \Y_u \right) d u
            \Biggr\} \prod_{k=1}^m h^*_k \left( \Y_{t_k} \right) \Biggr] \notag\\
        &  = \lim_{\eps \rightarrow 0} \E^{\eps}
            \Biggl[ \Biggl\{
            F\left(\Y_{t}\right) - F \left(\Y_{s} \right)
            \Biggr.\Biggr. \Biggl.\Biggl.
                - \int_{s }^{t} \left( \LL^* F \right) \left( \Y_u \right) d u
            \Biggr\} \prod_{k=1}^m h^*_k \left( \Y_{t_k} \right) \Biggr] \notag\\
        &  = \lim_{\eps \rightarrow 0}  \E^\eps
            \Biggl[ \Biggl\{
            F \circ \pi \left(\X_{t \wedge \tau}\right)
                - F \circ \pi \left(\X_{s \wedge \tau}\right)
            \Biggr.\Biggr. \Biggl.\Biggl.
                - \int_{s \wedge \tau}^{t \wedge \tau}
                \left( \LL^* F \right) \left( \Y_u \right) d u
            \Biggr\} \prod_{k=1}^m h^*_k \left( \Y_{t_k} \right) \Biggr]. \label{eqn:lim-mtg-pr}
    \end{align}

Notice that
    \begin{equation*}
    \begin{split}
    F \circ \pi
            \left(\X_{t \wedge \tau}\right)
                & - F \circ \pi \left(\X_{s \wedge \tau}\right)
                - \int_{s \wedge \tau}^{t \wedge \tau}
                \left( \LL^* F \right) \left( \Y_u \right) d u\\
        &  = \biggl(F \circ \pi \left(\X_{t \wedge \tau}\right)
                - F^\eps\left(\X_{t \wedge \tau}\right)\biggr)
                - \biggl(F \circ \pi\left(\X_{s \wedge \tau}\right)
                - F^\eps\left(\X_{s \wedge \tau}\right)\biggr)\\
        & + F^\eps \left(\X_{t \wedge \tau} \right)
                - F^\eps \left(\X_{s \wedge \tau} \right)
                - \int_{s \wedge \tau}^{t \wedge \tau}
                \LL^\eps_{u/\eps^2} F^\eps \left(\X_u \right) d u\\
        & +  \int_{s \wedge \tau}^{t \wedge \tau}
                \biggl( \left(\LL^\eps_{u/\eps^2} F^\eps\right) \left(\X_u \right)
                -\left(\M \left(\LL^\eps F^\eps \right)\right) \left(\X_u \right)
                \biggr)du\\
        & +  \int_{s \wedge \tau}^{t \wedge \tau}
                \biggl( \left(\M \left(\LL^\eps F^\eps \right)\right) \left(\X_u \right)
                - \left(\LL^* F\right) \left(\Y_u \right)\biggr) du
     \end{split}
    \end{equation*}
Keeping this expression in mind, we use the statement of the martingale problem (\ref{eqn:mpX}), Lemma \ref{lemma:f-Mf}
with $f(\x,t) = \LL^\eps_t F^\eps(\x)$ concerning the first-level averaging, and Eq.~(\ref{eqn:Eeps(Leps-L)}) to see
that (\ref{eqn:lim-mtg-pr}) implies (\ref{eqn:mpY}), and this completes the proof.
\end{proof}

\subsection{Uniqueness}
\noindent
 We already know that every limit point of the $\PP^\eps$-laws of $\left\{ \Y_t, t \geq 0 \right\}$
satisfies the martingale problem for $\left(\LL^*,\delta_{[\x]}\right)$. Our next step is to verify that the
$\PP^\eps$-laws have only one limit point.

\begin{proposition}
The operator $\LL^*$ generates a strongly continuous contraction semigroup on $C \left(\G\right).$
\end{proposition}

\begin{proof}[\bf{Proof}]
To prove this proposition using the Hille-Yosida theorem, we need to show that $\D^*$ is dense in $C(\G),$
that $\LL^*$ satisfies the positive maximum principle, and that the range of $\lambda I -\LL^*$ is dense in
$C(\G)$ for some $\lambda>0.$

\begin{lemma} \label{lemma:density}
$\D^*$ is dense in $C(\G).$
\end{lemma}

\begin{proof}[\bf{Proof}]
Let $F \in C(\G).$ To prove the statement of this lemma, we will approximate this function with functions
$F_n \in \D^*$ in the $\| \cdot \|_{C(\G)}$-norm. We will construct the $F_n$'s in two steps.

Given $[\x] \in \G$, take
\begin{equation*}
\begin{split}
F^0_n([\x]) & := F([\x]) + \biggl( F\left(\G_{\partV}\right) - F([\x])\biggr)v \left(2nK(\x)\right) \\
& + \biggl( F\left(\G_{\partI}\right) - F([\x])\biggr)v \left(2n\left(K(\x)-K^* \right)\right).
\end{split}
\end{equation*} Then,
\begin{equation} \label{eqn:F-Fn0}
\lim_{n \rightarrow 0} \left\| F- F^0_n \right\|_{C(\G)} = 0.
\end{equation}

Second, for every $n$, using the standard mollification, there exists a function   $F^1_n~\in~C^\infty(\V~\cup~\U)$
such that
    \begin{equation*}
    \sup_{\substack{|K(\x)| \geq \frac{1}{2n} \\
        \left|K(\x) - K^*\right| \geq \frac{1}{2n}}}
        \left| F([\x]) - F^1_n([\x])\right| \leq \frac{1}{n}.
    \end{equation*}
Now choose
\begin{equation*}
\begin{split}
F_n([\x]) & := F^1_n([\x]) + \biggl( F\left(\G_{\partV}\right) - F^1_n([\x])\biggr)v \left(2nK(\x)\right) \\
& + \biggl( F\left(\G_{\partU}\right) - F^1_n([\x])\biggr)v\left(2n\left(K(\x)-K^* \right)\right).
\end{split}
\end{equation*}
Then,
\begin{equation}\label{eqn:Fn-Fn0}
\lim_{n \rightarrow 0} \left\| F^0_n - F_n \right\|_{C(\G)} = 0.
\end{equation}
Each of the $F_n$s is a smooth function on $ \V \cup \U$ and is a constant in a small neighborhood of
$\partV$ and in a small neighborhood of $\partI$. Thus, $F_n \in \D^*.$ Finally, (\ref{eqn:F-Fn0}) and
(\ref{eqn:Fn-Fn0}) imply that
\begin{equation*}
\lim_{n \rightarrow 0} \left\| F- F_n \right\|_{C(\G)} = 0.
\end{equation*} Therefore, $\D^*$ is dense in $C(\G)$.
\end{proof}

\begin{lemma} \label{lemma:PMP}
$\LL^*$ satisfies the positive maximum principle.
\end{lemma}
\begin{proof}[\bf{Proof}]
By the definition of the positive maximum principle we need to show that whenever $F \in \D^*,$ $\left[\x_0 \right]\in
\G,$ and $ \sup_{\left[\x \right]\in \G}F(\left[\x\right]) = F\left(\left[\x_0\right]\right) \geq 0$, we have $\LL^*F
\left(\left[\x_0 \right]\right) \leq 0.$

Consider $F \in \D^*,$ and let it attain its supremum at a point $\left[\x_0\right] \in \G.$ If $\left[\x_0\right] =
\G_{\partI}$ then, by definition of $\D^*$, we have $\LL^* F\left(\left[\x_0\right]\right) = 0.$

Assume that $ \left[\x_0\right] \in \G_\V \cup \G_\U.$ Recall that the operator $\LL^*$ acts locally as a strongly
elliptic operator on both $\V$ and $\U$. Therefore, $\LL^* F\left(\left[\x_0\right]\right) \leq 0.$

Assume that $F$ attains its supremum at $ \left[\x_0\right] =\G_{\partV}$, then $\overline{G}F \geq 0$ and
$\underline{G}F \leq 0.$ Since $F \in \D^*$, we know that $\overline{G}F = \underline{G}F = 0.$ This implies that
    \begin{equation*}
    F_K'(0) = \frac{\overline{G}F}{\sg^2(0)} = 0.
    \end{equation*}
By Lemma \ref{lemma:extension-FK}, $F_K \in C^2 \left(\left[0,K^*\right]\right)$ and $h =0 $ is its local
maximum. Therefore, $F_K''(0) \leq 0.$ We notice that in this case
    \begin{equation*}
    \begin{split}
    \left(\LL^*F \right)\left(\left[\x_0\right]\right)
     = \lim_{h \downarrow 0}
        \left(\M \left( \LL^K F_K\right)\right)(h)
    = \bb(0)F_K'(0)
    + \frac{1}{2}\sg^2 \left(0\right) F_K''(0)
    \leq 0.
    \end{split}
    \end{equation*} This finishes the proof of Lemma~\ref{lemma:PMP}.
\end{proof}

\begin{lemma} \label{lemma:range}
The range of $\lambda I -\LL^*$ is dense in $C\left(\G\right)$ for some $\lambda>0$.
\end{lemma}

We choose to omit the proof of this lemma since it is very similar to that of Lemma 8.4 in \cite{S.1}.

Lemmas \ref{lemma:density}, \ref{lemma:PMP}, and \ref{lemma:range} above allow us to conclude that uniqueness
holds for the solution of the martingale problem for $\left(\LL^*,\delta_{[\x]}\right)$. Therefore, the
$\PP^\eps$-laws of $\left\{\Y_t, t \geq 0\right\}$ converge to the law of a $\G$-valued Markov process.
\end{proof}


\nonumsection{References} \noindent
\bibliographystyle{plain}
\bibliography{bibliography}

\end{document}